\documentclass[A4j,11pt]{article}
\usepackage{amssymb,amsmath,amsthm,amscd,amsfonts}
\usepackage{amsmath,amssymb,amsthm,amscd,ascmac,mathrsfs,bm,type1cm,multirow,array,enumerate,mathrsfs,tabularx}
\usepackage{graphicx}
\usepackage{graphics}

\begin{document}

\title{{\bf Regularized mean curvature flow for\\
invariant hypersurfaces in a Hilbert space\\
and its application to gauge theory}}
\author{{\bf Naoyuki Koike}}
\date{}
\maketitle

\begin{abstract}
In this paper, we investigate a regularized mean curvature flow starting from an invariant 
hypersurface in a Hilbert space equipped with an isometric and almost free action of a Hilbert Lie 
group whose orbits are minimal regularizable submanifolds.  
We prove that, if the initial invariant hypersurface satisfies a certain kind of horizontally convexity 
condition and some additional conditions, then it collapses to an orbit of the Hilbert Lie group action along 
the regularized mean curvature flow.  
In the final section, we state a vision for applying the study of the regularized mean curvature flow to the gauge theory.  
\end{abstract}

\section{Introduction}
C. L. Terng (\cite{Te1}) defined the notion of a proper Fredholm submanifold in a 
(separable) Hilbert space as a submanifold of finite codimension satisfying certain conditions for the normal 
exponential map.  Note that the shape operators of a proper Fredholm submanifold are compact operator.  
By using this fact, C. King and C. L. Terng (\cite{KiTe}) defined 
the regularized trace of the shape operator for each unit normal vector of a proper Fredholm submanifold.  
Later, E. Heintze, C. Olmos and X. Liu (\cite{HLO}) defined another 
regularized trace of the shape operator for each unit normal vector
of a proper Fredholm submanifold, which differs from one defined in \cite{KiTe}.  
They called the regularized trace defined in \cite{KiTe} $\zeta$-{\it reguralized trace}.  
The regularized trace in \cite{HLO} is easier to handle than one in \cite{KiTe}.  In almost all relevant cases, 
these regularized traces coincide.  In this paper, we adopt the regularized trace defined in \cite{HLO}.  
Let $M$ be a proper Fredholm submanifold in $V$ immersed by $f$.  
If, for each normal vector $\xi$ of $f$, the regularized trace ${\rm Tr}_r\,A_{\xi}$ of the shape operator 
$A_{\xi}$ of $f$ and the trace ${\rm Tr}\,A_{\xi}^2$ of $A_{\xi}^2$ exist, then $M$ (or $f$) is said to be 
{\it regularizable}.  See Section 2 about the definition of the regularized trace ${\rm Tr}_r\,A_{\xi}$.  
Let $M$ be a Hilbert manifold and $\{f_t\}_{t\in[0,T)}$ be a $C^{\infty}$-family of regularizable immersions of 
codimension one of $M$ into $V$ which admit a unit normal vector field $\xi_t$.  
The regularized mean curvature vector $H_t$ is defined by $H_t:=-{\rm Tr}_r((A_t)_{-\xi_t})\cdot\xi_t$, where 
$A_t$ denotes the shape tensor of $f_t$.  
Define a map $F:M\times[0,T)\to V$ by $F(x,t):=f_t(x)$ ($(x,t)\in M\times[0,T)$).  
We call $\{f_t\}_{t\in[0,T)}$ the {\it regularized mean curvature flow} 
if the following evolution equation holds:
$$\frac{\partial F}{\partial t}=H_t.\leqno{(1.1)}$$
This notion was introduced in \cite{Koi2}.  

R. S. Hamilton (\cite{Ha}) proved the existence and the uniqueness of (in short time) of solutions 
of a weakly parabolic equation for sections of a finite dimensional vector bundle.  
The evolution equation $(1.1)$ is regarded as the evolution equation for sections of 
the infinite dimensional trivial vector bundle $M\times V$ over $M$.  
Since $M$ is an infinite dimensional Hilbert manifold, we cannot apply the Hamilton's result to this evolution equation $(1.1)$.  
Also, since $M$ is of infinite dimension, $M$ is not locally compact.  
Thus, we cannot show the existence and the uniqueness of solutions of $(1.1)$ starting from $f$ for a general regularizable 
$C^{\infty}$-immersion $f$ of codimension one.  
Hence, at least, we must impose the finiteness of the cohomogeneity of $f(M)$ to the initial data $f$ in order to use the compactness.  Here 
``the finiteness of the cohomogeneity of $f(M)$'' means that there exists a closed subgroup $\mathcal G$ of the (full) isometry group 
$O(V)\ltimes V$ (which is not a Banach Lie group) of $V$ such that $f(M)$ is $\mathcal G$-invariant and that $f(M)/\mathcal G$ is a finite dimensional 
manifold with singularity.  Furthermore, in order to use the compactness, we must impose $f$ the condition that $f(M)/\mathcal G$ is compact.  

So, we (\cite{Koi2}) considered the following special case.  
Assume that an action of a Hilbert Lie group ${\mathcal G}$ on a Hilbert space $V$ satisifies the following conditions:

\vspace{0.15truecm}

(I)\ \ The action $\mathcal G\curvearrowright V$ is isometric and almost free, where ``almost free'' means that 
the isotropy group of the action at each point is finite;

(II)\ \ All $\mathcal G$-orbits are minimal regularizale submanifolds, that is, they are regularizable submanifold and their regularized 
mean curvature vectors vanish.  


\vspace{0.15truecm}

\noindent
Note that $V/\mathcal G$ is a finite dimensional orbifold by the condition (II).  
Denote by $N$ the orbit space $V/\mathcal G$.  Give $N$ the Riemannian orbimetric $g_N$ such that the orbit map 
$\phi:V\to N$ is a Riemannian orbisubmersion.  
Let $M(\subset V)$ be a ${\mathcal G}$-invariant hypersurface in $V$.  
Furthermore, we assume the following condition:

\vspace{0.15truecm}

(III)\ \ $\overline M:=\phi(M)$ is compact.  

\vspace{0.15truecm}

\noindent
Note that $\overline M$ is a compact hypersurface in $N$.  
Denote by $f$ the inclusion maps of $M$ into $V$ and $\overline f$ that of $\overline M$ into $N$.  
We (\cite{Koi2}) showed that the regularized mean curvature flow starting from $M$ exists uniquely in short time.  
However there were gaps in the statement and the proof.  
In this paper, we close the gaps (see the statement and the proof of Theorem 4.1 (also those of Theorem 3.1)).  
Here we note that the uniqueness of the flow is assured under the $\mathcal G$-invariance of the flow.  

In \cite{Koi2}, we mainly proved that the horizontally strongly convexity is preserved along 
the ($\mathcal G$-invariant) regularized mean curvature flow in the case where $M$ is a ${\mathcal G}$-invariant hypersurface in $V$ 
(see Thoerem 6.1 in \cite{Koi2}), where ``$\mathcal G$-invariance'' of a regularized mean curvature flow means that the regularized mean curvature 
flow consists of $\mathcal G$-invariant regularizable hypersurfaces.  In the statement of Theorem 6.1 of \cite{Koi2}, it is not specfied that 
the regularized mean curvature flow is $\mathcal G$-invariant but the assumption of the $\mathcal G$-invariance of the flow is needed in the statement.  

In this paper, we prove that the $\mathcal G$-invariant regularized mean curvature flow starting from 
a horizontally strongly convex ${\mathcal G}$-invariant hypersurface in $V$ collapses to 
a ${\mathcal G}$-orbit in finite time under some additional conditions (see Theorem A).  
We shall state this collapsing theorem in detail.  
Let $M$ be a hypersurface admitting a (global) unit normal vector field $\xi$.  
Denote by $\overline K$ the maximal sectional curvature of $(N,g_N)$, which is nonnegative because $V$ is flat.  Set $b:=\sqrt{\overline K}$.  
Let $\Sigma$ be the singular set of $(N,g_N)$ and $\{\Sigma_1,\cdots,\Sigma_k\}$ be the set of all connected components of $\Sigma$.  
Set $B^T_r(x):=\{v\in T_xN\,|\,\|v\|\leq r\}$.  
For $x\in N$ and $r>0$, denote by $B_r(x)$ the geodesic ball of radius $r$ centered at $x$.  
Here we note that, even if $x\in\Sigma$, the exponential map $\exp_x:T_xN\to N$ 
is defined in the same manner as the Riemannian manifold-case and $B_r(x)$ is defined by 
$B_r(x):=\exp_x(B^T_r(x))$.  
We assume the following:

\vspace{0.25truecm}

\noindent
$(\ast_1)$ $\overline M$ is included by $B_{\frac{\pi}{b}}(x_0)$ for some $x_0\in N$ and 
$\exp_{x_0}|_{B^T_{\frac{\pi}{b}}(x_0)}$ is injective.  

\vspace{0.25truecm}

Furthermore, we assume the following:
$$b^2(1-\alpha)^{-2/n}\left(\omega_n^{-1}\cdot{\rm Vol}_{\bar g}(\overline M)\right)^{2/n}\leq 1,\leqno{(\ast_2)}$$
where $\alpha$ is a positive constant smaller than one and $\overline g$ denotes the induced metric on 
$\overline M$ and $\omega_n$ denotes the volume of the unit ball in the Euclidean space ${\mathbb R}^n$.  
Denote by $f$ the inclusion map of $M$ to $V$ and $\overline f$ that of $\overline M$ into $N$.  
Let $\{f_t\}_{t\in[0,T)}$ be the 
$\mathcal G$-invariant regularized mean curvature flow starting from $f$.  
Denote by $H_t$ the regularized mean curvature vector of $f_t$ and set $H_t^s:=-\langle H_t,\xi_t\rangle
(={\rm Tr}_r\,((A_t)_{-\xi_t}))$, where 
$\xi_t$ is a unit normal vector field of $f_t$ such that $\xi_0=\xi$ and $t\mapsto\xi_t$ is continuous.  
Set $(H^s_t)_{\min}:=\min_M\,H^s_t$ and $(H^s_t)_{\max}:=\max_M\,H^s_t$.  
Define a constant $L$ by 
$$L:=\mathop{\sup}_{u\in V}\mathop{\max}_{(X_1,\cdots,X_5)\in({\widetilde{\mathcal H}}_1)_u^5}
\vert\langle{\mathcal A}^{\phi}_{X_1}((\widetilde{\nabla}_{X_2}{\mathcal A}^{\phi})_{X_3}X_4),\,X_5\rangle\vert,
\leqno{(\sharp)}$$
where $\widetilde{\mathcal H}$ denotes the horizontal distribution of $\phi$, $(\widetilde{\mathcal H}_1)_u$ denotes 
the set of all unit horizontal vectors of $\phi$ at $u$ and $\mathcal A^{\phi}$ denotes one of the O'Neill's tensors 
defined in \cite{O'N} (see Section 4 about the definition of $\mathcal A^{\phi}$).  
Note that the restriction $\mathcal A^{\phi}|_{\widetilde{\mathcal H}\times\widetilde{\mathcal H}}$ of 
$\mathcal A^{\phi}$ to $\widetilde{\mathcal H}\times\widetilde{\mathcal H}$ is the tensor indicating the obstruction 
of the integrabilty of $\widetilde{\mathcal H}$.  

\vspace{0.5truecm}

{\small
\centerline{
\unitlength 0.1in
\begin{picture}( 46.7700, 25.8900)(  0.5100,-31.2000)
%
\special{pn 8}%
\special{ar 1296 2196 266 232  6.2831853 6.2831853}%
\special{ar 1296 2196 266 232  0.0000000 3.1415927}%
%
\special{pn 8}%
\special{ar 1380 2196 350 634  3.1415927 4.4620977}%
%
\special{pn 8}%
\special{ar 1202 2196 362 634  4.9605235 6.2831853}%
%
\special{pn 20}%
\special{sh 1}%
\special{ar 1296 1590 10 10 0  6.28318530717959E+0000}%
\special{sh 1}%
\special{ar 1296 1590 10 10 0  6.28318530717959E+0000}%
%
\special{pn 8}%
\special{ar 3068 2206 268 232  6.2831853 6.2831853}%
\special{ar 3068 2206 268 232  0.0000000 3.1415927}%
%
\special{pn 8}%
\special{ar 3150 2206 350 634  3.1415927 4.4620977}%
%
\special{pn 8}%
\special{ar 2974 2206 362 634  4.9605235 6.2831853}%
%
\special{pn 8}%
\special{ar 4380 2214 266 232  6.2831853 6.2831853}%
\special{ar 4380 2214 266 232  0.0000000 3.1415927}%
%
\special{pn 8}%
\special{ar 4462 2214 350 634  3.1415927 4.4620977}%
%
\special{pn 8}%
\special{ar 4286 2214 362 634  4.9625019 6.2831853}%
%
\special{pn 8}%
\special{pa 1278 772}%
\special{pa 948 1174}%
\special{fp}%
%
\special{pn 8}%
\special{pa 1278 772}%
\special{pa 1608 1174}%
\special{fp}%
%
\special{pn 8}%
\special{ar 1278 1014 198 72  6.2831853 6.2831853}%
\special{ar 1278 1014 198 72  0.0000000 3.1415927}%
%
\special{pn 8}%
\special{ar 1278 1022 198 72  3.1415927 3.5893538}%
\special{ar 1278 1022 198 72  3.8580106 4.3057718}%
\special{ar 1278 1022 198 72  4.5744285 5.0221897}%
\special{ar 1278 1022 198 72  5.2908464 5.7386076}%
\special{ar 1278 1022 198 72  6.0072643 6.2831853}%
%
\special{pn 8}%
\special{ar 1296 1796 184 54  6.2831853 6.2831853}%
\special{ar 1296 1796 184 54  0.0000000 3.1415927}%
%
\special{pn 8}%
\special{ar 1296 1786 184 54  3.1415927 3.6479218}%
\special{ar 1296 1786 184 54  3.9517192 4.4580483}%
\special{ar 1296 1786 184 54  4.7618458 5.2681749}%
\special{ar 1296 1786 184 54  5.5719724 6.0783015}%
%
\special{pn 8}%
\special{pa 1296 1254}%
\special{pa 1296 1486}%
\special{fp}%
\special{sh 1}%
\special{pa 1296 1486}%
\special{pa 1316 1418}%
\special{pa 1296 1432}%
\special{pa 1276 1418}%
\special{pa 1296 1486}%
\special{fp}%
%
\special{pn 8}%
\special{pa 2206 1360}%
\special{pa 2212 1756}%
\special{pa 2646 1168}%
\special{pa 2638 816}%
\special{pa 2638 816}%
\special{pa 2638 816}%
\special{pa 2206 1360}%
\special{fp}%
%
\special{pn 8}%
\special{pa 2498 1160}%
\special{pa 2510 1188}%
\special{pa 2506 1220}%
\special{pa 2496 1250}%
\special{pa 2480 1278}%
\special{pa 2462 1304}%
\special{pa 2438 1326}%
\special{pa 2414 1348}%
\special{pa 2386 1362}%
\special{pa 2354 1362}%
\special{pa 2338 1336}%
\special{pa 2342 1304}%
\special{pa 2352 1274}%
\special{pa 2366 1246}%
\special{pa 2384 1220}%
\special{pa 2406 1196}%
\special{pa 2430 1176}%
\special{pa 2458 1160}%
\special{pa 2490 1156}%
\special{pa 2498 1160}%
\special{sp}%
%
\special{pn 8}%
\special{pa 3052 1698}%
\special{pa 3062 1726}%
\special{pa 3056 1758}%
\special{pa 3046 1788}%
\special{pa 3030 1816}%
\special{pa 3010 1842}%
\special{pa 2988 1864}%
\special{pa 2962 1882}%
\special{pa 2930 1886}%
\special{pa 2914 1860}%
\special{pa 2914 1828}%
\special{pa 2926 1798}%
\special{pa 2940 1770}%
\special{pa 2958 1744}%
\special{pa 2980 1720}%
\special{pa 3006 1702}%
\special{pa 3036 1692}%
\special{pa 3052 1698}%
\special{sp}%
%
\special{pn 20}%
\special{sh 1}%
\special{ar 2994 1788 10 10 0  6.28318530717959E+0000}%
\special{sh 1}%
\special{ar 2994 1788 10 10 0  6.28318530717959E+0000}%
%
\special{pn 20}%
\special{sh 1}%
\special{ar 2428 1260 10 10 0  6.28318530717959E+0000}%
\special{sh 1}%
\special{ar 2428 1260 10 10 0  6.28318530717959E+0000}%
%
\special{pn 20}%
\special{sh 1}%
\special{ar 1278 772 10 10 0  6.28318530717959E+0000}%
\special{sh 1}%
\special{ar 1278 772 10 10 0  6.28318530717959E+0000}%
%
\special{pn 8}%
\special{pa 4528 1752}%
\special{pa 4534 1784}%
\special{pa 4526 1814}%
\special{pa 4510 1842}%
\special{pa 4492 1868}%
\special{pa 4468 1890}%
\special{pa 4446 1912}%
\special{pa 4422 1934}%
\special{pa 4394 1950}%
\special{pa 4368 1968}%
\special{pa 4340 1984}%
\special{pa 4310 1996}%
\special{pa 4282 2008}%
\special{pa 4250 2016}%
\special{pa 4218 2022}%
\special{pa 4188 2020}%
\special{pa 4160 2004}%
\special{pa 4158 2002}%
\special{sp}%
%
\special{pn 8}%
\special{pa 4158 2004}%
\special{pa 4150 1974}%
\special{pa 4160 1944}%
\special{pa 4174 1916}%
\special{pa 4194 1890}%
\special{pa 4216 1866}%
\special{pa 4240 1846}%
\special{pa 4264 1824}%
\special{pa 4290 1806}%
\special{pa 4316 1788}%
\special{pa 4346 1774}%
\special{pa 4374 1760}%
\special{pa 4404 1748}%
\special{pa 4434 1742}%
\special{pa 4466 1736}%
\special{pa 4498 1738}%
\special{pa 4526 1754}%
\special{pa 4528 1756}%
\special{sp 0.070}%
%
\special{pn 8}%
\special{pa 3616 1198}%
\special{pa 3572 1718}%
\special{pa 4398 1002}%
\special{pa 4430 542}%
\special{pa 4430 542}%
\special{pa 4430 542}%
\special{pa 3616 1198}%
\special{fp}%
%
\special{pn 8}%
\special{pa 4188 932}%
\special{pa 4212 952}%
\special{pa 4224 982}%
\special{pa 4224 1014}%
\special{pa 4214 1046}%
\special{pa 4198 1074}%
\special{pa 4182 1102}%
\special{pa 4164 1126}%
\special{pa 4142 1150}%
\special{pa 4120 1174}%
\special{pa 4098 1196}%
\special{pa 4074 1216}%
\special{pa 4046 1234}%
\special{pa 4020 1252}%
\special{pa 3994 1270}%
\special{pa 3966 1286}%
\special{pa 3936 1298}%
\special{pa 3906 1308}%
\special{pa 3874 1314}%
\special{pa 3842 1320}%
\special{pa 3810 1318}%
\special{pa 3782 1304}%
\special{pa 3768 1276}%
\special{pa 3764 1244}%
\special{pa 3772 1214}%
\special{pa 3786 1184}%
\special{pa 3800 1154}%
\special{pa 3818 1128}%
\special{pa 3840 1106}%
\special{pa 3860 1082}%
\special{pa 3884 1058}%
\special{pa 3908 1038}%
\special{pa 3934 1018}%
\special{pa 3960 1000}%
\special{pa 3986 982}%
\special{pa 4014 966}%
\special{pa 4042 952}%
\special{pa 4072 940}%
\special{pa 4104 934}%
\special{pa 4136 928}%
\special{pa 4166 928}%
\special{pa 4188 932}%
\special{sp}%
%
\special{pn 20}%
\special{sh 1}%
\special{ar 3994 1128 10 10 0  6.28318530717959E+0000}%
\special{sh 1}%
\special{ar 3994 1128 10 10 0  6.28318530717959E+0000}%
%
\special{pn 20}%
\special{sh 1}%
\special{ar 4298 1664 10 10 0  6.28318530717959E+0000}%
\special{sh 1}%
\special{ar 4298 1664 10 10 0  6.28318530717959E+0000}%
%
\special{pn 8}%
\special{pa 4096 1352}%
\special{pa 4242 1592}%
\special{fp}%
\special{sh 1}%
\special{pa 4242 1592}%
\special{pa 4224 1526}%
\special{pa 4214 1546}%
\special{pa 4190 1546}%
\special{pa 4242 1592}%
\special{fp}%
\put(13.4200,-13.0700){\makebox(0,0)[lt]{$\exp_{x_0}$}}%
\put(27.2800,-15.4700){\makebox(0,0)[lb]{$\exp_{x_0}$}}%
%
\special{pn 8}%
\special{pa 2554 1432}%
\special{pa 2866 1708}%
\special{fp}%
\special{sh 1}%
\special{pa 2866 1708}%
\special{pa 2830 1650}%
\special{pa 2826 1674}%
\special{pa 2804 1680}%
\special{pa 2866 1708}%
\special{fp}%
\put(42.2400,-14.7600){\makebox(0,0)[lb]{$\exp_{x_0}$}}%
\put(11.6800,-7.0100){\makebox(0,0)[lb]{$T_{x_0}N$}}%
\put(27.2800,-10.9300){\makebox(0,0)[lb]{$T_{x_0}N$}}%
\put(45.3600,-8.5200){\makebox(0,0)[lb]{$T_{x_0}N$}}%
\put(16.7300,-22.1600){\makebox(0,0)[lt]{$N$}}%
\put(34.1600,-22.2500){\makebox(0,0)[lt]{$N$}}%
\put(47.2800,-22.3400){\makebox(0,0)[lt]{$N$}}%
%
\special{pn 4}%
\special{pa 1182 902}%
\special{pa 1286 800}%
\special{dt 0.027}%
\special{pa 1098 1038}%
\special{pa 1312 828}%
\special{dt 0.027}%
\special{pa 1132 1058}%
\special{pa 1340 856}%
\special{dt 0.027}%
\special{pa 1174 1070}%
\special{pa 1364 884}%
\special{dt 0.027}%
\special{pa 1222 1076}%
\special{pa 1390 914}%
\special{dt 0.027}%
\special{pa 1270 1084}%
\special{pa 1416 942}%
\special{dt 0.027}%
\special{pa 1334 1076}%
\special{pa 1442 970}%
\special{dt 0.027}%
\special{pa 1404 1060}%
\special{pa 1468 998}%
\special{dt 0.027}%
%
\special{pn 4}%
\special{pa 1146 1740}%
\special{pa 1284 1604}%
\special{dt 0.027}%
\special{pa 1296 1592}%
\special{pa 1302 1588}%
\special{dt 0.027}%
\special{pa 1128 1810}%
\special{pa 1334 1610}%
\special{dt 0.027}%
\special{pa 1164 1828}%
\special{pa 1362 1638}%
\special{dt 0.027}%
\special{pa 1210 1838}%
\special{pa 1392 1660}%
\special{dt 0.027}%
\special{pa 1256 1846}%
\special{pa 1416 1692}%
\special{dt 0.027}%
\special{pa 1306 1852}%
\special{pa 1438 1722}%
\special{dt 0.027}%
\special{pa 1370 1842}%
\special{pa 1458 1756}%
\special{dt 0.027}%
\special{pa 1440 1826}%
\special{pa 1480 1788}%
\special{dt 0.027}%
%
\special{pn 4}%
\special{pa 2388 1228}%
\special{pa 2444 1174}%
\special{dt 0.027}%
\special{pa 2344 1324}%
\special{pa 2504 1170}%
\special{dt 0.027}%
\special{pa 2370 1352}%
\special{pa 2498 1228}%
\special{dt 0.027}%
%
\special{pn 4}%
\special{pa 2924 1830}%
\special{pa 3056 1702}%
\special{dt 0.027}%
\special{pa 2928 1880}%
\special{pa 3058 1752}%
\special{dt 0.027}%
%
\special{pn 4}%
\special{pa 3774 1218}%
\special{pa 4050 950}%
\special{dt 0.027}%
\special{pa 3776 1270}%
\special{pa 4114 942}%
\special{dt 0.027}%
\special{pa 3800 1300}%
\special{pa 4176 934}%
\special{dt 0.027}%
\special{pa 3840 1314}%
\special{pa 4210 956}%
\special{dt 0.027}%
\special{pa 3908 1302}%
\special{pa 4218 1002}%
\special{dt 0.027}%
\special{pa 4014 1254}%
\special{pa 4208 1064}%
\special{dt 0.027}%
%
\special{pn 4}%
\special{pa 4206 1816}%
\special{pa 4404 1622}%
\special{dt 0.027}%
\special{pa 4164 1908}%
\special{pa 4436 1646}%
\special{dt 0.027}%
\special{pa 4146 1980}%
\special{pa 4464 1670}%
\special{dt 0.027}%
\special{pa 4176 2006}%
\special{pa 4492 1698}%
\special{dt 0.027}%
\special{pa 4216 2020}%
\special{pa 4512 1732}%
\special{dt 0.027}%
\special{pa 4276 2014}%
\special{pa 4532 1766}%
\special{dt 0.027}%
\special{pa 4390 1958}%
\special{pa 4520 1830}%
\special{dt 0.027}%
\put(17.5000,-28.8000){\makebox(0,0)[lt]{{\small(I),(II) $\,\,:\,$ $\exp_{x_0}|_{B^T_r(x_0)}:$ injective}}}%
\put(14.8000,-17.3500){\makebox(0,0)[lb]{$B_r(x_0)$}}%
\put(45.5400,-17.1700){\makebox(0,0)[lb]{$B_r(x_0)$}}%
\put(26.3600,-20.0200){\makebox(0,0)[rt]{$B_r(x_0)$}}%
%
\special{pn 8}%
\special{ar 2738 2092 276 250  3.6221163 3.6678306}%
\special{ar 2738 2092 276 250  3.8049734 3.8506877}%
\special{ar 2738 2092 276 250  3.9878306 4.0335449}%
\special{ar 2738 2092 276 250  4.1706877 4.2164020}%
\special{ar 2738 2092 276 250  4.3535449 4.3992592}%
\special{ar 2738 2092 276 250  4.5364020 4.5821163}%
%
\special{pn 8}%
\special{pa 2746 1842}%
\special{pa 2948 1842}%
\special{dt 0.045}%
\special{sh 1}%
\special{pa 2948 1842}%
\special{pa 2882 1822}%
\special{pa 2896 1842}%
\special{pa 2882 1862}%
\special{pa 2948 1842}%
\special{fp}%
\put(11.3100,-8.9700){\makebox(0,0)[rb]{$B_r^T(x_0)$}}%
\put(23.7900,-10.7500){\makebox(0,0)[rb]{$B_r^T(x_0)$}}%
\put(39.3000,-8.2000){\makebox(0,0)[rb]{$B_r^T(x_0)$}}%
%
\special{pn 8}%
\special{pa 2252 1112}%
\special{pa 2454 1218}%
\special{dt 0.045}%
\special{sh 1}%
\special{pa 2454 1218}%
\special{pa 2404 1170}%
\special{pa 2406 1194}%
\special{pa 2386 1204}%
\special{pa 2454 1218}%
\special{fp}%
%
\special{pn 8}%
\special{pa 3790 860}%
\special{pa 4020 1030}%
\special{dt 0.045}%
\special{sh 1}%
\special{pa 4020 1030}%
\special{pa 3978 974}%
\special{pa 3976 998}%
\special{pa 3954 1006}%
\special{pa 4020 1030}%
\special{fp}%
\put(12.0000,-25.8000){\makebox(0,0)[lt]{(I)}}%
\put(29.2500,-25.8800){\makebox(0,0)[lt]{(II)}}%
\put(42.3700,-25.8800){\makebox(0,0)[lt]{(III)}}%
\put(17.6000,-31.2000){\makebox(0,0)[lt]{{\small(III) $\,\,:\,$ $\exp_{x_0}|_{B^T_r(x_0)}:$ not injective}}}%
%
\special{pn 13}%
\special{ar 4180 2340 290 1040  5.4926414 5.8632361}%
%
\special{pn 13}%
\special{ar 4730 1510 460 410  2.0808587 2.7545726}%
%
\special{pn 13}%
\special{pa 4000 1130}%
\special{pa 4120 1180}%
\special{fp}%
%
\special{pn 13}%
\special{pa 4000 1120}%
\special{pa 4200 940}%
\special{fp}%
%
\special{pn 13}%
\special{ar 4440 1770 190 190  3.8239092 4.3616610}%
\end{picture}%
\hspace{2.5truecm}}
}

\vspace{0.5truecm}

\centerline{{\bf Figure 1.1$\,\,:\,\,$ The injectivity of $\exp_{x_0}$}}

\vspace{0.5truecm}

In this paper, we prove the following collapsing theorem.  

\vspace{0.5truecm}

\noindent
{\bf Theorem A.} {\sl 
Assume that the initial $\mathcal G$-invariant hypersurface $f:M\hookrightarrow V$ satisfies 
the above conditions $(\ast_1),\,(\ast_2)$ and 
$(H^s_0)^2(h_{\mathcal H})_0>2n^2L(g_{\mathcal H})_0$, where $(g_{\mathcal H})_0$ (resp. $(h_{\mathcal H})_0$) 
denotes the horizontal component of the induced metric (resp. the scalar-valued second fundamental form) of $f$.  
Then, for a $\mathcal G$-invariant regularized mean curvature flow $\{f_t\}_{t\in[0,T)}$ 
($T\,:\,$ the maximal time) starting from $f$, the following statements 
{\rm (i)} and {\rm (ii)} hold:

{\rm(i)} $T<\infty$ and 
$f_t(M)$ collapses to a $\mathcal G$-orbit as $t\to T$.  

{\rm (ii)} $\displaystyle{\lim_{t\to T}\frac{(H^s_t)_{\max}}
{(H^s_t)_{\min}}=1}$ holds.
}


\vspace{0.5truecm}

\noindent
{\it Remark 1.1.} ``$\displaystyle{\lim_{t\to T}\frac{(H^s_t)_{\max}}
{(H^s_t)_{\min}}=1}$'' implies that 
$f_t(M)$ converges to an infinitesimal constant tube over some $G$-orbit as $t\to T$ 
(or equivalently, $\phi(f_t(M))$ converges to a round point(=an infinitesimal round sphere) as $t\to T$) 
(see Figure 1.2).  

\vspace{0.5truecm}

In Section 2, we recall the definition of the the regularized mean curvature flow and, in Section 3, 
we discuss the existence and the uniqueness of mean curvature flows starting from a compact orbifold in a Riemannian orbifold.  
In the first-half part of Section 4, we give a new proof of the existence and the uniqueness of a $\mathcal G$-invariant regularized 
mean curvature flow starting from a $\mathcal G$-invariant regularizable hypersurface satisfying the condition (III) in a Hilbert space $V$ 
equipped with a Hilbert Lie group action $\mathcal G\curvearrowright V$ satisfying the conditions (I) and (II).  
In the second-half part of the section, we prepare the evolution equations for some basic geometric quantities along the $\mathcal G$-invariant 
regularized mean curvature flow.  In Section 5, we prove the Sobolev inequality for Riemannian suborbifolds.  
Sections 6-8 is devoted to prove Theorem A.  
In Section 9, we we state a vision for applying the study of the regularized mean curvature flow to the gauge theory.  
In more detail, we state a vision for find and study interesting flows of hypersurfaces in the Yang-Milles (or self-dual) muduli space from 
regularlized mean curvature flows in a Hilbert space.

\vspace{0.5truecm}

{\small
\centerline{
\unitlength 0.1in
\begin{picture}( 20.1400, 18.0000)(  9.9000,-25.9000)
%
\special{pn 8}%
\special{ar 2304 2418 582 172  0.0000000 6.2831853}%
%
\special{pn 8}%
\special{ar 2304 1810 582 172  0.0000000 6.2831853}%
%
\special{pn 8}%
\special{ar 2304 1202 582 172  0.0000000 6.2831853}%
%
\special{pn 8}%
\special{pa 2886 1202}%
\special{pa 2886 1810}%
\special{fp}%
%
\special{pn 8}%
\special{pa 1710 1202}%
\special{pa 1710 1810}%
\special{fp}%
%
\special{pn 8}%
\special{pa 2822 2024}%
\special{pa 2822 2246}%
\special{fp}%
\special{sh 1}%
\special{pa 2822 2246}%
\special{pa 2842 2180}%
\special{pa 2822 2194}%
\special{pa 2802 2180}%
\special{pa 2822 2246}%
\special{fp}%
\put(29.1800,-20.7300){\makebox(0,0)[lt]{$\phi$}}%
\put(29.9300,-24.0800){\makebox(0,0)[lt]{$V/\mathcal G$}}%
\put(30.0400,-14.0500){\makebox(0,0)[lt]{$V$}}%
%
\special{pn 8}%
\special{ar 2368 1212 260 62  1.5707963 4.7123890}%
%
\special{pn 8}%
\special{ar 2368 1212 108 62  4.7123890 6.2831853}%
\special{ar 2368 1212 108 62  0.0000000 1.5707963}%
%
\special{pn 8}%
\special{ar 2368 1800 260 62  1.5707963 4.7123890}%
%
\special{pn 8}%
\special{ar 2368 1800 108 62  4.7123890 6.2831853}%
\special{ar 2368 1800 108 62  0.0000000 1.5707963}%
%
\special{pn 8}%
\special{ar 2368 2408 260 60  1.5707963 4.7123890}%
%
\special{pn 8}%
\special{ar 2368 2408 108 60  4.7123890 6.2831853}%
\special{ar 2368 2408 108 60  0.0000000 1.5707963}%
%
\special{pn 8}%
\special{pa 2110 1212}%
\special{pa 2110 1790}%
\special{fp}%
%
\special{pn 8}%
\special{pa 2476 1222}%
\special{pa 2476 1790}%
\special{fp}%
%
\special{pn 8}%
\special{pa 2320 1210}%
\special{pa 2320 1798}%
\special{fp}%
%
\special{pn 13}%
\special{pa 2152 1486}%
\special{pa 2292 1486}%
\special{fp}%
\special{sh 1}%
\special{pa 2292 1486}%
\special{pa 2226 1466}%
\special{pa 2240 1486}%
\special{pa 2226 1506}%
\special{pa 2292 1486}%
\special{fp}%
%
\special{pn 13}%
\special{pa 2142 1678}%
\special{pa 2282 1678}%
\special{fp}%
\special{sh 1}%
\special{pa 2282 1678}%
\special{pa 2216 1658}%
\special{pa 2230 1678}%
\special{pa 2216 1698}%
\special{pa 2282 1678}%
\special{fp}%
%
\special{pn 13}%
\special{pa 2152 1334}%
\special{pa 2292 1334}%
\special{fp}%
\special{sh 1}%
\special{pa 2292 1334}%
\special{pa 2226 1314}%
\special{pa 2240 1334}%
\special{pa 2226 1354}%
\special{pa 2292 1334}%
\special{fp}%
%
\special{pn 13}%
\special{pa 2450 1490}%
\special{pa 2350 1490}%
\special{fp}%
\special{sh 1}%
\special{pa 2350 1490}%
\special{pa 2418 1510}%
\special{pa 2404 1490}%
\special{pa 2418 1470}%
\special{pa 2350 1490}%
\special{fp}%
%
\special{pn 13}%
\special{pa 2460 1330}%
\special{pa 2360 1330}%
\special{fp}%
\special{sh 1}%
\special{pa 2360 1330}%
\special{pa 2428 1350}%
\special{pa 2414 1330}%
\special{pa 2428 1310}%
\special{pa 2360 1330}%
\special{fp}%
%
\special{pn 13}%
\special{pa 2460 1680}%
\special{pa 2360 1680}%
\special{fp}%
\special{sh 1}%
\special{pa 2360 1680}%
\special{pa 2428 1700}%
\special{pa 2414 1680}%
\special{pa 2428 1660}%
\special{pa 2360 1680}%
\special{fp}%
%
\special{pn 8}%
\special{pa 1840 2160}%
\special{pa 2170 2370}%
\special{fp}%
\special{sh 1}%
\special{pa 2170 2370}%
\special{pa 2124 2318}%
\special{pa 2126 2342}%
\special{pa 2104 2352}%
\special{pa 2170 2370}%
\special{fp}%
\put(28.1000,-9.6000){\makebox(0,0)[lb]{$M$}}%
%
\special{pn 8}%
\special{pa 2800 1010}%
\special{pa 2420 1590}%
\special{fp}%
\special{sh 1}%
\special{pa 2420 1590}%
\special{pa 2474 1546}%
\special{pa 2450 1546}%
\special{pa 2440 1524}%
\special{pa 2420 1590}%
\special{fp}%
\put(18.0000,-21.9000){\makebox(0,0)[rb]{$\phi(M)$}}%
%
\special{pn 8}%
\special{pa 2110 1820}%
\special{pa 2110 2410}%
\special{dt 0.045}%
%
\special{pn 8}%
\special{pa 2480 1830}%
\special{pa 2480 2390}%
\special{dt 0.045}%
%
\special{pn 8}%
\special{pa 2320 1800}%
\special{pa 2320 2400}%
\special{dt 0.045}%
%
\special{pn 20}%
\special{sh 1}%
\special{ar 2320 2410 10 10 0  6.28318530717959E+0000}%
\special{sh 1}%
\special{ar 2320 2410 10 10 0  6.28318530717959E+0000}%
%
\special{pn 8}%
\special{ar 2320 2410 50 30  0.0000000 6.2831853}%
%
\special{pn 8}%
\special{pa 2370 2410}%
\special{pa 2370 1210}%
\special{dt 0.045}%
%
\special{pn 8}%
\special{pa 2280 2420}%
\special{pa 2280 1210}%
\special{dt 0.045}%
%
\special{pn 8}%
\special{ar 2320 1210 50 30  0.0000000 6.2831853}%
%
\special{pn 8}%
\special{ar 2320 1800 50 30  0.0000000 6.2831853}%
\end{picture}%
\hspace{2.75truecm}}
}

\vspace{0.5truecm}

\centerline{{\bf Figure 1.2$\,\,:\,\,$ Collapse to a $\mathcal G$-orbit}}


\section{The regularized mean curvature flow} 
Let $f$ be an immersion of an (infinite dimensional) Hilbert manifold $M$ 
into a Hilbert space $V$ and $A$ the shape tensor of $f$.  
If ${\rm codim}\,M\,<\,\infty$, 
if the differential of the normal exponential map $\exp^{\perp}$ of $f$ at each point of $M$ 
is a Fredholm operator and if the restriction $\exp^{\perp}$ to the unit normal ball bundle of $f$ 
is proper, then $M$ is called a {\it proper Fredholm submanifold}.
In this paper, we then call $f$ a {\it proper Fredholm immersion}.  
Then the shape operator $A_v$ is a compact operator for each normal vector $v$ of $M$.  
Furthermore, if, for each normal vector $v$ of $M$, the regularized trace 
${\rm Tr}_r\,A_v$ and ${\rm Tr}\,A_v^2$ exist, then $M$ is called 
{\it regularizable submanifold}, where 
${\rm Tr}_r\,A_v$ is defined by 
${\rm Tr}_r\,A_v:=\sum\limits_{i=1}^{\infty}(\mu^+_i+\mu^-_i)$ 
($\mu^-_1\leq\mu^-_2\leq\cdots\leq 0\leq\cdots\leq\mu^+_2\leq\mu^+_1\,:\,$
the spectrum of $A_v$).  
In this paper, we then call $f$ {\it regularizable immersion}.  
If ${\rm Tr}_r\,A_v=0$ holds for any $v\in T^{\perp}M$, then $f$ is said to be {\it minimal}.  
If $f$ is a regulalizable immersion and if $\rho_u:v\mapsto{\rm Tr}_r\,(A_u)_v$ ($v\in T^{\perp}_uM$) is linear 
for any $u\in M$, 
then the {\it regularized mean curvature vector} $H$ of $f$ is defined as 
the normal vector field satisfying 
$\langle H_u,v\rangle={\rm Tr}_r\,(A_u)_v\,\,(\forall\,v\in T_u^{\perp}M)$ ($u\in M$), where 
$\langle\,\,,\,\,\rangle$ denotes the inner product of $V$ and $T_u^{\perp}M$ denotes the normal space of $f$ 
at $u$.  

\hspace{0.5truecm}

\noindent
{\it Example 2.1.}\ 
We consider the case where $f$ is isoparametric.  Then 
First we recall the notion of an isoparametric submanifold in a Hilbert space.  
If the normal connection of $f$ is flat and if the principal curvatures of $f$ for $v$ are constant for 
any parallel normal vector field $v$, then it is called an {\it isoparametric submanifold}.  
Then, by analyzing the focal structure of $f$, we can show that the set $\Lambda$ of all the principal curvatures 
of $f$ is given by 
$$\Lambda=\mathop{\cup}_{a=1}^k\left\{\left.\frac{\lambda_a}{1+b_aj}\,\right|\,j\in\mathbb Z\right\},$$
where $\lambda_a$'s are parallel sections of the normal bundle $T^{\perp}M$ and $b_a$'s are positive constants 
greater than one.  
See 
the first-half part of the proof of Theorem A in \cite{Koi1} 
about the proof of this fact.  
Note that, even if Theorem A in \cite{Koi1} is a result for isoparametric submanifolds in a Hilbert space arising from 
equifocal submanifolds in symmetric space of compact type, the first-half part of the proof is discussed for general 
isoparametric submanifolds in a Hilbert space.  
Hence the spectrum ${\rm Spec}\,(A_u)_v$ of the shape operator $(A_u)_v$ for each normal vector $v$ of $f$ at $u\in M$ is 
given by 
$${\rm Spec}\,(A_u)_v=\mathop{\cup}_{a=1}^k\left\{\left.\frac{(\lambda_a)_u(v)}{1+b_aj}\,\right|\,
j\in\mathbb Z\right\}.$$
Hence the regularized trace ${\rm Tr}_r\,(A_u)_v$ is given by 
$${\rm Tr}_r\,(A_u)_v=\sum_{a=1}^k\sum_{j=1}^{\infty}\frac{(\lambda_a)_u(v)}{1+b_aj}.$$
From this fact, it directly follows that $\rho_u$ is linear.  

\vspace{0.5truecm}

\noindent
{\it Example 2.2.}\ 
We consider the case where $f$ is a hypersurface.  
Then, since the normal space of $M$ is of dimension one, $\rho_u$ is linear for each point $u\in M$.  

\vspace{0.5truecm}

We consider the case where $f$ is a hypersurface and it admits a global unit normal vector field.  
Fix a global unit normal vector field $\xi$.  Then we call 
${\rm Tr}_r\,A_{-\xi}(=-\langle H,\xi\rangle)$ the {\it regularized mean curvature} of $f$ and denote it by $H^s$.  
Also, we call $-A_{\xi}$ the {\it shape operator} and denote it by the same symbol $A$.  

\vspace{0.5truecm}

\noindent
{\it Remark 2.1.} 
In the research of the mean curvature flow starting from strictly convex hypersurfaces, 
it is general to take the outward unit normal vector field as the unit normal vector field $\xi$ 
and $-A_{\xi}$ as the shape operator and $-\langle H,\xi\rangle$ as the mean curvature.  
Hence we take the shape operator $A$ and the regularized mean curvature $H^s$ as above.  

\vspace{0.5truecm}

Let $\{f_t\}_{t\in[0,T)}$ be a $C^{\infty}$-family of 
regularizable immersions of $M$ into $V$.  Denote by $H_t$ the 
regularized mean curvature vector of $f_t$.  
Define a map $F:M\times[0,T)\to V$ by 
$F(x,t):=f_t(x)$ ($(x,t)\in M\times[0,T)$).  If 
$\frac{\partial F}{\partial t}=H_t$ holds, then we call $\{f_t\}_{t\in[0,T)}$ the {\it regularized mean curvature flow}.  
We cannot show that there uniquely exists a regularized mean curvature flow satrting from $f$  
for any $C^{\infty}$-regularizable immersion $f:M\hookrightarrow V$ because $M$ is not compact.  
However, for a $\mathcal G$-invariant regularizable immersion $f:M\hookrightarrow V$ with the compact quotient $f(M)/\mathcal G$ 
in a Hilbert space $V$ equipped with a special Hilbert Lie group action $\mathcal G\curvearrowright V$, 
it is shown that there uniquely exists a $\mathcal G$-invariant regularized mean curvature flow starting from $f$ (see Theorem 4.1).  

\section{The mean curvature flow in Riemannian orbifolds}
The basic notions for a Riemannian orbifold and a suborbifold were defined in [AK,BB,GKP,Sa,Sh,Th].  
In \cite{Koi2}, we introduced the notion of the mean curvaure flow starting from a suborbifold 
in a Riemannian orbifold.  We shall recall this notion shortly.  
Let $M$ be a paracompact Hausdorff space and 
${\mathcal O}:=\{(U_{\lambda},\varphi_{\lambda},\widehat U_{\lambda}/\Gamma_{\lambda})\,\vert\,\lambda\in\Lambda\}$ 
an $n$-{\it dimensional} $C^k$-{\it orbifold atlas} of $M$, that is, a family satisfying the following condition 
(i)-(iv):

\vspace{0.15truecm}

(i) $\{U_{\lambda}\,\vert\,\lambda\in\Lambda\}$ is an open covering of $M$, 

(ii) $\widehat U_{\lambda}$ is an open set of ${\Bbb R}^n$ and $\Gamma_{\lambda}$ is a finite subgroup of 
the $C^k$-diffeo-

morphism group ${\rm Diff}^k(\widehat U_{\lambda})$ of $\widehat U_{\lambda}$,

(iii) $\varphi_{\lambda}$ is a homeomorphism of $U_{\lambda}$ onto $\widehat U_{\lambda}/\Gamma_{\lambda}$, 

(iv) for $\lambda,\mu\in\Lambda$ with $U_{\lambda}\cap U_{\mu}\not=\emptyset$, 
there exists $(U_{\nu},\phi_{\nu},\widehat U_{\nu}/\Gamma_{\nu})\in\mathcal O$ 

such that $C^k$-embeddings $\rho_{\lambda}:\widehat U_{\nu}\hookrightarrow\widehat U_{\lambda}$ 
and $\rho_{\mu}:\widehat U_{\nu}\hookrightarrow\widehat U_{\mu}$ satisfying 

$\varphi_{\lambda}^{-1}\circ\pi_{\Gamma_{\lambda}}\circ\rho_{\lambda}
=\phi_{\nu}^{-1}\circ\pi_{\Gamma_{\nu}}$ and 
$\varphi_{\mu}^{-1}\circ\pi_{\Gamma_{\mu}}\circ\rho_{\mu}
=\phi_{\nu}^{-1}\circ\pi_{\Gamma_{\nu}}$, where 

$\pi_{\Gamma_{\lambda}},\,\pi_{\Gamma_{\mu}}$ and $\pi_{\Gamma_{\nu}}$ are 
the orbit maps of $\Gamma_{\lambda},\,\Gamma_{\mu}$ and $\Gamma_{\nu}$, respectively.  

\vspace{0.15truecm}

\noindent
The pair $(M,{\cal O})$ is called an $n$-{\it dimensional} $C^k$-{\it orbifold}.  
and each $(U_{\lambda},\varphi_{\lambda},\widehat U_{\lambda}/\Gamma_{\lambda})$ is called 
an {\it orbifold chart}.  
Let $(U_{\lambda},\varphi_{\lambda},\widehat U_{\lambda}/\Gamma_{\lambda})$ be 
an orbifold chart around $x\in M$.  Then the group 
$(\Gamma_{\lambda})_{\widehat x}:=\{b\in\Gamma_{\lambda}\,\vert\,b(\widehat x)=\widehat x\}$ is 
unique for $x$ up to the conjugation, where $\widehat x$ is a point of $\widehat U_{\lambda}$ 
with $(\varphi_{\lambda}^{-1}\circ\pi_{\Gamma_{\lambda}})(\widehat x)=x$.  
Denote by $(\Gamma_{\lambda})_x$ the conjugate class of this group $(\Gamma_{\lambda})_{\widehat x}$, 
This conjugate class is called the {\it local group at} $x$.  
If the local group at $x$ is not trivial, then $x$ is called a {\it singular point} of 
$(M,{\cal O})$.  
Denote by ${\rm Sing}(M,{\cal O})$ (or ${\rm Sing}(M)$) the set of all singular points of $(M,{\cal O})$.  
This set ${\rm Sing}(M,{\cal O})$ is called the {\it singular set} of $(M,{\cal O})$.  

{Let $(M,{\cal O}_M)$ and $(N,{\cal O}_N)$ be orbifolds, and $f$ a map from $M$ to $N$.  
If, for each $x\in M$ and each pair of an orbifold chart 
$(U_{\lambda},\varphi_{\lambda},\widehat U_{\lambda}/\Gamma_{\lambda})$ of $(M,{\cal O}_M)$ 
around $x$ and an orbifold chart 
$(V_{\mu},\psi_{\mu},\widehat V_{\mu}/\Gamma'_{\mu})$ of $(N,{\cal O}_N)$ 
around $f(x)$ ($f(U_{\lambda})\subset V_{\mu}$), 
there exists a $C^k$-map $\widehat f_{\lambda,\mu}:\widehat U_{\lambda}\to\widehat V_{\mu}$ with 
$f\circ\varphi_{\lambda}^{-1}\circ\pi_{\Gamma_{\lambda}}
=\psi_{\mu}^{-1}\circ\pi_{\Gamma'_{\mu}}\circ\widehat f_{\lambda,\mu}$, then 
$f$ is called a $C^k$-{\it orbimap} (or simply a $C^k$-{\it map}).  Also 
$\widehat f_{\lambda,\mu}$ is called a {\it local lift} of $f$ with respect to 
$(U_{\lambda},\varphi_{\lambda},\widehat U_{\lambda}/\Gamma_{\lambda})$ and 
$(V_{\mu},\psi_{\mu},\widehat V_{\mu}/\Gamma'_{\mu})$.  
Furthermore, if each local lift $\widehat f_{\lambda,\mu}$ is an immersion, then 
$f$ is called a $C^k$-{\it orbiimmersion} (or simply a $C^k$-{\it immersion}) and 
$(M,{\cal O}_M)$ is called a $C^k$-{\it (immersed) suborbifold} in $(N,{\cal O}_N)$.  
In the sequel, we assume that $r=\infty$.  
If a $(0,2)$-orbitensor field $g$ of class $C^k$ on $(M,{\cal O}_M)$ is positive definite and symmetric, then 
we call $g$ a $C^k$-{\it Riemannian orbimetric} and $(M,{\cal O}_M,g)$ a $C^k$-{\it Riemannian orbifold}.  
See Section 3 of \cite{Koi2} about the definition of $(0,2)$-orbitensor field of class $C^k$.  

Let $\{f_t\}_{t\in[0,T)}$ be a $C^{\infty}$-family of 
$C^{\infty}$-orbiimmersions of a $C^{\infty}$-orbifold $(M,{\cal O}_M)$ into 
a $C^{\infty}$-Riemannian orbifold $(N,{\cal O}_N,\widetilde g)$.  
Assume that, for each $(x_0,t_0)\in M\times[0,T)$ and each pair of an orbifold chart 
$(U_{\lambda},\varphi_{\lambda},\widehat U_{\lambda}/\Gamma_{\lambda})$ of $(M,{\cal O}_M)$ around 
$x_0$ and an orbifold chart $(V_{\mu},\varphi_{\mu},\widehat V_{\mu}/\Gamma'_{\mu})$ of 
$(N,{\cal O}_N)$ around $f_{t_0}(x_0)$ such that $f_t(U_{\lambda})\subset V_{\mu}$ for any 
$t\in[t_0,t_0+\varepsilon)$ ($\varepsilon:$ a sufficiently small positive number), 
there exists local lifts $(\widehat f_t)_{\lambda,\mu}:\widehat U_{\lambda}\to\widehat V_{\mu}$ of 
$f_t$ ($t\in[t_0,t_0+\varepsilon)$) such that they give the mean curvature flow in 
$(\widehat V_{\mu},\widetilde g^{\wedge}_{\mu})$, where $\widetilde g^{\wedge}_{\mu}$ is the local lift of $g$ to $\widehat V_{\mu}$.  
Then we call $f_t$ ($0\leq t<T$) the {\it mean curvature flow} in $(N,{\cal O}_N,\widetilde g)$.  

In \cite{Koi2}, we proved the existence and the uniqueness theorem of a mean curvature flow starting from 
a $C^{\infty}$-orbiimmersion $f$ of a compact $C^{\infty}$-orbifold into a $C^{\infty}$-Riemannian orbifold (see Theorem 3.1 in \cite{Koi2}).  
However, there was a gap in the proof.  Hence we shall close the gap.  
The mean curvature flow equation is the same kind of partial differential equation as the Ricci flow equation.  
For the existence and the uniqueness of solutions of the Ricci flow equation in a compact orbifold, the following fact is known.  
According to Subsection 5.2 of \cite{KL}, it is shown that, for any $C^{\infty}$-orbimetric $g$ on a compact orbifold $M$, 
there uniquely exists a Ricci flow starting from $g$ in short time.  The method of the proof is as follows.  
The existence and the uniqueness of solutions of the Ricci flow equation in short time is reduced to 
those of a standard quasi-linear parabolic partial differential equation called the {\it Ricci-de Turck equation} by the de Turck trick.  
Since the Ricci-de Turck equation is a standard quasi-linear partial differential equation, it is shown that, 
in the case where $M$ is a compact manifold (without boundary), there uniquely exists a solution of the Ricci-de Turck equation having $g$ 
as the initial data in short time for any $C^{\infty}$-Riemannian metric $g$.  
Hence, in this case, it is shown that, for any $C^{\infty}$-Riemannian metric $g$, there uniquely exists a Ricci flow starting from $g$ in short time.  
As stated in Subsection 5.2 of \cite{KL}, even if $M$ is a compact orbifold (not manifold), it is shown similary that, for any $C^{\infty}$-Riemannian 
orbimetric $g$ on $M$, there uniquely exists a Ricci flow starting from $g$ in short time.  

We prove the following statement by applying this method of the proof to the case of the mean curvature flow starting from 
a $C^{\infty}$-orbiimmersion of a compact orbifold $M$ into a Riemannian orbifold $(N,\widetilde g)$.  

\vspace{0.5truecm}

\noindent
{\bf Theorem 3.1.}\ {\sl Let $(M,\mathcal O_M)$ be a compact orbifold and $(N,\mathcal O_N,\widetilde g)$ a Riemannian orbifold.  
For any $C^{\infty}$-orbiimmersion $f$ of $M$ into $N$, there uniquely exists a mean curvature flow starting $f$ in short time.}  

\vspace{0.5truecm}

\noindent
{\it Proof.}\ 
First we consider the case where $M$ and $N$ are manifolds.  
Then the existence and the uniqueness of solutions of the mean curvature flow equation on a compact manifold (without boundary) in short time 
is reduced to those of a standard quasi-linear parabolic partial differential equation called the {\it mean curvature-de Turck equation} 
by the de Turck trick (see \cite{Z}, \cite{CY} and \cite{Koi3} etc. for example).  
Let's recall the definition of the mean curvature-de Turck equation.  
Let $\{f_t\}_{t\in[0,T)}$ be a $C^{\infty}$-family of immersions of $M$ into $(N,\widetilde g)$.  
Set $g_t:=f_t^{\ast}\widetilde g$ and denote by $\nabla^t$ the Riemannian connection of $g_t$.  
Fix a torsion-free connection $\mathring{\nabla}$ on $M$.  Define a $(1,2)$-tensor field $S_t$ on $M$ by $S_t:=\nabla^t-\mathring{\nabla}$ and 
a vector field $V(f_t)$ on $M$ by $V(f_t):={\rm Tr}_{g_t}S_t$, where ${\rm Tr}_{g_t}S_t$ is the trace of $S_t$ with respect to $g_t$.  
The mean curvature-de Turck equation is defined by 
$$\frac{\partial f_t}{\partial t}=H_t+df_t(V(f_t)).\leqno{(3.1)}$$
We give the local expressions of the mean curvature flow equation and the mean curvature-de Turck equation.  
Let $n:={\rm dim}\,M$ and $n+r:={\rm dim}\,N$.  
Take a local coordinate $(U,(x_1,\cdots,x_n))$ of $M$ and a local coordinate $(W,(y_1,$\newline
$\cdots,y_{n+r}))$ of $N$ with $f(U)\subset W$.  
The local expression of the mean curvature flow equation is given by 
$$\begin{array}{l}
\displaystyle{\frac{\partial(f_t)^{\gamma}}{\partial t}=\sum_{i,j=1}^n(g_t)^{ij}\left(\partial_i\partial_j(f_t)^{\gamma}
+\sum_{\alpha,\beta=1}^{n+r}\partial_i(f_t)^{\alpha}\partial_j(f_t)^{\beta}
(\widetilde{\Gamma}_{\alpha\beta}^{\gamma}\circ f_t)\right.}\\
\hspace{3.25truecm}\displaystyle{\left.-\sum_{k=1}^n\partial_k(f_t)^{\gamma}(\Gamma_t)_{ij}^k\right),}
\end{array}\leqno{(3.2)}$$
where $(f_t)^{\gamma}$'s are the components of $f_t$ with respect to $(U,(x_1,\cdots,x_n))$ and $(W,(y_1,\cdots,y_{n+r}))$, 
$((g_t)^{ij})$ is the inverse matrix of the matrix $((g_t)_{ij})$ consisting of the components $(g_t)_{ij}$'s of the induced metric $g_t$ 
with respect to $(U,(x_1,\cdots,x_n))$, $\widetilde{\Gamma}_{\alpha\beta}^{\gamma}$'s are 
the Christoffel's symbols of $\widetilde g$ with respect to $(W,(y_1,\cdots,y_{n+r}))$ and 
$(\Gamma_t)_{ij}^k$ is the Christoffel's symbols of $g_t$ with repspec to $(U,(x_1,\cdots,x_n))$.  Here we note that $(g_t)_{ij}$ is given by 
$$(g_t)_{ij}=\sum_{\alpha,\beta=1}^{n+r}(\partial_i(f_t)^{\alpha})(\partial_j(f_t)^{\beta})(\widetilde g_{\alpha\beta}\circ f_t)$$
and $(\Gamma_t)_{ij}^k$ is given by 
$$(\Gamma_t)_{ij}^k=\sum_{l=1}^n\frac{(g_t)^{kl}}{2}\left(\partial_i(g_t)_{lj}+\partial_j(g_t)_{il}-\partial_l(g_t)_{ij}\right).$$
By the existence of the final term of the right-hand side of $(3.2)$, the equation $(3.2)$ is a quasi-linear parabolic partial differential equation 
but not strongly parabolic.  
On the other hand, the local expression of the mean curvature-de Turck equation is given by 
$$\begin{array}{l}
\displaystyle{\frac{\partial(f_t)^{\gamma}}{\partial t}=
\sum_{i,j=1}^n(g_t)^{ij}\left(
\partial_i\partial_j(f_t)^{\gamma}
+\sum_{\alpha,\beta=1}^{n+r}\partial_i(f_t)^{\alpha}\partial_j(f_t)^{\beta}
(\widetilde{\Gamma}_{\alpha\beta}^{\gamma}\circ f_t)\right.}\\
\hspace{3.25truecm}\displaystyle{\left.-\sum_{k=1}^n\partial_k(f_t)^{\gamma}\mathring{{\Gamma}}_{ij}^k\right),}
\end{array}\leqno{(3.3)}$$
where $(f_t)^{\gamma}$'s, $((g_t)^{ij})$ and $\widetilde{\Gamma}_{\alpha\beta}^{\gamma}$ are as above, and 
${\mathring{\Gamma}}_{ij}^k$'s are the Christoffel's symbol of a fixed connection $\mathring{\nabla}$ of $M$.  
Since the final term of the right-hand side of $(3.2)$ is changed by 
$\sum\limits_{i,j,k=1}^n(g_t)^{ij}\partial_k(f_t)^{\gamma}{\mathring{\Gamma}}_{ij}^k$, the equation $(3.3)$ (hence (3.1)) is 
a strongly parabolic quasi-linear parabolic partial differential equation.  Hence, if $M$ is a compact manifold (without boundary), then 
it is shown that there uniquely exists a solution $\{f_t\}_{t\in[0,T)}$ of $(3.1)$ with $f_0=f$ in short time.  
As B. Kleiner and J. Lott state in Subsection 5.2 of \cite{KL} in the case of the Ricci flow on a compact orbifold, 
even if $M$ is a compact orbifold (not manifold), it is shown similarly that there uniquely exists a solution $\{f_t\}_{t\in[0,T)}$ of $(3.1)$ 
with $f_0=f$ in short time.  For the solution $\{f_t\}_{t\in[0,T)}$, we consider the ordinary differential equation 
$\displaystyle{\frac{\partial\psi_t}{\partial t}=-V(f_t)\circ\psi_t}$.  Let $\{\psi_t\}_{t\in[0,T')}$ ($T'<T$) be the solution of this ordinary 
differential equation with the initial condition $\psi_0={\rm id}_M$, where we note that each $\psi_t$ is a $C^{\infty}$-diffeomorphism of $M$ onto 
oneself.  Then it is shown that $\{f_t\circ\psi_t\}_{t\in[0,T')}$ is a mean curvature flow starting from $f$.  
Conversely, it is easy to show that any mean curvature flow starting from $f$ is given like this.  This completes the proof.  \qed

\vspace{0.5truecm}

\noindent
{\it Remark 3.1.}\ Let $f$ be a $C^{\infty}$-immersion of a manifold $M$ into a complete Riemannian manifold $(N,\widetilde g)$ and 
$D$ a relative compact domain of $M$.  The existence of a mean curvature flow starting from $f|_D$ is shown but its uniqueness is not shown.  
In fact, it is shown that there exist infinitely many mean curvature flows starting from $f|_D$ as follows.  
Let $\widetilde{f|_D}:M\hookrightarrow N$ be an immersion satisfying the following conditions:

\vspace{0.15truecm}

(i)\ $(\widetilde{f|_D})|_D=f|_D$;

(ii)\ $\widetilde{f|_D}(M)$ is an immersed complete Riemannian submanifold of bounded second fundamental form in $(N,\widetilde g)$, where 
we note that the condition of\newline
``bounded second fundamental form'' controls the behavior of the submanifold $\widetilde{f|_D}(M)$ near the infinity.  

\vspace{0.15truecm}

\noindent
It is clear that there exists infinitely many complete exrensions $f|_D$ satisfying the condition (ii).  
Assume that the norms of the curvature tensor of $(N,\widetilde g)$, its first derivative and its second derivative are bounded.  
Then there uniquely exists a mean curvature flow $\{(\widetilde{f|_D})_t\}_{t\in[0,T)}$ of bounded second fundamental form starting from 
$\widetilde{f|_D}$ by the uniqueness theorem in \cite{CY}.  It is clear that 
$\{(\widetilde{f|_D})_t|_D\}_{t\in[0,T)}$ gives a mean curvature flow starting from $f|_D$ and that this flow depends on the choice of 
the complete extension $\widetilde{f|_D}$.  
Thus we see that there exist infinitely many mean curvature flows starting from $f|_D$.  

\section{Evolution equations}
Let $\mathcal G\curvearrowright V$ be an isometric almost free action with minimal regularizable orbit of a Hilbert Lie group $\mathcal G$ 
on a Hilbert space $V$ equipped with an inner product $\langle\,\,,\,\,\rangle$.  
The orbit space $V/\mathcal G$ is a (finite dimensional) $C^{\infty}$-orbifold.  
Let $\phi:V\to V/\mathcal G$ be the orbit map and set $N:=V/\mathcal G$.  
Here we give an example of such an isometric almost free action of a Hilbert Lie group.  

\vspace{0.5truecm}

\noindent
{\it Example 4.1.}\ Let $G$ be a compact semi-simple Lie group, $K$ a closed subgroup of $G$ and $\Gamma$ a discrete subgroup of $G$.  
Denote by $\mathfrak g$ and $\mathfrak k$ the Lie algebras of $G$ and $K$, 
respectively.  Assume that a reductive decomposition 
$\mathfrak g=\mathfrak k+\mathfrak p$ exists.  
Let $B$ be the Killing form of $\mathfrak g$ and $g$ the bi-invariant metric of $G$ induced from $-B$.  
Also, let $H^0([0,1],\mathfrak g)$ be the Hilbert space of all paths in the Lie algebra $\mathfrak g$ of $G$ which are $L^2$-integrable with respect to 
$-B$, and $H^1([0,1],G)$ the Hilbert Lie group of all paths in $G$ which are of class $H^1$ with respect to $g$.  This group $H^1([0,1],G)$ acts on 
$H^0([0,1],\mathfrak g)$ isometrically and transitively as a gauge action:
$$\begin{array}{r}
\displaystyle{({\bf g}\cdot u)(s)={\rm Ad}_G({\bf g}(s))(u(s))-(R_{{\bf g}(s)})_{\ast}^{-1}({\bf g}'(s))
\,\,\,\,(s\in[0,1])}\\
\displaystyle{({\bf g}\in H^1([0,1],G),\,u\in H^0([0,1],\mathfrak g)).}
\end{array}$$
Set $P(G,\Gamma\times K):=\{{\bf g}\in H^1([0,1],G)\,\vert\,({\bf g}(0),{\bf g}(1))\in\Gamma\times K\}$.  
The group $P(G,\Gamma\times K)$ acts on $H^0([0,1],\mathfrak g)$ almost freely and isometrically, and the orbit space of this action is diffeomorphic 
to the orbifold $\Gamma\setminus G\,/\,K$.  
Furthermore, each orbit of this action is regularizable and minimal (see \cite{HLO}, \cite{PiTh}, \cite{Te1}, \cite{Te2}, \cite{TeTh}, 
\cite{Koi2}).  In particular, in the case of $K=\Gamma=\{e\}$, $\phi$ is a Riemannian submersion of 
$H^0([0,1].\mathfrak g)$ onto $(G,g)$ and is called the {\it parallel transport map for} $G$.  

\vspace{0.5truecm}

Let $g_N$ be the Riemannian orbimetric on $N$ such that $\phi$ is a Riemannian orbisubmersion of $(V,\langle\,\,,\,\,\rangle)$ onto $(N,g_N)$.  
By using Theorem 3.1, we prove the following unique existence theorem for a $\mathcal G$-invariant regularized mean curvature flow starting from 
a $\mathcal G$-invariant regularizable hypersurface with compact quotient.  

\vspace{0.5truecm}

\noindent
{\bf Theorem 4.1.}\ {\sl Let $f:M\hookrightarrow V$ be an immersion such that $f(M)$ is a $\mathcal G$-invariant regularizable hypesurface and that 
$(\phi\circ f)(M)$ is compact.  Then there uniquely exists a $\mathcal G$-invariant regularized mean curvature flow $\{f_t\}_{t\in[0,T)}$ 
starting from $f$.}  

\vspace{0.5truecm}

\noindent
{\it Proof.}\ 
Since $f(M)$ is $\mathcal G$-invariant and $\phi(f(M))$ is compact, we can take an orbiimmesion $\overline f$ of a compact orbifold $M'$ into $N$ and 
an orbifold submersion $\phi_M:M\to M'$ with $\phi\circ f=\overline f\circ\phi_M$.  
Let $H$ be the regularized mean curvature vector of $f$ and $\overline H$ the mean curvature vector of $\overline f$.  
Then, since the fibres of $\phi$ are minimal regularizable submanifolds, 
we have $\phi_{\ast}(H)=\overline H$.  In more detail, $H$ is the horizontal lift of 
$\overline H$ (with respect to the Riemannian orbisubmersion $\phi$).  
According to Theorem 3.1, there uniquely exists the mean curvature flow $\{\overline f_t\}_{t\in[0,T)}$ starting from $\overline f$.  
Define a curve $c_x:[0,T)\to N$ by $c_x(t):=\overline f_t(x)$ and let $(c_x)_u^L:[0,T)\to V$ be the horizontal lift of $c_x$ to $f(u)$.  
Define an immersion $f_t:M\hookrightarrow V$ by $f_t(u)=(c_x)_u^L(t)$ ($u\in\widetilde M$).  
Let $H_t$ be the regularized mean curvature vector of $f_t$ and $\overline H_t$ the mean curvature of $\overline f_t$.  
Then, it is clear that $f_t(M)=\phi^{-1}(\overline f_t(M'))$ holds.  Hence, since the fibres of $\phi$ are minimal regularizable submanifolds, 
$H_t$ is the horizontal lift of $\overline H_t$.  On the other hand, it follows from the construction of $f_t$ that 
$\displaystyle{\frac{\partial f_t}{\partial t}}$ is the horizontal lift of $\displaystyle{\frac{\partial\overline f_t}{\partial t}}$.  
Also, since $\{\overline f_t\}_{t\in [0,T)}$ is the mean curvature flow, 
we have $\displaystyle{\frac{\partial\overline f_t}{\partial t}=\overline H_t}$.  
From these facts, we can derive $\displaystyle{\frac{\partial f_t}{\partial t}=H_t}$.  
Hence $\{f_t\}_{t\in[0,T)}$ is a regularized mean curvature flow starting from $f$.  
It is clear that $f_t(M)$'s ($0\leq t<T$) are $\mathcal G$-invariant.  
Thus the existence of a $\mathcal G$-invariant regularized mean curvature flow starting from $f$ is shown.  

Next we shall show the uniqueness of a $\mathcal G$-invariant regularized mean curvature flow $\{f_t\}_{t\in[0,T)}$ starting from $f$.  
Let $\{f_t\}_{t\in[0,T)}$ be such a flow.  Then, since $f_t(M)$ is $\mathcal G$-invariant, 
we can take an orbiimmersion $\overline f_t$ of a compact orbifold $M'$ into $N$ satisfying $\phi\circ f_t=\overline f_t\circ\phi_M$.  
Then we have 
$$\frac{\partial\overline f_t}{\partial t}=\phi_{\ast}\left(\frac{\partial f_t}{\partial t}\right)
=\phi_{\ast}(H_t)=\overline H_t.$$
Hence $\{\overline f_t\}_{t\in[0,T)}$ is the mean curvature flow starting from $\overline f$.  
Therefore $\{f_t\}_{t\in[0,T)}$ is the flow constructed as above.  
Thus the uniqueness of a $\mathcal G$-invariant regularized mean curvature flow starting from $f$ also is shown.  \qed

\vspace{0.5truecm}

\noindent
{\it Remark 4.1.}\ We cannot conclude whether there uniquely exists a (not necessarily $\mathcal G$-invariant) regularized mean curvature flow 
starting from a $\mathcal G$-invariant reguralized immersion $f$ as in the statement of Theorem 4.1.  

\vspace{0.5truecm}

Let $f:M\hookrightarrow V$ be an immersion such that $f(M)$ is a $\mathcal G$-invariant regularizable hypesurface and that 
$(\phi\circ f)(M)$ is compact and $\{f_t\}_{t\in[0,T)}$ the $\mathcal G$-invariant regularized mean curvature flow starting from $f$.  
Define a map $F:M\times[0,T)\to V$ by 
$F(u,t):=f_t(u)$ ($(u,t)\in M\times[0,T)$) and 
a map $\overline F:\overline M\times[0,T)\to N$ by 
$\overline F(x,t):=\overline f_t(x)$ ($(x,t)\in \overline M\times[0,T)$).  
Denote by $H_t$ the regularized mean curvature vector of $f_t$ and $\overline H_t$ the 
mean curvature vector of $\overline f_t$.  Since $\phi$ has minimal regularizable fibres, 
$H_t$ is the horizontal lift of $\overline H_t$, we can show that $\phi\circ f_t=\overline f_t\circ\phi_M$ holds 
for all $t\in [0,T)$.  In the sequel, we consider the case where the codimension of $M$ is equal to one.  
Denote by $\widetilde{{\mathcal H}}$ (resp. $\widetilde{{\mathcal V}}$) the horizontal 
(resp. vertical) distribution of $\phi$.  
Denote by ${\rm pr}_{\widetilde{\mathcal H}}$ (resp. 
${\rm pr}_{\widetilde{\mathcal V}}$) the orthogonal projection of $TV$ onto 
$\widetilde{\mathcal H}$ (resp. $\widetilde{\mathcal V}$).  For simplicity, for 
$X\in TV$, we denote ${\rm pr}_{\widetilde{\mathcal H}}(X)$ (resp. 
${\rm pr}_{\widetilde{\mathcal V}}(X)$) by $X_{\widetilde{\mathcal H}}$ 
(resp. $X_{\widetilde{\mathcal V}}$).  
Define a distribution ${\mathcal H}_t$ on 
$M$ by $f_{t\ast}(({\mathcal H}_t)_u)=f_{t\ast}(T_uM)\cap
\widetilde{{\mathcal H}}_{f_t(u)}$ ($u\in M$) and a distribution ${\mathcal V}_t$ on 
$M$ by $f_{t\ast}(({\mathcal V}_t)_u)=\widetilde{{\mathcal V}}_{f_t(u)}$ ($u\in M$).  
Note that ${\mathcal V}_t$ is independent of the choice of $t\in[0,T)$.  
Fix a unit normal vector field $\xi_t$ of $f_t$.  Denote by $g_t,h_t,A_t,H_t$ and $H^s_t$ 
the induced metric, the second fundamental form (for $-\xi_t$), the shape operator (for $-\xi_t$) and 
the regularized mean curvature vector and the regularized mean curvature (for $-\xi_t$), respectively.  
The group $\mathcal G$ acts on $M$ through $f_t$.  
Since $\phi:V\to V/\mathcal G$ is a $\mathcal G$-orbibundle and $\widetilde{\mathcal H}$ is a connection 
of this orbibundle, it follows from Proposition 4.1 in \cite{Koi2} that this action $\mathcal G\curvearrowright M$ 
is independent of the choice of $t\in[0,T)$.  
It is clear that quantities $g_t,h_t,A_t,H_t$ and $H^s_t$ are $\mathcal G$-invariant.  
Also, let $\nabla^t$ 
be the Riemannian connection of $g_t$.  
Let $\pi_M$ be the projection of $M\times[0,T)$ onto $M$.  
For a vector bundle $E$ over $M$, denote by $\pi_M^{\ast}E$ the induced bundle 
of $E$ by $\pi_M$.  Also denote by $\Gamma(E)$ the space of all sections of 
$E$.  Define a section $g$ of $\pi_M^{\ast}(T^{(0,2)}M)$ by 
$g(u,t)=(g_t)_u$ ($(u,t)\in M\times[0,T)$), where $T^{(0,2)}M$ is the 
$(0,2)$-tensor bundle of $M$.  Similarly, we define a section 
$h$ of $\pi_M^{\ast}(T^{(0,2)}M)$, a section $A$ of $\pi_M^{\ast}
(T^{(1,1)}M)$, a section $H$ of $F^{\ast}TV$ and a section $\xi$ of $F^{\ast}TV$.  
We regard $H$ and $\xi$ as $V$-valued functions over 
$M\times[0,T)$ under the identification of $T_uV$'s ($u\in V$) and $V$.  
Define a subbundle ${\mathcal H}$ (resp. ${\mathcal V}$) of 
$\pi_M^{\ast}TM$ by ${\mathcal H}_{(u,t)}:=({\mathcal H}_t)_u$ (resp. 
${\mathcal V}_{(u,t)}:=({\mathcal V}_t)_u$).  Denote by ${\rm pr}_{\mathcal H}$ (resp. 
${\rm pr}_{\mathcal V}$) the orthogonal projection of $\pi_M^{\ast}(TM)$ onto 
${\mathcal H}$ (resp. ${\mathcal V}$).  For simplicity, for $X\in \pi_M^{\ast}(TM)$, 
we denote ${\rm pr}_{\mathcal H}(X)$ (resp. ${\rm pr}_{\mathcal V}(X)$) by 
$X_{\mathcal H}$ (resp. $X_{\mathcal V}$).  
For a section $B$ of $\pi_M^{\ast}(T^{(r,s)}M)$, 
we define $\displaystyle{\frac{\partial B}{\partial t}}$ by 
$\displaystyle{\left(\frac{\partial B}{\partial t}\right)_{(u,t)}
:=\frac{d B_{(u,t)}}{dt}}$, where the right-hand side of this relation is 
the derivative of the vector-valued function $t\mapsto B_{(u,t)}\,(\in 
T^{(r,s)}_uM)$.  Also, we define a section $B_{\mathcal H}$ of 
$\pi_M^{\ast}(T^{(r,s)}M)$ by 
$$B_{\mathcal H}=\mathop{({\rm pr}_{\mathcal H}\otimes\cdots
\otimes{\rm pr}_{\mathcal H})}_{(r{\rm -times})}\circ B\circ
\mathop{({\rm pr}_{\mathcal H}\otimes\cdots
\otimes{\rm pr}_{\mathcal H})}_{(s{\rm -times})}.$$
The restriction of $B_{\mathcal H}$ to 
${\mathcal H}\times\cdots\times{\mathcal H}$ ($s$-times) is regarded as a section of 
the $(r,s)$-tensor bundle ${\mathcal H}^{(r,s)}$ of $\mathcal H$.  This restriction 
also is denoted by the same symbol $B_{\mathcal H}$.  
Let $D_M$ (resp. $D_{[0,T)}$) be the subbundle of $T(M\times[0,T))$ defined by 
$(D_M)_{(u,t)}:=T_{(u,t)}(M\times\{t\})$ (resp. $(D_{[0,T)})_{(u,t)}:=T_{(x,t)}(\{u\}\times[0,T))$ 
for each $(u,t)\in M\times[0,T)$.  
Denote by $v^L_{(u,t)}$ the horozontal lift of $v\in T_uM$ to $(u,t)$ with respect to $\pi_M$ (i.e., 
$v^L_{(u,t)}$ is the element of $(D_M)_{(u,t)}$ with $(\pi_M)_{\ast(u,t)}(v^L_{(u,t)})=v$).  
Under the identification of $((u,t),v)(=v)\in(\pi^{\ast}TM)_{(u,t)}$ with 
$v^L_{(u,t)}\in(D_M)_{(u,t)}$, we identify $\pi_M^{\ast}TM$ with $D_M$.  
For a tangent vector field $X$ on $M$ (or an open set $U$ of $M$), we define 
$\overline X\in\Gamma(\pi_M^{\ast}TM)(=\Gamma(D_M))\,\,{\rm(or}\,\,\Gamma((\pi_M^{\ast}TM)\vert_U)
(=\Gamma((D_M)\vert_U)){\rm)}$ by 
$\overline X_{(u,t)}:=((u,t),X_u)(=(X_u)_{(u,t)}^L)$ ($(u,t)\in M\times[0,T)$).  
Denote by $\widetilde{\nabla}$ the Riemannian connection of $V$.  
Let $\nabla$ be the connection of $\pi_M^{\ast}TM$ defined by 
$$(\nabla_XY)_{(u,t)}:=\nabla^t_{X_{(u,t)}}Y_{(\cdot,t)}\,\,\,{\rm and}\,\,\,
(\nabla_{\frac{\partial}{\partial t}}Y)_{(u,t)}:=\frac{dY_{(u,\cdot)}}{dt}$$
for $X\in\Gamma(D_M)$ and $Y\in\Gamma(\pi_M^{\ast}TM)$, where 
$X_{(u,t)}$ is identified with $(\pi_M)_{\ast(u,t)}$\newline
$(X_{(u,t)})\in T_uM$, $Y_{(\cdot,t)}$ is identified with 
$(\pi_M)_{\ast}(Y_{(\cdot,t)})\in\Gamma(TM)$ and $Y_{(u,\cdot)}$ is identified with 
$(\pi_M)_{\ast}(Y_{u,\cdot})\in C^{\infty}([0,T),T_uM)$.  
Note that $\nabla_{\frac{\partial}{\partial t}}\overline X=0$.  
Denote by the same symbol $\nabla$ the connection of $\pi_M^{\ast}T^{(r,s)}M$ defined in terms of 
$\nabla^t$'s similarly.  
Define a connection $\nabla^{\mathcal H}$ of ${\mathcal H}$ by 
$\nabla^{\mathcal H}_XY:=(\nabla_XY)_{\mathcal H}$ for $X\in\Gamma(M\times[0,T))$ and 
$Y\in\Gamma({\mathcal H})$.  Similarly, define a connection 
$\nabla^{\mathcal V}$ of ${\mathcal V}$ by $\nabla^{\mathcal V}_XY:=(\nabla_XY)_{\mathcal V}$ 
for $X\in\Gamma(M\times[0,T))$ and $Y\in\Gamma({\mathcal V})$.  

Now we shall recall the evolution equations for some geometric quantities given in \cite{Koi2}.  
By the same calculation as the proof of Lemma 4.2 of \cite{Koi2} 
(where we replace $H=\|H\|\xi$ in the proof to $H=\langle H,\xi\rangle\xi=-H^s\xi$), we can derive the following 
evolution equation.  

\vspace{0.5truecm}

\noindent
{\bf Lemma 4.2.} {\sl The sections $(g_{\mathcal H})_t$'s of $\pi_M^{\ast}(T^{(0,2)}M)$ 
satisfy the following evolution equation:
$$\frac{\partial g_{\mathcal H}}{\partial t}=-2H^s h_{\mathcal H}.$$
}

\vspace{0.5truecm}

According to the proof of Lemma 4.3 in \cite{Koi2}, we obtain the following evolution equation.  

\vspace{0.5truecm}

\noindent
{\bf Lemma 4.3.} {\sl The unit normal vector fields $\xi_t$'s 
satisfy the following evolution equation:
$$\frac{\partial\xi}{\partial t}=-F_{\ast}({\rm grad}_gH^s),$$
where ${\rm grad}_g(H^s)$ is the element of $\pi_M^{\ast}(TM)$ such that $dH^s(X)=g({\rm grad}_gH^s,X)$ for 
any $X\in\pi_M^{\ast}(TM)$.}

\vspace{0.5truecm}

Let $S_t$ ($0\leq t<T$) be a $C^{\infty}$-family of a $(r,s)$-tensor fields on 
$M$ and $S$ a section of $\pi_M^{\ast}(T^{(r,s)}M)$ defined by 
$S_{(u,t)}:=(S_t)_u$.  We define a section $\triangle_{{\mathcal H}}S$ 
of $\pi_M^{\ast}(T^{(r,s)}M)$ by 
$$(\triangle_{{\mathcal H}}S)_{(u,t)}:=\sum_{i=1}^n\nabla_{e_i}\nabla_{e_i}S,$$
where $\nabla$ is the connection of $\pi_M^{\ast}(T^{(r,s)}M)$ (or 
$\pi_M^{\ast}(T^{(r,s+1)}M)$) induced from $\nabla$ and $(e_1,\cdots,e_n)$ 
is an orthonormal base of ${\mathcal H}_{(u,t)}$ with respect to 
$(g_{\mathcal H})_{(u,t)}$.  
Also, we define a section $\bar{\triangle}_{{\mathcal H}}S_{\mathcal H}$ 
of ${\mathcal H}^{(r,s)}$ by 
$$(\triangle_{\mathcal H}^{\mathcal H}S_{\mathcal H})_{(u,t)}:=
\sum_{i=1}^n\nabla^{\mathcal H}_{e_i}\nabla^{\mathcal H}_{e_i}S_{\mathcal H},$$
where $\nabla^{\mathcal H}$ is the connection of ${\mathcal H}^{(r,s)}$ (or 
${\mathcal H}^{(r,s+1)}$) induced from $\nabla^{\mathcal H}$ and $\{e_1,\cdots,e_n\}$ 
is as above.  Let ${\mathcal A}^{\phi}$ be the section of $T^{\ast}V\otimes 
T^{\ast}V\otimes TV$ defined by 
$${\mathcal A}^{\phi}_XY:=(\widetilde{\nabla}_{X_{\widetilde{\mathcal H}}}
Y_{\widetilde{\mathcal H}})_{\widetilde{\mathcal V}}
+(\widetilde{\nabla}_{X_{\widetilde{\mathcal H}}}Y_{\widetilde{\mathcal V}})_{\widetilde{\mathcal H}}
\,\,\,\,(X,Y\in TV).$$
Also, let ${\mathcal T}^{\phi}$ be the section of 
$T^{\ast}V\otimes T^{\ast}V\otimes TV$ defined by 
$${\mathcal T}^{\phi}_XY:=(\widetilde{\nabla}_{X_{\widetilde{\mathcal V}}}
Y_{\widetilde{\mathcal H}})_{\widetilde{\mathcal V}}
+(\widetilde{\nabla}_{X_{\widetilde{\mathcal V}}}
Y_{\widetilde{\mathcal V}})_{\widetilde{\mathcal H}}
\,\,\,\,(X,Y\in TV).$$
Also, let ${\mathcal A}_t$ be the section of $T^{\ast}M\otimes 
T^{\ast}M\otimes TM$ defined by 
$$({\mathcal A}_t)_XY:=(\nabla^t_{X_{{\mathcal H}_t}}
Y_{{\mathcal H}_t})_{{\mathcal V}_t}+(\nabla^t_{X_{{\mathcal H}_t}}
Y_{{\mathcal V}_t})_{{\mathcal H}_t}\,\,\,\,(X,Y\in TM).$$
Also let ${\mathcal A}$ be the section of $\pi_M^{\ast}(T^{\ast}M\otimes T^{\ast}M
\otimes TM)$ defined in terms of ${\mathcal A}_t$'s ($t\in[0,T)$).  
Also, let ${\mathcal T}_t$ be the section of $T^{\ast}M\otimes 
T^{\ast}M\otimes TM$ defined by 
$$({\mathcal T}_t)_XY:=(\nabla^t_{X_{{\mathcal V}_t}}Y_{{\mathcal V}_t})_{{\mathcal H}_t}
+(\nabla^t_{X_{{\mathcal V}_t}}Y_{{\mathcal H}_t})_{{\mathcal V}_t}\,\,\,\,(X,Y\in TM).$$
Also let ${\mathcal T}$ be the section of $\pi_M^{\ast}(T^{\ast}M\otimes T^{\ast}M
\otimes TM)$ defined in terms of ${\mathcal T}_t$'s ($t\in[0,T)$).  
Clearly we have 
$$F_{\ast}({\mathcal A}_XY)={\mathcal A}^{\phi}_{F_{\ast}X}F_{\ast}Y$$
for $X,Y\in{\mathcal H}$ and 
$$F_{\ast}({\mathcal T}_WX)={\mathcal T}^{\phi}_{F_{\ast}W}F_{\ast}X$$
for $X\in{\mathcal H}$ and $W\in{\mathcal V}$.  
Let $E$ be a vector bundle over $M$.  For a section $S$ of 
$\pi_M^{\ast}(T^{(0,r)}M\otimes E)$, we define 
$\displaystyle{{\rm Tr}_{g_{\mathcal H}}^{\bullet}\,S(\cdots,\mathop{\bullet}^j,
\cdots,\mathop{\bullet}^k,\cdots)}$ by 
$$({\rm Tr}_{g_{\mathcal H}}^{\bullet}\,S(\cdots,\mathop{\bullet}^j,\cdots,
\mathop{\bullet}^k,\cdots))
_{(u,t)}=\sum_{i=1}^nS_{(u,t)}(\cdots,\mathop{e_i}^j,\cdots,\mathop{e_i}^k,
\cdots)$$
$((u,t)\in M\times[0,T))$, where $\{e_1,\cdots,e_n\}$ is an 
orthonormal base of ${\mathcal H}_{(u,t)}$ with respect to 
$(g_{\mathcal H})_{(u,t)}$, $\displaystyle{S(\cdots,\mathop{\bullet}^j,
\cdots,\mathop{\bullet}^k,\cdots)}$ means that $\bullet$ is entried into the 
$j$-th component and the $k$-th component of $S$ and 
$\displaystyle{S_{(u,t)}(\cdots,\mathop{e_i}^j,\cdots,\mathop{e_i}^k,\cdots)}$ 
means that $e_i$ is entried into the $j$-th component and the $k$-th component 
of $S_{(u,t)}$.  

In \cite{Koi2}, we derived the following relation.  

\vspace{0.5truecm}

\noindent
{\bf Lemma 4.4.} {\sl Let $S$ be a section of 
$\pi_M^{\ast}(T^{(0,2)}M)$ which is symmetric with respect to $g$.  
Then we have 
$$\begin{array}{l}
\displaystyle{(\triangle_{\mathcal H}S)_{\mathcal H}(X,Y)=
(\triangle_{\mathcal H}^{\mathcal H}S_{\mathcal H})(X,Y)}\\
\hspace{3.2truecm}
\displaystyle{-2{\rm Tr}^{\bullet}_{g_{\mathcal H}}
((\nabla_{\bullet}S)({\mathcal A}_{\bullet}X,Y))
-2{\rm Tr}^{\bullet}_{g_{\mathcal H}}((\nabla_{\bullet}S)({\mathcal A}_{\bullet}Y,X))
}\\
\hspace{3.2truecm}
\displaystyle{-{\rm Tr}^{\bullet}_{g_{\mathcal H}}
S({\mathcal A}_{\bullet}({\mathcal A}_{\bullet}X),Y)
-{\rm Tr}^{\bullet}_{g_{\mathcal H}}S({\mathcal A}_{\bullet}({\mathcal A}_{\bullet}Y),X)}\\
\hspace{3.2truecm}
\displaystyle{-{\rm Tr}^{\bullet}_{g_{\mathcal H}}
S((\nabla_{\bullet}{\mathcal A})_{\bullet}X,Y)
-{\rm Tr}^{\bullet}_{g_{\mathcal H}}S((\nabla_{\bullet}{\mathcal A})_{\bullet}Y,X)}\\
\hspace{3.2truecm}
\displaystyle{-2{\rm Tr}^{\bullet}_{g_{\mathcal H}}
S({\mathcal A}_{\bullet}X,{\mathcal A}_{\bullet}Y)}
\end{array}$$
for $X,Y\in{\mathcal H}$, where $\nabla$ is the connection of 
$\pi_M^{\ast}(T^{(1,2)}M)$ induced from $\nabla$.}

\vspace{0.5truecm}

According to the proof of Lemma 4.5 in \cite{Koi2}, we obtain the following Simons-type identity.  

\vspace{0.5truecm}

\noindent
{\bf Lemma 4.5.} {\sl We have 
$$\triangle_{{\mathcal H}}h=\nabla dH^s+H^s(A^2)_{\sharp}-({\rm Tr}\,(A^2)_{{\mathcal H}})h,$$
where $(A^2)_{\sharp}$ is the element of $\Gamma(\pi_M^{\ast}T^{(0,2)}M)$ 
defined by $(A^2)_{\sharp}(X,Y):=g(A^2X,Y)$ ($X,Y\in\pi_M^{\ast}TM$).}

\vspace{0.5truecm}

\noindent
{\bf Note.} In the sequel, we omit the notation $F_{\ast}$ 
for simplicity.  

\vspace{0.5truecm}

Define a section ${\mathcal R}$ of $\pi_M^{\ast}({\mathcal H}^{(0,2)})$ by 
$$\begin{array}{l}
\displaystyle{{\mathcal R}(X,Y):=
{\rm Tr}^{\bullet}_{g_{\mathcal H}}h({\mathcal A}_{\bullet}({\mathcal A}_{\bullet}X),Y)
+{\rm Tr}^{\bullet}_{g_{\mathcal H}}h({\mathcal A}_{\bullet}({\mathcal A}_{\bullet}Y),X)}\\
\hspace{2.25truecm}\displaystyle{
+{\rm Tr}^{\bullet}_{g_{\mathcal H}}h((\nabla_{\bullet}{\mathcal A})_{\bullet}X,
Y)+{\rm Tr}^{\bullet}_{g_{\mathcal H}}h((\nabla_{\bullet}{\mathcal A})_{\bullet}Y,X)}\\
\hspace{2.25truecm}\displaystyle{
+2{\rm Tr}^{\bullet}_{g_{\mathcal H}}(\nabla_{\bullet}h)({\mathcal A}_{\bullet}X,Y)
+2{\rm Tr}^{\bullet}_{g_{\mathcal H}}(\nabla_{\bullet}h)({\mathcal A}_{\bullet}Y,X)}\\
\hspace{2.25truecm}\displaystyle{
+2{\rm Tr}^{\bullet}_{g_{\mathcal H}}h({\mathcal A}_{\bullet}X,{\mathcal A}_{\bullet}Y)
\qquad\qquad\quad(X,Y\in{\mathcal H}).}
\end{array}$$
According to Theorem 4.6 in \cite{Koi2}, we obtain the following evolution equation from from 
Lemmas 4.3, 4.4 and 4.5.  

\vspace{0.5truecm}

\noindent
{\bf Lemma 4.6.} {\sl The sections $(h_{\mathcal H})_t$'s of $\pi_M^{\ast}(T^{(0,2)}M)$ 
satisfies the following evolution equation:
$$\begin{array}{l}
\displaystyle{\frac{\partial h_{\mathcal H}}{\partial t}(X,Y)
=(\triangle_{\mathcal H}^{\mathcal H}h_{\mathcal H})(X,Y)-2H^s
((A_{\mathcal H})^2)_{\sharp}(X,Y)-2H^s
(({\mathcal A}^{\phi}_{\xi})^2)_{\sharp}(X,Y)}\\
\hspace{2.4truecm}\displaystyle{
+{\rm Tr}\left((A_{\mathcal H})^2-(({\mathcal A}^{\phi}_{\xi})^2)_{\mathcal H}\right)
h_{\mathcal H}(X,Y)
-{\mathcal R}(X,Y)}
\end{array}$$
for $X,Y\in{\mathcal H}$.}

\vspace{0.5truecm}

According to the proof of Lemma 4.8 of \cite{Koi2}, we obatin the following relation from Lemma 4.2.  

\vspace{0.5truecm}

\noindent
{\bf Lemma 4.7.} {\sl Let $X$ and $Y$ be local sections 
of ${\mathcal H}$ such that $g(X,Y)$ is constant.  
Then we have 
$g(\nabla_{\frac{\partial}{\partial t}}X,Y)
+g(X,\nabla_{\frac{\partial}{\partial t}}Y)
=2H^s h(X,Y)$.}

\vspace{0.5truecm}

According to Lemmas 4.8 and 4.10 in \cite{Koi2}, we obtain the following relation .  

\vspace{0.5truecm}

\noindent
{\bf Lemma 4.8.} {\sl For $X,Y\in{\mathcal H}$, we have 
$$\begin{array}{l}
\displaystyle{{\mathcal R}(X,Y)=2{\rm Tr}^{\bullet}_{g_{\mathcal H}}
\left(\langle({\mathcal A}^{\phi}_{\bullet}X,{\mathcal A}^{\phi}_{\bullet}
(A_{\mathcal H}Y)\rangle
+\langle({\mathcal A}^{\phi}_{\bullet}Y,{\mathcal A}^{\phi}_{\bullet}(A_{\mathcal H}X)
\rangle\right)}\\
\hspace{1.95truecm}\displaystyle{+2{\rm Tr}^{\bullet}_{g_{\mathcal H}}
\left(\langle({\mathcal A}^{\phi}_{\bullet}X,{\mathcal A}^{\phi}_Y(A_{\mathcal H}\bullet)
\rangle+\langle({\mathcal A}^{\phi}_{\bullet}Y,{\mathcal A}^{\phi}_X
(A_{\mathcal H}\bullet)\rangle\right)}\\
\hspace{1.95truecm}\displaystyle{+2{\rm Tr}^{\bullet}_{g_{\mathcal H}}
\left(\langle(\widetilde{\nabla}_{\bullet}{\mathcal A}^{\phi})_{\xi}Y,
{\mathcal A}^{\phi}_{\bullet}
X\rangle+\langle(\widetilde{\nabla}_{\bullet}{\mathcal A}^{\phi})_{\xi}X,
{\mathcal A}^{\phi}_{\bullet}Y
\rangle\right)}\\
\hspace{1.95truecm}\displaystyle{+{\rm Tr}^{\bullet}_{g_{\mathcal H}}
\left(\langle(\widetilde{\nabla}_{\bullet}{\mathcal A}^{\phi})_{\bullet}X,
{\mathcal A}^{\phi}_{\xi}Y\rangle
+\langle(\widetilde{\nabla}_{\bullet}{\mathcal A}^{\phi})_{\bullet}Y,
{\mathcal A}^{\phi}_{\xi}X\rangle\right)}\\
\hspace{1.95truecm}\displaystyle{+2{\rm Tr}^{\bullet}_{g_{\mathcal H}}
\langle{\mathcal T}^{\phi}_{{\mathcal A}^{\phi}_{\bullet}X}\xi,
{\mathcal A}^{\phi}_{\bullet}Y\rangle,}
\end{array}\leqno{(4.1)}$$
where we omit $F_{\ast}$.
In particular, we have 
$$\begin{array}{l}
\displaystyle{{\mathcal R}(X,X)=4{\rm Tr}^{\bullet}_{g_{\mathcal H}}
\langle{\mathcal A}^{\phi}_{\bullet}X,{\mathcal A}^{\phi}_{\bullet}(A_{\mathcal H}X)
\rangle
+4{\rm Tr}^{\bullet}_{g_{\mathcal H}}\langle{\mathcal A}^{\phi}_{\bullet}X,
{\mathcal A}^{\phi}_X(A_{\mathcal H}\bullet)\rangle}\\
\hspace{1.95truecm}\displaystyle{
+3{\rm Tr}^{\bullet}_{g_{\mathcal H}}\langle(\widetilde{\nabla}_{\bullet}
{\mathcal A}^{\phi})_{\xi}X,{\mathcal A}^{\phi}_{\bullet}X\rangle
+2{\rm Tr}^{\bullet}_{g_{\mathcal H}}\langle(\widetilde{\nabla}_{\bullet}
{\mathcal A}^{\phi})_{\bullet}X,{\mathcal A}^{\phi}_{\xi}X
\rangle}\\
\hspace{1.95truecm}\displaystyle{+{\rm Tr}^{\bullet}_{g_{\mathcal H}}
\langle{\mathcal A}^{\phi}_{\bullet}X,
(\widetilde{\nabla}_X{\mathcal A}^{\phi})_{\xi}\bullet\rangle}
\end{array}$$
and hence 
$${\rm Tr}^{\bullet}_{g_{\mathcal H}}{\mathcal R}(\bullet,\bullet)=0.$$
}

\vspace{0.5truecm}

\noindent
{\it Simple proof of the third relation.} We give a simple proof of 
${\rm Tr}^{\bullet}_{g_{\mathcal H}}{\mathcal R}(\bullet,\bullet)=0$.  
Take any $(u,t)\in M\times[0,T)$ and an orthonormal base $(e_1,\cdots,e_n)$ of ${\mathcal H}_{(u,t)}$ with respect 
to $g_{(u,t)}$.  
According to Lemma 4.4 and the definiton of ${\mathcal R}$, we  have 
$$\begin{array}{l}
\displaystyle{({\rm Tr}^{\bullet}_{g_{\mathcal H}}{\mathcal R}(\bullet,\bullet))_{(u,t)}
=({\rm Tr}^{\bullet}_{g_{\mathcal H}}(\triangle_{\mathcal H}h)_{\mathcal H}(\bullet,\bullet))_{(u,t)}
-({\rm Tr}^{\bullet}_{g_{\mathcal H}}(\triangle^{\mathcal H}_{\mathcal H}h_{\mathcal H})(\bullet,\bullet))_{(u,t)}}\\
\displaystyle{=(\triangle_{\mathcal H}H^s)_{(u,t)}
-(\triangle^{\mathcal H}_{\mathcal H}H^s)_{(u,t)}
=\sum_{i=1}^n\left((\nabla dH^s)(e_i,e_i)-(\nabla^{\mathcal H}dH^s)(e_i,e_i)
\right)}\\
\displaystyle{=-\sum_{i=1}^n({\mathcal A}_{e_i}e_i)H^s=0,}
\end{array}$$
where we use $H^s=\sum\limits_{i=1}^nh(e_i,e_i)$ (which holds because the fibres of $\phi$ is 
regularized minimal).  \qed

\vspace{0.5truecm}

According to the proof of Corollary 4.11 in \cite{Koi2}, we obatin the following evolution equation.  

\vspace{0.5truecm}

\noindent
{\bf Lemma 4.9.} {\sl The norms $H^s_t$'s of $H_t$ 
satisfy the following evolution equation:
$$
\frac{\partial H^s}{\partial t}=
\triangle_{\mathcal H}H^s+H^s
\vert\vert A_{\mathcal H}\vert\vert^2-3H^s
{\rm Tr}(({\mathcal A}^{\phi}_{\xi})^2)_{\mathcal H}.$$
}

\vspace{0.5truecm}

According to the proof of Corollary 4.12 in \cite{Koi2}, we obtain the following evolution equation.  

\vspace{0.5truecm}

\noindent
{\bf Lemma 4.10.} {\sl The quantities $\vert\vert(A_{\mathcal H})_t\vert\vert^2$'s 
satisfy the following evolution equation:
$$\begin{array}{l}
\displaystyle{\frac{\partial\vert\vert A_{\mathcal H}\vert\vert^2}{\partial t}
=\triangle_{\mathcal H}(\vert\vert A_{\mathcal H}\vert\vert^2)
-2\vert\vert\nabla^{\mathcal H}A_{\mathcal H}\vert\vert^2}\\
\hspace{2.2truecm}\displaystyle{
+2\vert\vert A_{\mathcal H}\vert\vert^2\left(
\vert\vert A_{\mathcal H}\vert\vert^2-{\rm Tr}(({\mathcal A}^{\phi}_{\xi})^2)_{\mathcal H}\right)}\\
\hspace{2.2truecm}\displaystyle{
-4H^s{\rm Tr}\left((({\mathcal A}^{\phi}_{\xi})^2)\circ 
A_{\mathcal H}\right)
-2{\rm Tr}^{\bullet}_{g_{\mathcal H}}{\mathcal R}(A_{\mathcal H}\bullet,\bullet).}
\end{array}$$
}

\vspace{0.5truecm}

From Lemmas 4.9 and 4.10, we obtain the following evolution equation.  

\vspace{0.5truecm}

\noindent
{\bf Lemma 4.11.} 
{\sl The quantities $\vert\vert(A_{\mathcal H})_t\vert\vert^2
-\frac{(H^s_t)^2}{n}$'s satisfy the following evolution equation:
$$\begin{array}{l}
\displaystyle{
\frac{\partial(\vert\vert(A_{\mathcal H})_t\vert\vert^2
-\frac{(H^s)^2}{n})}{\partial t}=
\triangle_{\mathcal H}\left(\vert\vert A_{\mathcal H}\vert\vert^2
-\frac{(H^s)^2}{n}\right)
+\frac2n\left\vert\left\vert{\rm grad}
H^s\right\vert\right\vert^2}\\
\hspace{3.8truecm}\displaystyle{
+2\vert\vert A_{\mathcal H}\vert\vert^2\times\left(\vert\vert A_{\mathcal H}\vert\vert^2
-\frac{(H^s)^2}{n}\right)-2\vert\vert\nabla^{\mathcal H}A_{\mathcal H}\vert\vert^2}\\
\hspace{3.8truecm}\displaystyle{
-2{\rm Tr}(({\mathcal A}^{\phi}_{\xi})^2)_{\mathcal H}\times
\left(\vert\vert A_{\mathcal H}\vert\vert^2
-\frac{(H^s)^2}{n}\right)}\\
\hspace{3.8truecm}\displaystyle{
-4H^s\left({\rm Tr}\left(({\mathcal A}^{\phi}_{\xi})^2\circ
\left(A_{\mathcal H}-\frac{H^s}{n}{\rm id}\right)\right)\right)
}\\
\hspace{3.8truecm}\displaystyle{
-2{\rm Tr}^{\bullet}_{g_{\mathcal H}}{\mathcal R}\left(\left(
A_{\mathcal H}-\frac{H^s}{n}{\rm id}\right)\bullet,\bullet
\right),}
\end{array}$$
where ${\rm grad}H^s$ is the gradient vector field of 
$H^s$ with respect to $g$ and 
$\left\vert\left\vert{\rm grad}H^s\right\vert\right\vert$ 
is the norm of ${\rm grad}H^s$ with respect to $g$.}

\vspace{0.5truecm}

Set $n:={\rm dim}\,{\mathcal H}={\rm dim}\,\overline M$ and denote by $\bigwedge^n{\mathcal H}^{\ast}$ 
the exterior product bundle of degree $n$ of ${\mathcal H}^{\ast}$.  
Let $d\mu_{g_{\mathcal H}}$ be the section of $\pi_M^{\ast}(\bigwedge^n{\mathcal H}^{\ast})$ 
such that $(d\mu_{g_{\mathcal H}})_{(u,t)}$ is the volume element of $(g_{\mathcal H})_{(u,t)}$ for any 
$(u,t)\in M\times[0,T)$.  Then we can derive the following evolution equation for 
$\{(d\mu_{g_{\mathcal H}})_{(\cdot,t)}\}_{t\in[0,T)}$.  

\vspace{0.5truecm}

\noindent
{\bf Lemma 4.12.} {\sl The family $\{(d\mu_{g_{\mathcal H}})_{(\cdot,t)}\}_{t\in[0,T)}$ satisfies 
$$\frac{\partial\mu_{g_{\mathcal H}}}{\partial t}=-(H^s)^2\cdot d\mu_{g_{\mathcal H}}.$$
}

\vspace{0.5truecm}

\noindent
{\it Proof.} Let $(e_1,\cdots,e_n)$ be a local orthonormal base of $\mathcal H_{(u_0,t_0)}$ with respect to 
$(g_{\mathcal H})_{(u_0,t_0)}$ and $(E_1,\cdots,E_n)$ a local frame field of $\mathcal H_{t_0}|_U$ 
($U\,:\,$ an open set of $M$) with $(E_i)_{u_0}=e_i$ ($i=1,\cdots,n$) 
and $(\omega_1,\cdots,\omega_n)$ the dual frame field of $((\overline{E}_1)_{\mathcal H},\cdots,
(\overline{E}_n)_{\mathcal H})$, where $\overline E_i$ is the local section of $\pi_M^{\ast}TM|_{U\times[0,T)}$ 
defined by $(\overline{E}_i)_{(u,t)}:=(E_i)_u$ ($(u,t)\in U\times[0,T)$) and 
$(\overline{E}_i)_{\mathcal H}$ denotes the $\mathcal H$-component of $\overline{E}_i$.  
Here we note that $((\overline{E}_i)_{\mathcal H})_{(u,t_0)}=(\overline{E}_i)_{(u,t_0)}$ ($u\in U$) but 
$((\overline{E}_i)_{\mathcal H})_{(u,t)}\not=(\overline{E}_i)_{(u,t)}$ ($u\in U,t\not=t_0$).  
Then $d\mu_{g_{\mathcal H}}$ is described as 
$$d\mu_{g_{\mathcal H}}=\sqrt{{\rm det}(g_{\mathcal H}((\overline{E}_i)_{\mathcal H},(\overline{E}_j)_{\mathcal H}))}
\,\omega_1\wedge\cdots\wedge\omega_n.\leqno{(4.2)}$$
For simplicity, set 
$(g_{\mathcal H})_{ij}:=g_{\mathcal H}((\overline{E}_i)_{\mathcal H},(\overline{E}_j)_{\mathcal H})$ and 
$(h_{\mathcal H})_{ij}:=h_{\mathcal H}((\overline{E}_i)_{\mathcal H},(\overline{E}_j)_{\mathcal H})$.  
By using Lemma 4.2, we can show 
\begin{align*}
\frac{\partial(g_{\mathcal H})_{ij}}{\partial t}
=&\frac{\partial g_{\mathcal H}}{\partial t}((\overline{E}_i)_{\mathcal H},(\overline{E}_j)_{\mathcal H})
+g_{\mathcal H}(\nabla_{\frac{\partial}{\partial t}}(\overline{E}_i)_{\mathcal H},(\overline{E}_j)_{\mathcal H})\\
&+g_{\mathcal H}((\overline{E}_i)_{\mathcal H},\nabla_{\frac{\partial}{\partial t}}(\overline{E}_i)_{\mathcal H})\\
=&\frac{\partial g_{\mathcal H}}{\partial t}((\overline{E}_i)_{\mathcal H},(\overline{E}_j)_{\mathcal H})
=-2H^sh_{\mathcal H}((\overline{E}_i)_{\mathcal H},(\overline{E}_j)_{\mathcal H}),
\end{align*}
where we used $\nabla_{\frac{\partial}{\partial t}}(\overline{E}_i)_{\mathcal H}\in\mathcal V$ 
in the second equality.  By using this relation, we can derive 
\begin{align*}
&\frac{\partial\,d\mu_{g_{\mathcal H}}}{\partial t}
=\frac{\partial}{\partial t}\left(\sqrt{{\rm det}((g_{\mathcal H})_{ij})}\,\omega_1\wedge\cdots\omega_n\right)\\
=&\frac{1}{2\sqrt{{\rm det}((g_{\mathcal H})_{ij})}}\cdot\sum_{j=1}^n\left|
\begin{array}{ccccc}
(g_{\mathcal H})_{11} & \cdots & \displaystyle{\frac{\partial(g_{\mathcal H})_{1j}}{\partial t}} & \cdots & 
(g_{\mathcal H})_{1n}\\
\vdots & \vdots & \vdots & \vdots & \vdots\\
(g_{\mathcal H})_{n1} & \cdots & \displaystyle{\frac{\partial(g_{\mathcal H})_{nj}}{\partial t}} & \cdots & 
(g_{\mathcal H})_{nn}
\end{array}
\right|\,\,\omega_1\wedge\cdots\omega_n\\
&+\sqrt{{\rm det}((g_{\mathcal H})_{ij})}\cdot\sum_{j=1}^n\omega_1\wedge\cdots\wedge
\frac{\partial\omega_j}{\partial t}\wedge\cdots\wedge\omega_n\\
=&\frac{1}{2\sqrt{{\rm det}((g_{\mathcal H})_{ij})}}\cdot\sum_{j=1}^n\left|
\begin{array}{ccccc}
(g_{\mathcal H})_{11} & \cdots & -2g_{\mathcal H}((h_{\mathcal H})_{1j},H) & \cdots & (g_{\mathcal H})_{1n}\\
\vdots & \vdots & \vdots & \vdots & \vdots\\
(g_{\mathcal H})_{n1} & \cdots & -2g_{\mathcal H}((h_{\mathcal H})_{nj},H) & \cdots & (g_{\mathcal H})_{nn}
\end{array}
\right|\,\,\omega_1\wedge\cdots\wedge\omega_n\\
&+\sqrt{{\rm det}((g_{\mathcal H})_{ij})}\cdot\sum_{j=1}^n\omega_1\wedge\cdots\wedge
\frac{\partial\omega_j}{\partial t}\wedge\cdots\wedge\omega_n\\
=&-(H^s)^2\,d\mu_{g_{\mathcal H}}+\sqrt{{\rm det}((g_{\mathcal H})_{ij})}\cdot\sum_{j=1}^n\omega_1\wedge\cdots\wedge
\frac{\partial\omega_j}{\partial t}\wedge\cdots\wedge\omega_n.
\end{align*}
On the other hand, by using $\nabla_{\frac{\partial}{\partial t}}(\overline{E}_i)_{\mathcal H}\in\mathcal V$ agian, 
we can show 
$$\frac{\partial\omega_j}{\partial t}((\overline{E}_i)_{\mathcal H})
=-\omega_i\left(\nabla_{\frac{\partial}{\partial t}}((\overline{E}_i)_{\mathcal H}\right)=0.$$
Therefore we obtain the desired evolution equation.  \qed

\section{Sobolev inequality for Riemannian suborbifolds} 
In this section, we prove the divergence theorem for a compact Riemannian orbifold and 
Sobolev inequality for a compact Riemannian suborbifold, which may have the boundary.  
Let $(\overline M,\overline g)$ be an $n$-dimensional compact Riemannian orbifold, $\Sigma$ the singular set of 
$(\overline M,\overline g)$ and $\{\Sigma_1,\cdots,\Sigma_k\}$ be the set of all connected components of $\Sigma$.   
Set $\overline M':=\overline M\setminus\Sigma$.  
For a function $\rho$ over $\overline M$, we call $\displaystyle{\int_{\overline M'}\rho\,dv_{\overline g}}$ 
the {\it integral of} $\rho$ {\it over} $\overline M$ and denote it by 
$\displaystyle{\int_{\overline M}\rho\,dv_{\overline g}}$.  
First we prove the divergence theorem for orbitangent vector fields on a compact Riemannian orbifold.  

\vspace{0.5truecm}

\noindent
{\bf Theorem 5.1.} {\sl For any $C^1$-orbitangent vector field $X$ on $(\overline M,\overline g)$, 
the relation $\displaystyle{\int_{\overline M}{\rm div}_{\overline g}\,X\,dv_{\overline g}=0}$ holds.}

\vspace{0.5truecm}

\noindent
{\it Proof.} 
Let $U_i$ ($i=1,\cdots,k$) be a sufficiently small tubular neighborhood of $\Sigma_i$ with 
$U_i\cap U_j=\emptyset$ ($i\not=j$).  Set $\displaystyle{W:=\overline M\setminus
\left(\mathop{\cup}_{i=1}^kU_i\right)}$.  
Take families $\{(U_{ij},\varphi_{ij},\widehat U_{ij}/\Gamma_{ij})\,|\,j=1,\cdots,m_i\}$ of orbifold charts of 
$\overline M$ such that the family $\{{\rm cl}(U_{ij})\}_{j=1}^{m_i}$ of the clusure ${\rm cl}(U_{ij})$ of $U_{ij}$ 
gives a division of ${\rm cl}(U_i)$ ($i=1,\cdots,k$).  
Denote by $\pi_{ij}$ the projection $\pi_{\Gamma_{ij}}:\widehat U_{ij}\to\widehat U_{ij}/\Gamma_{ij}$ and 
$l_i$ the cardinal number of $\Gamma_{ij}$, which depends only on $i$.  
Let $\xi_i$ be the outward unit normal vector field of $\partial U_i$ and 
$\iota_i$ is the inclusion map of $\partial U_i$ into $\overline M$.  
Also, let $\xi_{ij}$ be the outward unit normal vector field of $\partial U_{ij}$ satisfying 
$\xi_{ij}|_{\partial U_i\cap\partial U_{ij}}=\xi_i|_{\partial U_i\cap\partial U_{ij}}$.  
Also, let $\widehat{\xi}_{ij}$ be the unit normal vector field of $\partial\widehat U_{ij}$ satisfying 
$(\varphi_{ij}^{-1}\circ\pi_{ij})_{\ast}(\widehat{\xi}_{ij})=\xi_{ij}\vert_{\partial U_{ij}}$ 
and let $\widehat X_{ij}$ be the vector field on $\widehat U_{ij}$ satisfying 
$(\varphi_{ij}^{-1}\circ\pi_{ij})_{\ast}(\widehat X_{ij})=X\vert_{U_{ij}}$.  
Denote by $\iota_{ij}$ the inclusion map of $\partial\widehat U_{ij}$ into $\mathbb R^n$ and 
$\widehat{\overline g}_{ij}$ the local lift of $\overline g$ with respect to 
$(U_{ij},\varphi_{ij},\widehat U_{ij}/\Gamma_{ij})$.  
Then, by using the divergence theorem (for a compact Riemannian manifold with boundary), we have 
$$
\begin{array}{l}
\displaystyle{\int_{\overline M}{\rm div}_{\overline g}\,X\,dv_{\overline g}
=\int_W{\rm div}_{\overline g}\,X\,dv_{\overline g}+\sum_{i=1}^k\sum_{j=1}^{m_i}
\int_{U_{ij}}{\rm div}_{\overline g}\,X\,dv_{\overline g}}\\
\hspace{2.5truecm}\displaystyle{=-\sum_{i=1}^k\int_{\partial U_i}\overline g(X,\xi_i)\,
dv_{\overline{\iota}_i^{\ast}\overline g}
+\sum_{i=1}^k\sum_{j=1}^{m_i}\int_{U_{ij}}{\rm div}_{\overline g}\,X\,dv_{\overline g}.}
\end{array}\leqno{(5.1)}$$
Also, by using the divergence theorem, we can show 
$$%
%
\hspace{4.75truecm}}

\vspace{0.5truecm}

\centerline{{\bf Figure 5.1$\,\,:\,\,$ The relation of $\xi_i$ and $\xi_{ij}$}}

\vspace{1truecm}

In 1973, J.H. Michael and L.M. Simon (\cite{MS}) proved the Sobolev inequality for compact Riemannian submanifolds 
(which may have boundary) in a Euclidean space, where we note that the integrand must vanishes on 
the boundary of the submanifold.  In 1974, D. Hoffman and J. Spruck (\cite{HoSp}) proved the same Sobolev inequality 
in a general Riemannian manifold, where we note that the integrand must vanishes on the boundary of 
the submanifold and furthemore, the volume of the support of the integrand must satisfy some estimate from above 
related to the curvature and the injective radius of the ambient space.  
We shall show the following Sobolev inequality for compact Riemannian suborbifolds.  

\vspace{0.5truecm}

\noindent
{\bf Theorem 5.2.} {\sl Let $(\overline M,\overline g)$ be a compact Riemannian suborbifold isometrically immersed 
into $(N,g_N)$ by $f$.  Assume that the sectional curvature $K$ of a complete Riemannian orbifold $(N,g_N)$ 
satisfies $K\leq b^2$ ($b$ is a non-negative real number or the purely imaginary number), 
$\overline M$ satisfies the following condition $(\ast_1)$:

\vspace{0.15truecm}

\noindent
$(\ast_1)$ $f(\overline M)$ is included by $B_{\frac{\pi}{b}}(x_0)$ for some $x_0\in N$ and 
$\exp_{x_0}|_{B^T_{\frac{\pi}{b}}(x_0)}$ is injective.  

\vspace{0.15truecm}

Let $\rho$ be any non-negative $C^1$-function on $\overline M$ satisfying 
$$b^2(1-\alpha)^{-2/n}\left(\omega_n^{-1}\cdot l\cdot{\rm Vol}_{\overline g}({\rm supp}\,\rho)
\right)^{2/n}\leq 1,
\leqno{(\ast_2)}$$
where $\overline g$ denotes the induced metric on $\overline M$, $\omega_n$ denotes the volume of the unit ball 
in the Euclidean space ${\mathbb R}^n$, $l$ denotes the cardinality of the local group at $x_0$ and 
$\alpha$ is any fixed positive constant smaller than one.  
Then the following inequality for $\rho$ holds:
$$\left(\int_{\overline M}\rho^{\frac{n}{n-1}}dv_{\overline g}\right)^{\frac{n-1}{n}}
\leq C(n,\alpha)\int_{\overline M}\left(\vert\vert d\rho\vert\vert+\rho\cdot\|H\|\right)\,dv_{\overline g},
\leqno{(5.3)}$$
where $H$ denotes the mean curvature vector of $f$ and $C(n,\alpha)$ is the positive constant 
depending only on $n$ and $\alpha$.  
}

\vspace{0.15truecm}

\noindent
{\it Remark 5.1.} In the case where $(\overline M,\overline g)$ is a compact Riemannian manifold, 
the statement of this theorem follows from the Sobolev's inequality in [HoSp] because the condition $(\ast_1)$ 
assures the condition $(2.3)$ in Theorem 2.1 of [HoSp].  

\vspace{0.3truecm}

We shall prepare some lemmas to prove this theorem.  
Let $\Gamma$ be the local group at $x_0$ and set $\widehat x_0:=(\exp_{x_0}\circ\pi)^{-1}(x_0)$.  
The orbitangent space $T_{x_0}M$ is identified with 
the orbit space $\mathbb R^{n+1}/\Gamma$.  Let $\pi:\mathbb R^{n+1}\to T_{x_0}N$ be the orbit map.  
Set $\widehat B_{\frac{\pi}{b}}(\widehat x_0):=(\exp_{x_0}\circ\pi)^{-1}(B_{\frac{\pi}{b}}(x_0))$ and 
$\widehat{g_N}:=(\exp_{x_0}\circ\pi)^{\ast}g_N$.  
Also, set 
$$\widehat{\overline M}:=\{(x,\widehat y)\in\overline M\times\widehat B_{\frac{\pi}{b}}(\widehat x_0)\,|\,
f(x)=(\exp_{x_0}\circ\pi)(\widehat y)\}$$
and define a map $\widehat f:\widehat{\overline M}\hookrightarrow\widehat B_{\frac{\pi}{b}}(\widehat x_0)$ by 
$\widehat f(x,\widehat y)=\widehat y\,\,\,((x,\widehat y)\in\widehat{\overline M})$.  
Also, define a map $\pi_{\overline M}:\widehat{\overline M}\to\overline M$ by 
$\pi_{\overline M}(x,\widehat y):=x\,\,\,((x,\widehat y)\in\widehat{\overline M})$.  
It is clear that $\widehat{\overline M}$ is a $C^{\infty}$-manifold and that $\widehat f$ is 
a $C^{\infty}$-immersion.  Also, it is clear that $\pi_{\overline M}$ is a $C^{\infty}$-orbisubmersion and that 
$\exp_{x_0}\circ\pi\circ\widehat f=f\circ\pi_{\overline M}$ holds.  
Let $\widehat{\overline g}$ be the Riemannian metric on $\widehat{\overline M}$ such that 
$\pi_{\overline M}:(\widehat{\overline M},\widehat{\overline g})\to(\overline M,\overline g)$ is a Riemannian 
orbisubmersion.  Let $\widehat{\nabla^N}$ be the Riemannian connection of $\widehat{g_N}$ and $\widehat{\overline{\nabla}}$ that of 
$\widehat{\overline g}$.  Let $X\in\Gamma^{\infty}({\widehat f}^{\ast}T(\widehat B_{\frac{\pi}{b}}(\widehat x_0)))$.  
Let $X^T$ (resp. $X^{\perp}$) be the tangential (resp. the normal) component of $X$, that is, 
$$X_{\widehat x}=f_{\ast\widehat x}(X^T_{\widehat x})+X^{\perp}_{\widehat x}\qquad
(X^T_{\widehat x}\in T_{\widehat x}\widehat{\overline M},\,\,X^{\perp}_{\widehat x}\in T^{\perp}_{\widehat x}
\widehat{\overline M}).$$
Define ${\rm div}_{\widehat f}X\in C^{\infty}(\widehat{\overline M})$ by 
$$({\rm div}_{\widehat f}X)_{\widehat x}
:=({\rm Tr}_{\widehat{\overline g}}((\widehat{\nabla^N})^{\widehat f}X)^T)_{\widehat x}
=\sum_{i=1}^n\widehat{\overline g}((\widehat{\nabla^N})^{\widehat f}_{e_i}X,\widehat f_{\ast\widehat x}(e_i))
\quad(\widehat x\in\widehat{\overline M}),$$
where $(e_1,\cdots,e_n)$ is an orthonormal base of $T_{\widehat x}\widehat{\overline M}$ with respect to 
$\widehat{\overline g}_{\widehat x}$ and $(\widehat{\nabla^N})^{\widehat f}$ denotes the induced connection of 
$\widehat{\nabla^N}$ by $\widehat f$.  

\vspace{0.25truecm}

\centerline{
\unitlength 0.1in
%
\hspace{4truecm}}

\vspace{0.25truecm}

\centerline{{\bf Figure 5.2$\,\,:\,\,$ The lift of an orbisubmanifold in a singular geodesic ball (I)}}


\vspace{0.5truecm}

{\small
\centerline{
\unitlength 0.1in
%
%
\hspace{3.5truecm}}
}

\vspace{0.35truecm}

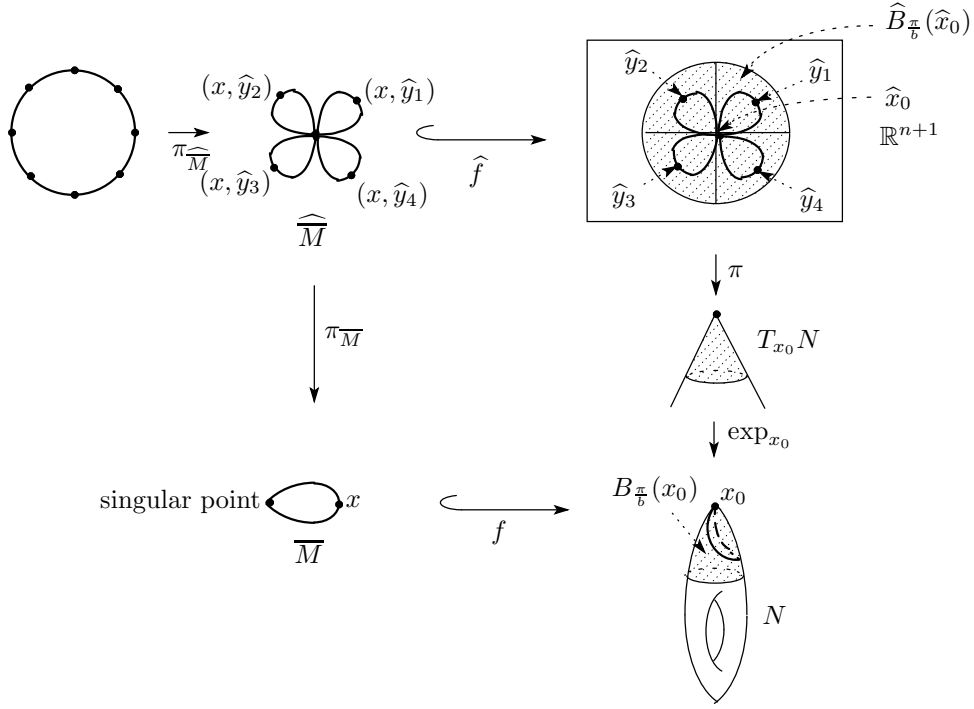
\centerline{{\bf Figure 5.3$\,\,:\,\,$ The lift of an orbisubmanifold in a singular geodesic ball (II)}}

\vspace{0.5truecm}

First we prepare the following lemma.  


\vspace{0.5truecm}

\noindent
{\bf Lemma 5.3.} {\sl {\rm(i)} ${\rm div}_{\widehat{\overline g}}X^T={\rm div}_{\widehat f}X+\widehat{g_N}(X,H)$,
where $H$ denotes the mean curvature vector of $\widehat f$.  

{\rm (ii)}\ For $\rho\in C^{\infty}(\widehat{\overline M})$, we have 
$${\rm div}_{\widehat f}(\rho X)=\rho{\rm div}_{\widehat f}X
+\widehat{\overline g}(X^T,{\rm grad}_{\widehat{\overline g}}\,\rho),$$
where ${\rm grad}_{\widehat{\overline g}}(\bullet)$ denotes the gradient vector field of $(\bullet)$ with respect to 
$\widehat{\overline g}$.}

\vspace{0.5truecm}

See the proof of Lemma 3.2 in \cite{HoSp} about the proof of this lemma.  
Let $r_{\widehat x_0}:\widehat B_{\frac{\pi}{b}}(\widehat x_0)\to[0,\infty)$ be the distance 
function from $\widehat x_0$ with respect to $\widehat{g_N}$, that is, 
$r_{\widehat x_0}(\widehat x):=d_{\widehat{g_N}}(\widehat x_0,\widehat x)\,\,\,
(\widehat x\in\widehat B_{\frac{\pi}{b}}(\widehat x_0))$.  
Define a $C^{\infty}$-vector field $\mathbb P$ on $\widehat B_{\frac{\pi}{b}}(\widehat x_0)$ by 
$${\mathbb P}_{\widehat x}:=r_{\widehat x_0}(\widehat x)\cdot
({\rm grad}_{\widehat{g_N}}r_{\widehat x_0})_{\widehat x}\qquad(\widehat x\in\widehat B_{\frac{\pi}{b}}
(\widehat x_0)),$$
where ${\rm grad}_{\widehat{g_N}}(\bullet)$ denotes the gradient vector field of $\bullet$ with respect to 
$\widehat{g_N}$.  Also, define a $C^{\infty}$-vector field $\widetilde{\mathbb P}$ over the tangent space 
$T_{\widehat x_0}\mathbb R^{n+1}$ by 
$$\widetilde{\mathbb P}_v:=v\quad(v\in T_{\widehat x_0}\mathbb R^{n+1}).$$
By using the discussion in the proof Lemmas 3.5 and 3.6 in \cite{HoSp}, we can show the following fact.  

\vspace{0.5truecm}

\noindent
{\bf Lemma 5.4.} {\sl {\rm (i)} For any unit vector $v$ of $\widehat B_{\frac{\pi}{b}}(\widehat x_0)$ at any 
$\widehat x\in\widehat B_{\frac{\pi}{b}}(\widehat x_0)$, the following inequality holds:
$$\widehat{g_N}(\widehat{\nabla^N}_v{\mathbb P},v)\geq b\cdot r_{\widehat x_0}(\widehat x)\cdot
\cot(b\cdot r_{\widehat x_0}(\widehat x)).$$

{\rm (ii)} For any $(x,\widehat y)\in\widehat{\overline M}$, the following inequality holds:
$$({\rm div}_{\widehat f}({\mathbb P}\circ\widehat f))_{(x,\widehat y)}\geq n\cdot b\cdot r_{\widehat x_0}
(\widehat y)\cdot\cot(b\cdot r_{\widehat x_0}(\widehat y)).$$
}

\vspace{0.5truecm}

{\small
\centerline{
\unitlength 0.1in
\begin{picture}( 70.4000, 26.1200)(-22.5000,-28.5200)
%
\special{pn 8}%
\special{ar 4520 2448 120 230  1.8325098 4.4572066}%
%
\special{pn 8}%
\special{ar 4370 2438 130 210  5.3697781 6.2831853}%
\special{ar 4370 2438 130 210  0.0000000 0.9758854}%
%
\special{pn 8}%
\special{ar 4590 2354 310 690  3.1415927 4.2698110}%
%
\special{pn 8}%
\special{ar 4320 2354 310 690  5.1549669 6.2831853}%
%
\special{pn 8}%
\special{ar 4560 2354 280 520  1.9223130 3.1415927}%
%
\special{pn 8}%
\special{ar 4350 2364 280 520  6.2831853 6.2831853}%
\special{ar 4350 2364 280 520  0.0000000 1.2192797}%
%
\special{pn 20}%
\special{sh 1}%
\special{ar 4450 1744 10 10 0  6.28318530717959E+0000}%
\special{sh 1}%
\special{ar 4450 1744 10 10 0  6.28318530717959E+0000}%
%
\special{pn 8}%
\special{ar 4450 2134 160 40  6.2831853 6.2831853}%
\special{ar 4450 2134 160 40  0.0000000 3.1415927}%
%
\special{pn 8}%
\special{ar 4440 2134 160 40  3.1415927 3.2615927}%
\special{ar 4440 2134 160 40  3.6215927 3.7415927}%
\special{ar 4440 2134 160 40  4.1015927 4.2215927}%
\special{ar 4440 2134 160 40  4.5815927 4.7015927}%
\special{ar 4440 2134 160 40  5.0615927 5.1815927}%
\special{ar 4440 2134 160 40  5.5415927 5.6615927}%
\special{ar 4440 2134 160 40  6.0215927 6.1415927}%
%
\special{pn 20}%
\special{sh 1}%
\special{ar 3170 1240 10 10 0  6.28318530717959E+0000}%
\special{sh 1}%
\special{ar 3170 1240 10 10 0  6.28318530717959E+0000}%
%
\special{pn 8}%
\special{ar 3170 1588 170 40  6.2831853 6.2831853}%
\special{ar 3170 1588 170 40  0.0000000 3.1415927}%
%
\special{pn 8}%
\special{ar 3170 1598 170 40  3.1415927 3.2558784}%
\special{ar 3170 1598 170 40  3.5987355 3.7130212}%
\special{ar 3170 1598 170 40  4.0558784 4.1701641}%
\special{ar 3170 1598 170 40  4.5130212 4.6273069}%
\special{ar 3170 1598 170 40  4.9701641 5.0844498}%
\special{ar 3170 1598 170 40  5.4273069 5.5415927}%
\special{ar 3170 1598 170 40  5.8844498 5.9987355}%
%
\special{pn 8}%
\special{pa 2014 528}%
\special{pa 1460 1028}%
\special{pa 2614 1028}%
\special{pa 3120 528}%
\special{pa 3120 528}%
\special{pa 2014 528}%
\special{fp}%
%
\special{pn 8}%
\special{ar 2304 778 400 150  0.0364703 6.2831853}%
%
\special{pn 20}%
\special{sh 1}%
\special{ar 2304 778 10 10 0  6.28318530717959E+0000}%
\special{sh 1}%
\special{ar 2304 778 10 10 0  6.28318530717959E+0000}%
%
\special{pn 8}%
\special{pa 1904 778}%
\special{pa 2704 778}%
\special{fp}%
%
\special{pn 8}%
\special{pa 2044 888}%
\special{pa 2574 658}%
\special{fp}%
%
\special{pn 8}%
\special{pa 2444 918}%
\special{pa 2174 638}%
\special{fp}%
%
\special{pn 20}%
\special{pa 2440 778}%
\special{pa 2560 778}%
\special{fp}%
\special{sh 1}%
\special{pa 2560 778}%
\special{pa 2494 758}%
\special{pa 2508 778}%
\special{pa 2494 798}%
\special{pa 2560 778}%
\special{fp}%
%
\special{pn 20}%
\special{pa 2630 778}%
\special{pa 2850 778}%
\special{fp}%
\special{sh 1}%
\special{pa 2850 778}%
\special{pa 2784 758}%
\special{pa 2798 778}%
\special{pa 2784 798}%
\special{pa 2850 778}%
\special{fp}%
%
\special{pn 20}%
\special{pa 2170 778}%
\special{pa 2040 768}%
\special{fp}%
\special{sh 1}%
\special{pa 2040 768}%
\special{pa 2106 792}%
\special{pa 2094 772}%
\special{pa 2108 752}%
\special{pa 2040 768}%
\special{fp}%
%
\special{pn 20}%
\special{pa 1980 778}%
\special{pa 1720 778}%
\special{fp}%
\special{sh 1}%
\special{pa 1720 778}%
\special{pa 1788 798}%
\special{pa 1774 778}%
\special{pa 1788 758}%
\special{pa 1720 778}%
\special{fp}%
%
\special{pn 20}%
\special{pa 2380 728}%
\special{pa 2480 688}%
\special{fp}%
\special{sh 1}%
\special{pa 2480 688}%
\special{pa 2412 694}%
\special{pa 2430 708}%
\special{pa 2426 730}%
\special{pa 2480 688}%
\special{fp}%
%
\special{pn 20}%
\special{pa 2530 678}%
\special{pa 2740 568}%
\special{fp}%
\special{sh 1}%
\special{pa 2740 568}%
\special{pa 2672 580}%
\special{pa 2694 592}%
\special{pa 2690 616}%
\special{pa 2740 568}%
\special{fp}%
%
\special{pn 20}%
\special{pa 2220 818}%
\special{pa 2140 858}%
\special{fp}%
\special{sh 1}%
\special{pa 2140 858}%
\special{pa 2210 846}%
\special{pa 2188 834}%
\special{pa 2192 810}%
\special{pa 2140 858}%
\special{fp}%
%
\special{pn 20}%
\special{pa 2090 868}%
\special{pa 1880 978}%
\special{fp}%
\special{sh 1}%
\special{pa 1880 978}%
\special{pa 1948 964}%
\special{pa 1928 952}%
\special{pa 1930 928}%
\special{pa 1880 978}%
\special{fp}%
%
\special{pn 20}%
\special{pa 2270 738}%
\special{pa 2210 688}%
\special{fp}%
\special{sh 1}%
\special{pa 2210 688}%
\special{pa 2248 746}%
\special{pa 2252 722}%
\special{pa 2274 714}%
\special{pa 2210 688}%
\special{fp}%
%
\special{pn 20}%
\special{pa 2190 658}%
\special{pa 2070 538}%
\special{fp}%
\special{sh 1}%
\special{pa 2070 538}%
\special{pa 2104 598}%
\special{pa 2108 576}%
\special{pa 2132 570}%
\special{pa 2070 538}%
\special{fp}%
%
\special{pn 20}%
\special{pa 2350 828}%
\special{pa 2400 878}%
\special{fp}%
\special{sh 1}%
\special{pa 2400 878}%
\special{pa 2368 816}%
\special{pa 2362 840}%
\special{pa 2340 844}%
\special{pa 2400 878}%
\special{fp}%
%
\special{pn 20}%
\special{pa 2430 908}%
\special{pa 2560 1028}%
\special{fp}%
\special{sh 1}%
\special{pa 2560 1028}%
\special{pa 2526 968}%
\special{pa 2522 992}%
\special{pa 2498 996}%
\special{pa 2560 1028}%
\special{fp}%
\put(28.0000,-11.9000){\makebox(0,0)[rt]{$\pi$}}%
\put(38.7000,-18.2000){\makebox(0,0)[rt]{$\exp_{x_0}$}}%
\put(34.0000,-13.2700){\makebox(0,0)[lt]{$T_{x_0}N$}}%
\put(47.2000,-23.0000){\makebox(0,0)[lt]{$N$}}%
%
\special{pn 8}%
\special{ar 1460 728 570 250  4.8754032 4.9046715}%
\special{ar 1460 728 570 250  4.9924764 5.0217447}%
\special{ar 1460 728 570 250  5.1095495 5.1388178}%
\special{ar 1460 728 570 250  5.2266227 5.2558910}%
\special{ar 1460 728 570 250  5.3436959 5.3729642}%
\special{ar 1460 728 570 250  5.4607690 5.4900373}%
\special{ar 1460 728 570 250  5.5778422 5.6071105}%
\special{ar 1460 728 570 250  5.6949154 5.7241837}%
\special{ar 1460 728 570 250  5.8119886 5.8412568}%
\special{ar 1460 728 570 250  5.9290617 5.9583300}%
%
\special{pn 8}%
\special{pa 2020 668}%
\special{pa 2060 728}%
\special{dt 0.045}%
\special{sh 1}%
\special{pa 2060 728}%
\special{pa 2040 660}%
\special{pa 2030 684}%
\special{pa 2006 684}%
\special{pa 2060 728}%
\special{fp}%
%
\special{pn 8}%
\special{pa 4740 1810}%
\special{pa 4480 1930}%
\special{dt 0.045}%
\special{sh 1}%
\special{pa 4480 1930}%
\special{pa 4550 1920}%
\special{pa 4528 1908}%
\special{pa 4532 1884}%
\special{pa 4480 1930}%
\special{fp}%
\put(47.9000,-18.5000){\makebox(0,0)[lb]{$B_{\frac{\pi}{b}}(x_0)$}}%
\put(15.3000,-5.6000){\makebox(0,0)[rb]{$\widehat B_{\frac{\pi}{b}}(\widehat x_0)$}}%
\put(30.5000,-7.3000){\makebox(0,0)[lt]{$\mathbb R^{n+1}$}}%
%
\special{pn 20}%
\special{pa 3380 450}%
\special{pa 3560 450}%
\special{fp}%
\special{sh 1}%
\special{pa 3560 450}%
\special{pa 3494 430}%
\special{pa 3508 450}%
\special{pa 3494 470}%
\special{pa 3560 450}%
\special{fp}%
\put(36.2000,-3.8000){\makebox(0,0)[lt]{'s : $\mathbb P$}}%
%
\special{pn 8}%
\special{ar 2730 870 430 550  3.4041391 3.4286289}%
\special{ar 2730 870 430 550  3.5020982 3.5265880}%
\special{ar 2730 870 430 550  3.6000574 3.6245472}%
\special{ar 2730 870 430 550  3.6980166 3.7225064}%
\special{ar 2730 870 430 550  3.7959758 3.8204656}%
\special{ar 2730 870 430 550  3.8939350 3.9184248}%
\special{ar 2730 870 430 550  3.9918942 4.0163840}%
\special{ar 2730 870 430 550  4.0898533 4.1143431}%
\special{ar 2730 870 430 550  4.1878125 4.2123023}%
\special{ar 2730 870 430 550  4.2857717 4.3102615}%
%
\special{pn 8}%
\special{pa 2320 720}%
\special{pa 2300 770}%
\special{dt 0.045}%
\special{sh 1}%
\special{pa 2300 770}%
\special{pa 2344 716}%
\special{pa 2320 720}%
\special{pa 2306 702}%
\special{pa 2300 770}%
\special{fp}%
\put(26.2000,-4.1000){\makebox(0,0)[lb]{$\widehat x_0$}}%
%
\special{pn 4}%
\special{pa 4390 1854}%
\special{pa 4474 1770}%
\special{dt 0.027}%
\special{pa 4370 1934}%
\special{pa 4498 1806}%
\special{dt 0.027}%
\special{pa 4330 2034}%
\special{pa 4522 1842}%
\special{dt 0.027}%
\special{pa 4300 2124}%
\special{pa 4540 1884}%
\special{dt 0.027}%
\special{pa 4334 2150}%
\special{pa 4550 1934}%
\special{dt 0.027}%
\special{pa 4378 2166}%
\special{pa 4568 1976}%
\special{dt 0.027}%
\special{pa 4430 2174}%
\special{pa 4584 2020}%
\special{dt 0.027}%
\special{pa 4490 2174}%
\special{pa 4590 2074}%
\special{dt 0.027}%
\special{pa 4560 2164}%
\special{pa 4598 2126}%
\special{dt 0.027}%
%
\special{pn 4}%
\special{pa 1920 790}%
\special{pa 2040 670}%
\special{dt 0.027}%
\special{pa 1942 830}%
\special{pa 2122 648}%
\special{dt 0.027}%
\special{pa 1970 860}%
\special{pa 2190 640}%
\special{dt 0.027}%
\special{pa 2012 878}%
\special{pa 2254 636}%
\special{dt 0.027}%
\special{pa 2054 898}%
\special{pa 2320 630}%
\special{dt 0.027}%
\special{pa 2102 908}%
\special{pa 2380 630}%
\special{dt 0.027}%
\special{pa 2150 920}%
\special{pa 2432 640}%
\special{dt 0.027}%
\special{pa 2210 920}%
\special{pa 2484 646}%
\special{dt 0.027}%
\special{pa 2264 926}%
\special{pa 2532 658}%
\special{dt 0.027}%
\special{pa 2330 920}%
\special{pa 2578 674}%
\special{dt 0.027}%
\special{pa 2390 920}%
\special{pa 2616 694}%
\special{dt 0.027}%
\special{pa 2464 908}%
\special{pa 2652 718}%
\special{dt 0.027}%
\special{pa 2538 894}%
\special{pa 2680 752}%
\special{dt 0.027}%
\special{pa 2620 870}%
\special{pa 2670 820}%
\special{dt 0.027}%
%
\special{pn 8}%
\special{pa 3160 1240}%
\special{pa 2910 1760}%
\special{fp}%
%
\special{pn 8}%
\special{pa 3170 1250}%
\special{pa 3440 1770}%
\special{fp}%
%
\special{pn 4}%
\special{pa 3144 1278}%
\special{pa 3168 1254}%
\special{dt 0.027}%
\special{pa 3094 1388}%
\special{pa 3188 1294}%
\special{dt 0.027}%
\special{pa 3042 1498}%
\special{pa 3208 1332}%
\special{dt 0.027}%
\special{pa 3008 1594}%
\special{pa 3228 1372}%
\special{dt 0.027}%
\special{pa 3050 1610}%
\special{pa 3248 1412}%
\special{dt 0.027}%
\special{pa 3100 1620}%
\special{pa 3270 1452}%
\special{dt 0.027}%
\special{pa 3150 1630}%
\special{pa 3290 1492}%
\special{dt 0.027}%
\special{pa 3210 1630}%
\special{pa 3310 1532}%
\special{dt 0.027}%
\special{pa 3290 1610}%
\special{pa 3330 1570}%
\special{dt 0.027}%
%
\special{pn 8}%
\special{pa 2730 1080}%
\special{pa 2950 1270}%
\special{fp}%
\special{sh 1}%
\special{pa 2950 1270}%
\special{pa 2914 1212}%
\special{pa 2910 1236}%
\special{pa 2886 1242}%
\special{pa 2950 1270}%
\special{fp}%
%
\special{pn 8}%
\special{pa 3710 1660}%
\special{pa 4070 1930}%
\special{fp}%
\special{sh 1}%
\special{pa 4070 1930}%
\special{pa 4030 1874}%
\special{pa 4028 1898}%
\special{pa 4006 1906}%
\special{pa 4070 1930}%
\special{fp}%
\end{picture}%
\hspace{9truecm}}
}

\vspace{0.5truecm}

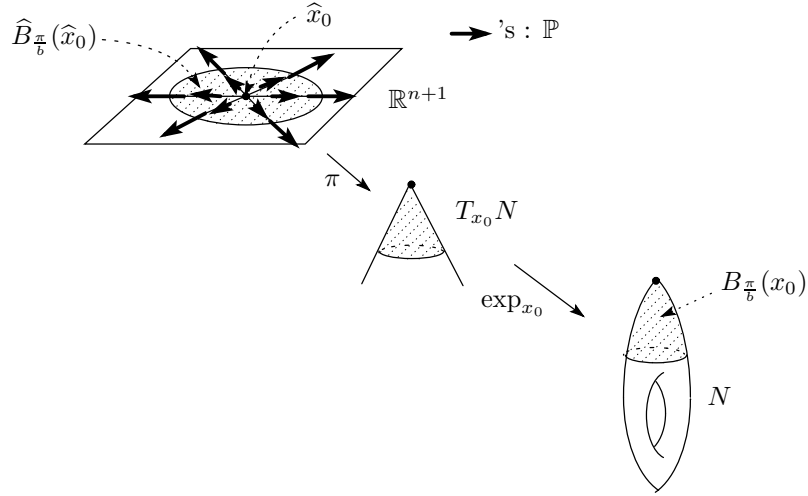
\centerline{{\bf Figure 5.4$\,\,:\,\,$ The position vector field on the orbicovering of a geodesic ball}}

\vspace{0.5truecm}

Let $\lambda$ be a $C^1$-function over ${\mathbb R}$ satisfying the following condition:
$$\lambda(t)=0\,\,\,(t\leq 0),\quad\,\,\lambda'(t)\geq 0\,\,\,(t\geq0).\leqno{(5.4)}$$
Set $\widehat B^{\widehat{\overline M}}_s(\widehat x_0):={\widehat f}^{-1}(\widehat B_s(\widehat x_0))$, where 
$\widehat B_s(\widehat x_0)$ denotes the geodesic ball of radius $s$ centered at $\widehat x_0$.  
Let $\rho$ be a $C^1$-function as in the statement of Theorem 5.2 and set 
$\widehat{\rho}:=\rho\circ\pi_{\overline M}$.  
Define functions $\phi_{\widehat{\rho},\widehat x_0,\lambda},\,\eta_{\widehat{\rho},\widehat x_0},\,
\overline{\phi}_{\widehat{\rho},\widehat x_0},\,\overline{\eta}_{\widehat{\rho},\widehat x_0}$ over 
$[0,\frac{\pi}{b})$ by 
\begin{align*}
&\phi_{\widehat{\rho},\widehat x_0,\lambda}(s):=\int_{\widehat{\overline M}}
\lambda(s-(r_{\widehat x_0}\circ\widehat f))\cdot\widehat{\rho}\,dv_{\widehat{\overline g}},\\
&\eta_{\widehat{\rho},\widehat x_0,\lambda}(s):=\int_{\widehat{\overline M}}
\lambda(s-(r_{\widehat x_0}\circ\widehat f))\cdot(|{\rm grad}_{\widehat{\overline g}}\rho|+\rho|H|)\,
dv_{\widehat{\overline g}},\\
&\overline{\phi}_{\widehat{\rho},\widehat x_0}(s):=\int_{\widehat B^{\widehat{\overline M}}_s(\widehat x_0)}
\widehat{\rho}\,dv_{\widehat{\overline g}},\\
&\overline{\eta}_{\widehat{\rho},\widehat x_0,\lambda}(s)
:=\int_{\widehat B^{\widehat{\overline M}}_s(\widehat x_0)}\lambda(s-(r_{\widehat x_0}\circ\widehat f))
\cdot(|{\rm grad}_{\widehat{\overline g}}\rho|+\rho|H|)\,dv_{\widehat{\overline g}}.
\end{align*}

According to the proof of Lemma 4.1 in \cite{HoSp}, we can derive the following fact by using Lemmas 5.3 and 5.4.  

\vspace{0.5truecm}

\noindent
{\bf Lemma 5.5.} {\sl For all $s\in[0,\frac{\pi}{b})$, the following inequality hold:
$$\left\{\begin{array}{ll}
\displaystyle{-\frac{d}{ds}\left((\sin(bs))^{-n}\phi_{\widehat{\rho},\widehat x_0,\lambda}(s)\right)\leq
(\sin(bs))^{-n}\eta_{\widehat{\rho},\widehat x_0,\lambda}(s)} & (b\,:\,{\rm real})\\
\displaystyle{-\frac{d}{ds}\left(s^{-n}\phi_{\widehat{\rho},\widehat x_0,\lambda}(s)\right)\leq 
s^{-n}\eta_{\widehat{\rho},\widehat x_0,\lambda}(s)} & (b\,:\,{\rm purely}\,\,{\rm imaginary}).
\end{array}\right.$$
}

\vspace{0.5truecm}

According to the proof of Lemma 4.2 in \cite{HoSp}, we can show the following result by using Lemma 5.5.  

\vspace{0.5truecm}

\noindent
{\bf Lemma 5.6.} {\sl Let $\alpha$ and $\hat{\alpha}$ be constants with $0<\alpha<1\leq\hat{\alpha}$, and 
$\{\lambda_{\varepsilon}\}_{\varepsilon>0}$ be $C^1$-functions over $\mathbb R$ satisfying the above condition 
$(5.4)$ and the following condition:
$$\lambda_{\varepsilon}\leq 1,\quad\lambda_{\varepsilon}^{-1}(1)=[\varepsilon,\infty).\leqno{(5.5)}$$
Assume that the following conditions hold:
$${\rm (i)}\quad\widehat{\rho}(\widehat x_0)\geq 1,\qquad\,\,{\rm (ii)}\quad 
b^2\left(\frac{1}{(1-\alpha)\omega_n}\int_{\widehat{\overline M}}\widehat{\rho}\,dv_{\widehat{\overline g}}
\right)^{\frac{2}{n}}\leq 1.$$
Set 
$$s_{\widehat{\rho}}:=\left\{\begin{array}{ll}
\displaystyle{\frac{1}{b}\cdot\arcsin
\left\{b\left(\frac{1}{(1-\alpha)\omega_n}\int_{\widehat{\overline M}}\rho\,dv_{\widehat{\overline g}}
\right)^{\frac{1}{n}}\right\}} & (b\,:\,{\rm real})\\
\displaystyle{\left(\frac{1}{(1-\alpha)\omega_n}\int_{\widehat{\overline M}}\widehat{\rho}\,
dv_{\widehat{\overline g}}\right)^{\frac{1}{n}}} & (b\,:\,{\rm purely}\,\,{\rm imaginary}).
\end{array}\right.$$
which is defined by the above condition (ii).  
Furthermore, assume the following:
$$\hat{\alpha}s_{\rho}\leq\frac{\pi}{b}.\leqno{({\rm iii})}$$
Then there exists $s_1\in(0,s_{\rho})$ satisfying 
$$\overline{\phi}_{\widehat{\rho},\widehat x_0}(\hat{\alpha}s_1)\leq\alpha^{-1}\cdot\hat{\alpha}^{n-2}\cdot 
s_{\rho}\cdot\lim_{\varepsilon\to 0}\overline{\eta}_{\widehat{\rho},\widehat x_0,\lambda_{\varepsilon}}(s_1).$$
}

\vspace{0.5truecm}

By using Lemma 5.6, we prove Theorem 5.2.  

\vspace{0.5truecm}

\noindent
{\it Proof of Theorem 5.2.} 
We shall prove the statement in the case where $b$ is real (similar also the case where $b$ is purely imaginary).  
Let $\alpha,\,\hat{\alpha}$ be constants with $0<\alpha<1\leq\hat{\alpha}$ and 
$\lambda_{\varepsilon}\,\,\,(\varepsilon>0)$ be $C^1$-funnction over $\mathbb R$ satisfying $(5.4)$ and $(5.5)$.  
Define a function $\overline{\lambda}_{\varepsilon}\,\,\,(\varepsilon>0)$ over $\mathbb R$ by 
$\overline{\lambda}_{\varepsilon}(s):=\lambda_{\varepsilon}(s+\varepsilon)$ and define 
a function $\widehat{\rho}_{\varepsilon,t}\,\,\,(\varepsilon>0)$ over $\widehat{\overline M}$ by 
$\widehat{\rho}_{\varepsilon,t}:=\overline{\lambda}_{\varepsilon}(\widehat{\rho}-t)$.  
Since $\rho$ satisfies the condition $(\ast_2)$ in Theorem 5.2, $\widehat{\rho}_{\varepsilon,t}$ satisfies 
the conditions (ii) and (iii) in Lemma 5.6.  
By using Lemma 5.6 and discussing as in the proof of Theorem 2.1 in \cite{HoSp}, we can derive 
$$\begin{array}{l}
\displaystyle{\left(\int_{\widehat{\overline M}}\widehat{\rho}^{\frac{n}{n-1}}\,dv_{\widehat{\overline g}}
\right)^{\frac{n-1}{n}}
\leq\frac{n}{n-1}\cdot\frac{\pi}{2}\cdot\alpha^{-1}\cdot\hat{\alpha}^{n-2}\cdot(1-\alpha)^{-\frac{1}{n}}
\cdot\omega_n^{-\frac{1}{n}}}\\
\hspace{3.75truecm}\displaystyle{\times\int_{\widehat{\overline M}}
\left(\|{\rm grad}_{\widehat{\overline g}}\widehat{\rho}\|+\widehat{\rho}\cdot\|H^s\|\right)\,
dv_{\widehat{\overline g}}.}
\end{array}$$
Clearly we have 
$$\int_{\widehat{\overline M}}\widehat{\rho}^{\frac{n}{n-1}}\,dv_{\widehat{\overline g}}
=l\cdot\int_{\overline M}\rho^{\frac{n}{n-1}}\,dv_{\overline g}$$ 
and 
$$\int_{\widehat{\overline M}}\left(\|{\rm grad}_{\widehat{\overline g}}\widehat{\rho}\|+\widehat{\rho}
\cdot\|H^s\|\right)\,dv_{\widehat{\overline g}}
=l\cdot\int_{\overline M}\left(\|{\rm grad}_{\overline g}\rho\|+\rho\cdot\|H^s\|\right)\,dv_{\overline g}.$$
Hence we obtain 
$$\begin{array}{l}
\displaystyle{\left(\int_{\overline M}\rho^{\frac{n}{n-1}}\,dv_{\overline g}\right)^{\frac{n-1}{n}}
\leq l^{\frac{1}{n}}\cdot\frac{n}{n-1}\cdot\frac{\pi}{2}\cdot\alpha^{-1}\cdot\hat{\alpha}^{n-2}\cdot
(1-\alpha)^{-\frac{1}{n}}\cdot\omega_n^{-\frac{1}{n}}}\\
\hspace{3.75truecm}\displaystyle{\times\int_{\overline M}\left(\|{\rm grad}_{\overline g}\rho\|
+\rho\cdot\|H^s\|\right)\,dv_{\overline g}.}
\end{array}$$
\qed

\section{Approach to horizontally totally umbilicity} 
In this section, we recall the preservability of horizontally strongly convexity along the mean curvature flow.  
Let $\mathcal G\curvearrowright V$ be an isometric almost free action with minimal regularizable 
orbit of a Hilbert Lie group $\mathcal G$ on a Hilbert space $V$ equipped with an inner 
product $\langle\,\,,\,\,\rangle$ and $\phi:V\to N:=V/\mathcal G$ the orbit map.  
Denote by $\widetilde{\nabla}$ the Riemannian connection of $V$.  
Set $n:={\rm dim}\,\,N-1$.  
Let $M(\subset V)$ be a $\mathcal G$-invariant hypersurface in $V$ such that $\phi(M)$ is compact.  
Let $f$ be an inclusion map of $M$ into $V$ and 
$f_t\,\,(0\leq t<T)$ the $\mathcal G$-invariant regularized mean curvature flow starting from $f$.  
We use the notations in Sections 4.  In the sequel, we omit the notation $f_{t\ast}$ for simplicity.  
As stated in Introduction, set 
$$L:=\mathop{\sup}_{u\in V}\mathop{\max}_{(X_1,\cdots,X_5)\in({\widetilde{\mathcal H}}_1)_u^5}
\vert\langle{\mathcal A}^{\phi}_{X_1}
((\widetilde{\nabla}_{X_2}{\mathcal A}^{\phi})_{X_3}X_4),\,X_5\rangle\vert,\leqno{(6.1)}$$
where $\widetilde{\mathcal H}_1:=\{X\in\widetilde{\mathcal H}\,\vert\,\,\,
\vert\vert X\vert\vert=1\}$.  
Assume that $L<\infty$.  Note that $L<\infty$ in the case where $N$ is compact.  
In \cite{Koi2}, we proved the following horizontally strongly convexity preservability theorem 
by using evolution equations stated in Section 4 and the discussion in the proof of Theorem 5.1.  

\vspace{0.5truecm}

\noindent
{\bf Theorem 6.1(\cite{Koi2}).} {\sl 
If $M$ satisfies $(H^s_0)^2(h_{\mathcal H})_{(\cdot,0)}
>2n^2L(g_{\mathcal H})_{(\cdot,0)}$, then $T<\infty$ holds and 
$(H^s_t)^2(h_{\mathcal H})_{(\cdot,t)}>2n^2L
(g_{\mathcal H})_{(\cdot,t)}$ holds for all $t\in[0,T)$.}

\vspace{0.5truecm}

In this section, we shall prove the following result for the approach to the horizontally totally 
umbilicity of $f_t$ as $t\to T$.  

\vspace{0.5truecm}

\noindent
{\bf Proposition 6.2.} {\sl 
Under the hypothesis of Theorem 6.1, there exist positive constants $\delta$ and $C_0$ depending on only 
$f,\,L,\,\overline K$ and the injective radius $i(N)$ of $N$ such that 
$$\vert\vert(A_{\mathcal H})_{(\cdot,t)}\vert\vert^2-\frac{(H^s_t)^2}{n}\,<\,
C_0(H^s_t)^{2-\delta}$$
holds for all $t\in[0,T)$.}

\vspace{0.5truecm}

We prepare some lemmas to show this proposition.  
In the sequel, we denote the fibre metric of ${\mathcal H}^{(r,s)}$ induced from 
$g_{\mathcal H}$ by the same symbol $g_{\mathcal H}$, and set $\vert\vert S\vert\vert:=\sqrt{g_{\mathcal H}(S,S)}$ 
for $S\in\Gamma({\mathcal H}^{(r,s)})$.  
Define a function $\psi_{\delta}$ over $M$ by 
$$\psi_{\delta}:=\frac{1}{(H^s)^{2-\delta}}
\left(\vert\vert A_{\mathcal H}\vert\vert^2-\frac{(H^s)^2}{n}\right).$$

\vspace{0.5truecm}

\noindent
{\bf Lemma 6.2.1.} {\sl Set $\alpha:=2-\delta$.  Then we have 
$$\begin{array}{l}
\displaystyle{\frac{\partial\psi_{\delta}}{\partial t}=\triangle_{\mathcal H}
\psi_{\delta}+(2-\alpha)\vert\vert A_{\mathcal H}\vert\vert^2\psi_{\delta}
+\frac{(\alpha-1)(\alpha-2)}{(H^s)^2}
\vert\vert\,dH^s\,
\vert\vert^2\psi_{\delta}}\\
\hspace{1.2truecm}\displaystyle{+\frac{2(\alpha-1)}{H^s}
g_{\mathcal H}(dH^s,d\psi_{\delta})
-\frac{2}{(H^s)^{\alpha+2}}\left\vert\left\vert\,\vert\vert H
\vert\vert\nabla^{\mathcal H}A_{\mathcal H}-dH^s\otimes A_{\mathcal H}\,\right\vert\right\vert^2}\\
\hspace{1.2truecm}\displaystyle{+3(\alpha-2){\rm Tr}(({\mathcal A}^{\phi}_{\xi})^2)
_{\mathcal H})\psi_{\delta}-\frac{6}{(H^s)^{\alpha-1}}
{\rm Tr}(({\mathcal A}^{\phi}_{\xi})^2\circ A_{\mathcal H})}\\
\hspace{1.2truecm}\displaystyle{
+\frac{6}{n(H^s)^{\alpha-2}}{\rm Tr}
(({\mathcal A}^{\phi}_{\xi})^2)-\frac{4}{(H^s)^{\alpha}}
{\rm Tr}^{\cdot}_{g_{\mathcal H}}{\rm Tr}^{\bullet}_{g_{\mathcal H}}
h(((\nabla_{\bullet}{\mathcal A})_{\bullet}\circ A_{\mathcal H})\cdot,\cdot)}\\
\hspace{1.2truecm}\displaystyle{
-\frac{4}{(H^s)^{\alpha}}
{\rm Tr}^{\cdot}_{g_{\mathcal H}}{\rm Tr}^{\bullet}_{g_{\mathcal H}}
h(({\mathcal A}_{\bullet}\circ A_{\mathcal H})\cdot,{\mathcal A}_{\bullet}\cdot)
+\frac{4}{(H^s)^{\alpha}}
{\rm Tr}^{\cdot}_{g_{\mathcal H}}{\rm Tr}^{\bullet}_{g_{\mathcal H}}
h(({\mathcal A}_{\bullet}\circ{\mathcal A}_{\bullet})\cdot,A_{\mathcal H}\cdot).}
\end{array}$$
}

\vspace{0.5truecm}

\noindent
{\it Proof.} By using Lemmas 4.9 and 4.11, we have 
$$\begin{array}{l}
\displaystyle{\frac{\partial\psi_{\delta}}{\partial t}
=(2-\alpha)\vert\vert A_{\mathcal H}\vert\vert^2\psi_{\delta}
+\frac{1}{(H^s)^{\alpha}}
\triangle_{\mathcal H}(\vert\vert A_{\mathcal H}\vert\vert^2)}\\
\hspace{1.2truecm}\displaystyle{
-\frac{1}{(H^s)^{\alpha+1}}\left(
\alpha\vert\vert A_{\mathcal H}\vert\vert^2-\frac{(\alpha-2)(H^s)^2}{n}
\right)\triangle_{\mathcal H}H^s}\\
\hspace{1.2truecm}\displaystyle{-\frac{2}{(H^s)^{\alpha}}
\vert\vert\nabla^{\mathcal H}A_{\mathcal H}\vert\vert^2+(3\alpha-2){\rm Tr}
(({\mathcal A}^{\phi}_{\xi})^2)_{\mathcal H}\cdot\psi_{\delta}}\\
\hspace{1.2truecm}\displaystyle{-\frac{6}{(H^s)^{\alpha-1}}
{\rm Tr}\left(({\mathcal A}^{\phi}_{\xi})^2\circ(A_{\mathcal H}
-\frac{H^s}{n}\,{\rm id})\right)}\\
\hspace{1.2truecm}\displaystyle{-\frac{4}{(H^s)^{\alpha}}
{\rm Tr}^{\cdot}_{g_{\mathcal H}}{\rm Tr}^{\bullet}_{g_{\mathcal H}}
h(({\mathcal A}_{\bullet}\circ A_{\mathcal H})\cdot,{\mathcal A}_{\bullet}\cdot)}\\
\hspace{1.2truecm}\displaystyle{+\frac{4}{(H^s)^{\alpha}}
{\rm Tr}^{\cdot}_{g_{\mathcal H}}{\rm Tr}^{\bullet}_{g_{\mathcal H}}
h(({\mathcal A}_{\bullet}\circ{\mathcal A}_{\bullet})\cdot,A_{\mathcal H}\cdot)}\\
\hspace{1.2truecm}\displaystyle{-\frac{4}{(H^s)^{\alpha}}
{\rm Tr}^{\cdot}_{g_{\mathcal H}}{\rm Tr}^{\bullet}_{g_{\mathcal H}}
h((\nabla_{\bullet}{\mathcal A})_{\bullet}\circ A_{\mathcal H})\cdot,\cdot).}
\end{array}\leqno{(6.2)}$$
Also we have 
$$\begin{array}{l}
\displaystyle{\triangle_{\mathcal H}\psi_{\delta}
=\frac{1}{(H^s)^{\alpha}}\triangle_{\mathcal H}
\vert\vert A_{\mathcal H}\vert\vert^2-\frac{2\alpha}{(H^s)^{\alpha+1}}
g_{\mathcal H}(dH^s,d(\vert\vert A_{\mathcal H}\vert\vert^2))}\\
\hspace{1truecm}\displaystyle{
-\frac{1}{(H^s)^{\alpha+1}}\left(
\alpha\vert\vert A_{\mathcal H}\vert\vert^2-\frac{(\alpha-2)(H^s)^2}{n}
\right)\triangle_{\mathcal H}H^s}\\
\hspace{1truecm}\displaystyle{
+\frac{1}{(H^s)^{\alpha+2}}
\left(\alpha(\alpha+1)\vert\vert A_{\mathcal H}\vert\vert^2
-\frac{(\alpha-1)(\alpha-2)(H^s)^2}{n}
\right)\left\vert\left\vert dH^s\,\right\vert\right\vert^2}
\end{array}\leqno{(6.3)}$$
From $(6.2)$ and $(6.3)$, we obtain the desired relation.  \qed

\vspace{0.5truecm}

Then we have the following inequalities.  

\vspace{0.5truecm}

By using the Codazzi equation, we can derive the following relation.  

\vspace{0.5truecm}

\noindent
{\bf Lemma 6.2.2.} {\sl For any $X,Y,Z\in{\mathcal H}$, we have 
$$(\nabla^{\mathcal H}_Xh_{\mathcal H})(Y,Z)=(\nabla^{\mathcal H}_Yh_{\mathcal H})(X,Z)
+2h({\mathcal A}_XY,Z)-h({\mathcal A}_YZ,X)+h({\mathcal A}_XZ,Y)$$
or equivalently, 
$$(\nabla^{\mathcal H}_XA_{\mathcal H})(Y)=(\nabla^{\mathcal H}_YA_{\mathcal H})(X)
+2(A\circ{\mathcal A}_X)Y+({\mathcal A}_Y\circ A)(X)-({\mathcal A}_X\circ A)(Y).$$
}

\vspace{0.5truecm}

\noindent
{\it Proof.} Let $(x,t)$ be the base point of $X,Y$ and $Z$ and extend these vectors to sections 
$\widetilde X,\widetilde Y$ and $\widetilde Z$ of ${\mathcal H}_t$ with 
$(\nabla^{\mathcal H}\widetilde X)_{(x,t)}=(\nabla^{\mathcal H}\widetilde Y)_{(x,t)}
=(\nabla^{\mathcal H}\widetilde Z)_{(x,t)}=0$.  
Since $\nabla h$ is symmetric with respect to $g$ by the Codazzi equation and the flatness of $V$, 
we have 
$$\begin{array}{l}
\displaystyle{(\nabla^{\mathcal H}_Xh_{\mathcal H})(Y,Z)=X(h(\widetilde Y,\widetilde Z))}\\
\displaystyle{=(\nabla_Xh)(Y,Z)+h({\mathcal A}_XY,Z)+h({\mathcal A}_XZ,Y)}\\
\displaystyle{=(\nabla_Yh)(X,Z)+h({\mathcal A}_XY,Z)+h({\mathcal A}_XZ,Y)}\\
\displaystyle{=Y(h(\widetilde X,\widetilde Z))-h({\mathcal A}_YX,Z)-h({\mathcal A}_YZ,X)
+h({\mathcal A}_XY,Z)+h({\mathcal A}_XZ,Y)}\\
\displaystyle{=(\nabla^{\mathcal H}_Yh)(X,Z)+2h({\mathcal A}_XY,Z)-h({\mathcal A}_YZ,X)
+h({\mathcal A}_XZ,Y).}
\end{array}$$
\qed

\vspace{0.5truecm}

Set 
$$K:=\max_{(e_1,e_2):{\rm o.n.s.}\,\,{\rm of}\,\,TV}\vert\vert{\mathcal A}^{\phi}_{e_1}e_2\vert\vert^2,$$
where "${\rm o.n.s.}$" means "orthonormal system".  
Assume that $K<\infty$.  Note that $K<\infty$ if $N=V/\mathcal G$ is compact.  
For a section $S$ of 
${\mathcal H}^{(r,s)}$ 
and a permutaion $\sigma$ of $s$-symbols, we define a section $S_{\sigma}$ of 
${\mathcal H}^{(r,s)}$ by 
$$S_{\sigma}(X_1,\cdots,X_s):=S(X_{\sigma(1)},\cdots,X_{\sigma(s)})\,\,\,\,
(X_1,\cdots,X_s\in{\mathcal H})$$
and ${\rm Alt}(S)$ by 
$${\rm Alt}(S):=\frac{1}{s!}\sum_{\sigma}{\rm sgn}\,\sigma\,S_{\sigma},$$
where $\sigma$ runs over the symmetric group of degree $s$.  
Also, denote by $(i,j)$ the transposition exchanging $i$ and $j$.  
Since $(H^s_t)^2(h_{\mathcal H})_{(\cdot,t)}
>n^2L(g_{\mathcal H})_{(\cdot,t)}\,\,(t\in[0,T))$
and $\phi(M)$ is compact, there exists a positive constant $\varepsilon$ satisfying 
$$\vert\vert H_{(\cdot,0)}\vert\vert^2(h_{\mathcal H})_{(\cdot,0)}
\geq n^2L(g_{\mathcal H})_{(\cdot,0)}+\varepsilon\vert\vert H_{(\cdot,0)}\vert\vert^3
(g_{\mathcal H})_{(\cdot,0)}.\leqno{(\sharp)}$$
Then we can show that 
$$\vert\vert H_{(\cdot,t)}\vert\vert^2(h_{\mathcal H})_{(\cdot,t)}
\geq n^2L(g_{\mathcal H})_{(\cdot,t)}+\varepsilon\vert\vert H_{(\cdot,t)}\vert\vert^3
(g_{\mathcal H})_{(\cdot,t)}$$
holds for all $t\in[0,T)$.  
Without loss of generality, we may assume that $\varepsilon\leq 1$.  
Then we have the following inequalities.  

\vspace{0.5truecm}

\noindent
{\bf Lemma 6.2.3.} {\sl Let $\varepsilon$ be as above.  Then we have 
the following inequalities:
$$H^s{\rm Tr}_{\mathcal H}(A_{\mathcal H})^3
-\vert\vert(A_{\mathcal H})_t\vert\vert^4\geq n\varepsilon^2(H^s)^2
\left(\vert\vert(A_{\mathcal H})_t\vert\vert^2-\frac{(H^s)^2}{n}
\right),\leqno{(6.4)}$$
and 
$$\begin{array}{l}
\displaystyle{
\left\vert\left\vert\,\,
H^s\nabla^{\mathcal H}A_{\mathcal H}
-dH^s\otimes A_{\mathcal H}\,\,\right\vert\right\vert^2}\\
\displaystyle{
\geq -8\varepsilon^{-2}K(H^s_{(u,t)})^2
+\frac{1}{8}\vert\vert(dH^s)_{(u,t)}\vert\vert^2
\varepsilon^2(H^s_{(u,t)})^2.}
\end{array}
\leqno{(6.5)}$$
}

\noindent
{\it Proof.} First we shall show the inequality $(6.4)$.  
Fix $(u,t)\in M\times[0,T)$.  
Take an orthonormal base $\{e_1,\cdots,e_n\}$ of ${\mathcal H}_{(u,t)}$ with 
respect to $g_{(u,t)}$ consisting of the eigenvectors of 
$(A_{\mathcal H})_{(u,t)}$.  Let $(A_{\mathcal H})_{(u,t)}(e_i)=\lambda_ie_i\,\,
(i=1,\cdots,n)$.  Note that 
$\lambda_i>\varepsilon H^s(>0)\,\,(i=1,\cdots,n)$.  
Then we have 
$$H^s{\rm Tr}_{\mathcal H}(A_{\mathcal H})^3
-\vert\vert(A_{\mathcal H})_t\vert\vert^4
=\sum_{1\leq i<j\leq n}\lambda_i\lambda_j(\lambda_i-\lambda_j)^2
>\varepsilon^2(H^s)^2\sum_{1\leq i<j\leq n}
(\lambda_i-\lambda_j)^2.$$
On the other hand, we have 
$$\vert\vert(A_{\mathcal H})_t\vert\vert^2-\frac{(H^s)^2}{n}
=\frac1n\sum_{1\leq i<j\leq n}(\lambda_i-\lambda_j)^2.$$
From these inequalities, we can derive the inequality $(6.4)$.  

Next we shall show the inequality $(6.5)$.  
By using Lemma 6.2.2, we can show 
$$\begin{array}{l}
\displaystyle{
\left\vert\left\vert\,\,H^s\nabla^{\mathcal H}A_{\mathcal H}
-dH^s\otimes A_{\mathcal H}\,\,\right\vert\right\vert^2
\geq\left\vert\left\vert{\rm Alt}\left(H^s\nabla^{\mathcal H}A_{\mathcal H}
-dH^s\otimes A_{\mathcal H}\right)\right\vert\right\vert^2}\\
\hspace{0truecm}\displaystyle{
\geq\left\vert\left\vert A\circ{\mathcal A}-\frac{1}{2}{\mathcal A}\circ(A\times{\rm id})
-\frac{1}{2}{\mathcal A}\circ({\rm id}\times A)
-{\rm Alt}\left(dH^s\otimes A_{\mathcal H}\right)\right\vert\right\vert^2.}
\end{array}$$
For simplicity, we set 
$$S:=A\circ{\mathcal A}-\frac{1}{2}{\mathcal A}\circ(A\times{\rm id})
-\frac{1}{2}{\mathcal A}\circ({\rm id}\times A).$$
It is clear that $(6.5)$ holds at $(u,t)$ if $(dH^s)_{(u,t)}=0$.  
Assume that $(dH^s)_{(u,t)}\not=0$.  
Take an orthonormal base $(e_1,\cdots,e_n)$ of ${\mathcal H}_{(u,t)}$ with respect to 
$(g_{\mathcal H})_{(u,t)}$ with $e_1=\frac{(dH^s)_{(u,t)}}
{\vert\vert(dH^s)_{(u,t)}\vert\vert}$.  
Then we have 
$$\begin{array}{l}
\displaystyle{
\left\vert\left\vert S_{(u,t)}-{\rm Alt}\left(dH^s\otimes A_{\mathcal H}\right)_{(u,t)}
\right\vert\right\vert^2}\\
\displaystyle{\geq
\left\vert\left\vert S-{\rm Alt}\left(dH^s\otimes A_{\mathcal H}\right)(e_1,e_2)
\right\vert\right\vert^2}\\
\displaystyle{\geq\vert\vert S(e_1,e_2)\vert\vert^2
-\vert\vert(dH^s)_{(u,t)}\vert\vert g(S(e_1,e_2),A_{\mathcal H}e_2)
+\frac{1}{4}\vert\vert(dH^s)_{(u,t)}\vert\vert^2\cdot\vert\vert A_{\mathcal H}e_2\vert\vert^2}\\
\displaystyle{\geq\vert\vert S(e_1,e_2)\vert\vert^2-\vert\vert(dH^s)_{(u,t)}\vert\vert 
g(S(e_1,e_2),A_{\mathcal H}e_2)
+\frac{1}{4}\vert\vert(dH^s)_{(u,t)}\vert\vert^2
\varepsilon^2(H^s_{(u,t)})^2}\\
\displaystyle{\geq
(1-2\varepsilon^{-2})\vert\vert S(e_1,e_2)\vert\vert^2
+\left(\sqrt 2\varepsilon^{-1}\vert\vert S(e_1,e_2)\vert\vert
-\frac{1}{2\sqrt 2}\vert\vert(dH^s)_{(u,t)}\vert\vert\varepsilon H^s_{(u,t)}\right)^2}\\
\hspace{0.5truecm}\displaystyle{+\frac{1}{8}\vert\vert(dH^s)_{(u,t)}\vert\vert^2
\varepsilon^2(H^s_{(u,t)})^2}\\
\displaystyle{\geq
-2\varepsilon^{-2}\vert\vert S(e_1,e_2)\vert\vert^2
+\frac{1}{8}\vert\vert(dH^s)_{(u,t)}\vert\vert^2
\varepsilon^2(H^s_{(u,t)})^2}\\
\displaystyle{\geq
-8\varepsilon^{-2}K(H^s_{(u,t)})^2
+\frac{1}{8}\vert\vert(dH^s)_{(u,t)}\vert\vert^2
\varepsilon^2(H^s_{(u,t)})^2,}
\end{array}$$
where we use $\vert\vert A_{\mathcal H}e\vert\vert\leq H^s$ holds for any unit vector $e$ of 
$\mathcal H$.  Thus we see that $(6.5)$ holds at $(u,t)$.  
This completes the proof.  \qed

\vspace{0.5truecm}

From Lemma 6.2.1 and (6.5), we obtain the following lemma.  

\vspace{0.5truecm}

\noindent
{\bf Lemma 6.2.4.} {\sl Assume that $\delta<1$.  Then we have the following inequality:
$$\begin{array}{l}
\displaystyle{\frac{\partial\psi_{\delta}}{\partial t}\leq\triangle_{\mathcal H}
\psi_{\delta}+(2-\alpha)\vert\vert A_{\mathcal H}\vert\vert^2\psi_{\delta}
+\frac{2(\alpha-1)}{H^s}g_{\mathcal H}(dH^s,d\psi_{\delta})}\\
\hspace{1.2truecm}\displaystyle{-\frac{2}{(H^s)^{\alpha+2}}
\left(\frac{1}{8}\vert\vert\,dH^s\,\vert\vert^2
\varepsilon^2(H^s)^2-8\varepsilon^{-2}K(H^s)^2\right)}\\
\hspace{1.2truecm}\displaystyle{+3(\alpha-2){\rm Tr}(({\mathcal A}^{\phi}_{\xi})^2)_{\mathcal H})\psi_{\delta}
-\frac{6}{(H^s)^{\alpha-2}}
{\rm Tr}(({\mathcal A}^{\phi}_{\xi})^2\circ A_{\mathcal H})}\\
\hspace{1.2truecm}\displaystyle{
+\frac{6}{n(H^s)^{\alpha-2}}{\rm Tr}
(({\mathcal A}^{\phi}_{\xi})^2)-\frac{4}{(H^s)^{\alpha}}
{\rm Tr}^{\cdot}_{g_{\mathcal H}}{\rm Tr}^{\bullet}_{g_{\mathcal H}}
h(((\nabla_{\bullet}{\mathcal A})_{\bullet}\circ A_{\mathcal H})\cdot,\cdot)}\\
\hspace{1.2truecm}\displaystyle{
-\frac{4}{(H^s)^{\alpha}}
{\rm Tr}^{\cdot}_{g_{\mathcal H}}{\rm Tr}^{\bullet}_{g_{\mathcal H}}
h(({\mathcal A}_{\bullet}\circ A_{\mathcal H})\cdot,{\mathcal A}_{\bullet}\cdot)}\\
\hspace{1.2truecm}\displaystyle{
+\frac{4}{(H^s)^{\alpha}}
{\rm Tr}^{\cdot}_{g_{\mathcal H}}{\rm Tr}^{\bullet}_{g_{\mathcal H}}
h(({\mathcal A}_{\bullet}\circ{\mathcal A}_{\bullet})\cdot,A_{\mathcal H}\cdot).}
\end{array}$$
}

\vspace{0.5truecm}

On the other hand, we can show the following fact for $\psi_{\delta}$.  

\vspace{0.5truecm}

\noindent
{\bf Lemma 6.2.5.} {\sl We have 
$$\begin{array}{l}
\displaystyle{\triangle_{\mathcal H}\psi_{\delta}
=\frac{2}{(H^s)^{\alpha+2}}\times
\left\vert\left\vert\,\,
H^s\cdot\nabla^{\mathcal H}A_{\mathcal H}-dH^s\cdot A_{\mathcal H}\,\,
\right\vert\right\vert^2}\\
\hspace{1.5truecm}\displaystyle{+\frac{2}{(H^s)^{\alpha-1}}\left(
{\rm Tr}((A_{\mathcal H})^3)-{\rm Tr}(({\mathcal A}^{\phi}_{\xi})^2\circ A_{\mathcal H})\right)}\\
\hspace{1.5truecm}\displaystyle{-\frac{2}{(H^s)^{\alpha}}\left(
{\rm Tr}((A_{\mathcal H})^2-({\mathcal A}^{\phi}_{\xi})^2\vert_{\mathcal H}\right)\vert\vert A_{\mathcal H}\vert\vert^2}\\
\hspace{1.5truecm}\displaystyle{+\frac{2}{(H^s)^{\alpha}}
{\rm Tr}^{\bullet}_{g_{\mathcal H}}
\left((\nabla^{\mathcal H}dH^s)
\left(\left(A_{\mathcal H}-\frac{H^s}{n}\,{\rm id}\right)(\bullet),\bullet\right)\right)}\\
\hspace{1.5truecm}\displaystyle{-\frac{\alpha}{H^s}\psi_{\delta}
\triangle_{\mathcal H}H^s-\frac{(\alpha-1)(\alpha-2)}{(H^s)^2}
\vert\vert\,dH^s\,\,\vert\vert^2\psi_{\delta}}\\
\hspace{1.5truecm}\displaystyle{-\frac{2(\alpha-1)}{H^s}g_{\mathcal H}
(dH^s,d\psi_{\delta})+\frac{2}{(H^s)^{\alpha}}\times
{\rm Tr}_{g_{\mathcal H}}^{\bullet}{\mathcal R}(A_{\mathcal H}\bullet,\bullet).}
\end{array}$$
}

\vspace{0.5truecm}

\noindent
{\it Proof.} According to $(4.16)$ in \cite{Koi2}, we have 
$${\rm Tr}^{\bullet}_{g_{\mathcal H}}(\triangle_{\mathcal H}^{\mathcal H}h_{\mathcal H})(A_{\mathcal H}\bullet,
\bullet)
=\frac12\triangle_{\mathcal H}\vert\vert A_{\mathcal H}\vert\vert^2
-\vert\vert\nabla^{\mathcal H}A_{\mathcal H}\vert\vert^2.\leqno{(6.6)}$$
Also we have 
$$(A^2)_{\mathcal H}=(A_{\mathcal H})^2-({\mathcal A}^{\phi}_{\xi})^2.\leqno{(6.7)}$$
By using Lemmas 4.4, 4.5 and these relations, we can derive 
$$\begin{array}{l}
\displaystyle{\frac12\triangle_{\mathcal H}\vert\vert A_{\mathcal H}\vert\vert^2
={\rm Tr}^{\bullet}_{g_{\mathcal H}}(\nabla^{\mathcal H}dH^s)(A_{\mathcal H}\bullet,\bullet)
+H^s{\rm Tr}((A_{\mathcal H})^3)}\\
\displaystyle{-H^s{\rm Tr}(({\mathcal A}^{\phi}_{\xi})^2\circ A_{\mathcal H})
-{\rm Tr}\left((A_{\mathcal H})^2-({\mathcal A}^{\phi}_{\xi})^2\vert_{\mathcal H}\right)
\vert\vert A_{\mathcal H}\vert\vert^2}\\
\displaystyle{+{\rm Tr}_{g_{\mathcal H}}^{\bullet}{\mathcal R}(A_{\mathcal H}\bullet,\bullet)
+\vert\vert\nabla^{\mathcal H}A_{\mathcal H}\vert\vert^2.}
\end{array}\leqno{(6.8)}$$
By substituting this relation into $(6.3)$, we obtain 
$$\begin{array}{l}
\displaystyle{\triangle_{\mathcal H}\psi_{\delta}
=\frac{2}{(H^s)^{\alpha}}\times
\left\{
{\rm Tr}^{\bullet}_{g_{\mathcal H}}(\nabla^{\mathcal H}dH^s)(A_{\mathcal H}\bullet,\bullet)
+H^s{\rm Tr}((A_{\mathcal H})^3)\right.}\\
\hspace{2.8truecm}\displaystyle{-H^s{\rm Tr}(({\mathcal A}^{\phi}_{\xi})^2\circ A_{\mathcal H})
-{\rm Tr}\left((A_{\mathcal H})^2-({\mathcal A}^{\phi}_{\xi})^2\vert_{\mathcal H}\right)
\vert\vert A_{\mathcal H}\vert\vert^2}\\
\hspace{2.8truecm}\displaystyle{+{\rm Tr}_{g_{\mathcal H}}^{\bullet}{\mathcal R}(A_{\mathcal H}\bullet,\bullet)
+\vert\vert\nabla^{\mathcal H}A_{\mathcal H}\vert\vert^2\}}\\
\hspace{1.5truecm}\displaystyle{-\frac{2\alpha}{(H^s)^{\alpha+1}}
g_{\mathcal H}(dH^s,d(\vert\vert A_{\mathcal H}\vert\vert^2))}\\
\hspace{1.5truecm}\displaystyle{
-\frac{1}{(H^s)^{\alpha+1}}\left(
\alpha\vert\vert A_{\mathcal H}\vert\vert^2-\frac{(\alpha-2)(H^s)^2}{n}
\right)\triangle_{\mathcal H}H^s}\\
\hspace{1.5truecm}\displaystyle{
+\frac{1}{(H^s)^{\alpha+2}}
\left(\alpha(\alpha+1)\vert\vert A_{\mathcal H}\vert\vert^2
-\frac{(\alpha-1)(\alpha-2)(H^s)^2}{n}
\right)\left\vert\left\vert dH^s\,\right\vert\right\vert^2.}
\end{array}$$
From this relation, we can derive the desired relation.  \qed

\vspace{0.5truecm}

From this lemma, we can derive the following inequality for $\psi_{\delta}$ directly.  

\vspace{0.5truecm}

\noindent
{\bf Lemma 6.2.6.} {\sl We have 
$$\begin{array}{l}
\displaystyle{\triangle_{\mathcal H}\psi_{\delta}
\geq\frac{2}{(H^s)^{\alpha-1}}\left(
{\rm Tr}((A_{\mathcal H})^3)-{\rm Tr}(({\mathcal A}^{\phi}_{\xi})^2\circ A_{\mathcal H})\right)}\\
\hspace{1.5truecm}\displaystyle{-\frac{2}{(H^s)^{\alpha}}\left(
\vert\vert A_{\mathcal H}\vert\vert^2-{\rm Tr}({\mathcal A}^{\phi}_{\xi})^2\vert_{\mathcal H}\right)
{\rm Tr}((A_{\mathcal H})^2)}\\
\hspace{1.5truecm}\displaystyle{+\frac{2}{(H^s)^{\alpha}}
{\rm Tr}^{\bullet}_{g_{\mathcal H}}\left(
(\nabla^{\mathcal H}dH^s)\left(\left(A_{\mathcal H}-\frac{H^s}{n}\,{\rm id}
\right)(\bullet),\bullet\right)\right)}\\
\hspace{1.5truecm}\displaystyle{
-\frac{\alpha}{H^s}\psi_{\delta}\triangle_{\mathcal H}H^s
-\frac{2(\alpha-1)}{H^s}g_{\mathcal H}(dH^s,d\psi_{\delta})}\\
\hspace{1.5truecm}\displaystyle{+\frac{2}{(H^s)^{\alpha}}\times
{\rm Tr}_{g_{\mathcal H}}^{\bullet}{\mathcal R}(A_{\mathcal H}\bullet,\bullet).}
\end{array}$$
}

\vspace{0.5truecm}

For a function $\rho$ over $M\times[0,T)$ such that $\rho(\cdot,t)$ ($t\in[0,T))$ are $\mathcal G$-invariant, 
define a function $\rho_B$ over $\overline M\times[0,T)$ by 
$\rho_B\circ(\phi_M\times{\rm id}_{[0,T)})=\rho$.  We call this function the {\it function over} 
$\overline M\times[0,T)$ {\it associated with} $\rho$.  Denote by $g_N$ the Riemannian orbimetric of $N$ and set 
$\bar g_t:=\bar f_t^{\ast}g_N$.  Also, denote by $d\bar v_t$ the orbivolume element of $\bar g_t$.  
Define a section $\bar g$ of $\pi_{\overline M}^{\ast}(T^{(0,2)}\overline M)$ by 
$\bar g(x,t)=(g_t)_x\,\,\,((x,t)\in\overline M\times[0,T))$ and a section $d\bar v$ of 
$\pi_{\overline M}^{\ast}(\wedge^nT^{\ast}\overline M)$ by 
$d\bar v(x,t)=(d\bar v_t)_x\,\,\,((x,t)\in\overline M\times[0,T))$, where $\pi_{\overline M}$ is the natural 
projection of $\overline M\times[0,T)$ onto $\overline M$ and $\pi_{\overline M}^{\ast}(\bullet)$ denotes the 
induced bundle of $(\bullet)$ by $\pi_{\overline M}$.  
Denote by $\overline{\nabla}^t$ the Riemannian orbiconnection of $\bar g_t$ and by $\overline{\triangle}_t$ 
the Laplace operator of $\overline{\nabla}^t$.  Define an orbiconnection $\overline{\nabla}$ of 
$\pi^{\ast}_{\overline M}(T\overline M)$ by using $\overline{\nabla}^t$'s (see the definition of $\nabla$ in 
Section 4).  Also, let $\overline{\triangle}$ be the differential operator of $\pi_{\overline M}^{\ast}
(\overline M\times{\Bbb R})$ defined by using $\overline{\triangle}_t$'s.  
Denote by $\displaystyle{\int_{\overline M}\rho_B\,d\bar v}$ the function over $[0,T)$ defined by assigning 
$\displaystyle{\int_{\overline M}\rho_B(\cdot,t)\,d\bar v_t}$ to each $t\in[0,T)$.  

\vspace{0.5truecm}

Let $\rho$ and $\rho_B$ be as above.  According to Theorem 5.1, we have 
$$\int_{\overline M}({\rm div}_{\nabla^{\mathcal H}}\rho)_Bd\bar v
=\int_{\overline M}{\rm div}_{\overline{\nabla}}(\rho_B)d\bar v=0\leqno{(6.9)}$$ 
and 
$$\int_{\overline M}({\triangle}_{\mathcal H}\rho)_Bd\bar v
=\int_{\overline M}\overline{\triangle}(\rho_B)d\bar v=0.\leqno{(6.10)}$$

From the inequlaity in Lemma 6.2.6 and $(6.9)$, we can derive the following integral inequality.  

\vspace{0.5truecm}

\noindent
{\bf Lemma 6.2.7.} {\sl Assume that $0\leq\delta\leq\frac{1}{2}$.  
Then, for any $\beta\geq 2$, we have 
$$\begin{array}{l}
\hspace{0.5truecm}\displaystyle{n\varepsilon^2\int_{\overline M}
(H^s_B)^2(\psi_{\delta})_B^{\beta}d\bar v}\\
\displaystyle{\leq\frac{3\beta\eta+6}{2}\int_{\overline M}
(H^s_B)^{-\alpha}(\psi_{\delta})_B^{\beta-1}
\vert\vert\,\,dH^s\,\,\vert\vert_B^2d\bar v+\frac{3\beta}{2\eta}\int_{\overline M}
(\psi_{\delta})_B^{\beta-2}\vert\vert d\psi_{\delta}\vert\vert_B^2d\bar v}\\
\hspace{0.5truecm}\displaystyle{+C_1\int_{\overline M}
(H^s_B)^{-\alpha}(\psi_{\delta})_B^{\beta-1}
\vert\vert A_{\mathcal H}\vert\vert^2_Bd\bar v+C_2\int_{\overline M}
(H^s_B)^{-\alpha}(\psi_{\delta})_B^{\beta-1}d\bar v,}\\
\end{array}$$
where $C_i$ ($i=1,2$) are positive constants depending only on $K$ and $L$ 
($L$ is the constant defined in the previous section).}

\vspace{0.5truecm}

\noindent
{\it Proof.} By using $\displaystyle{\int_{\overline M}{\rm div}_{\nabla^{\mathcal H}}
\left((H^s)^{-\alpha}(\psi_{\delta})^{\beta-1}(A_{\mathcal H}-(H^s/n){\rm id})
({\rm grad}\,H^s)\right)_Bd\bar v=0}$ and Lemma 6.2.2, we can show 
$$\begin{array}{l}
\displaystyle{\int_{\overline M}(H^s_B)^{-\alpha}(\psi_{\delta})_B^{\beta-1}
\left({\rm Tr}^{\bullet}_{g_{\mathcal H}}\left(
(\nabla^{\mathcal H}dH^s)\left(\left(A_{\mathcal H}-(H^s/n)\,{\rm id}
\right)(\bullet),\bullet\right)\right)\right)_Bd\bar v}\\
\displaystyle{=\alpha\int_{\overline M}(H^s_B)^{-\alpha-1}(\psi_{\delta})_B^{\beta-1}
g_{\mathcal H}((dH^s\otimes dH^s,h_{\mathcal H}-(H^s/n)
g_{\mathcal H})_Bd\bar v}\\
\hspace{0.4truecm}\displaystyle{-(\beta-1)\int_{\overline M}
(H^s_B)^{-\alpha}(\psi_{\delta})_B^{\beta-2}
g_{\mathcal H}((dH^s\otimes d\psi_{\delta},h_{\mathcal H}-(H^s/n)
g_{\mathcal H})_Bd\bar v}\\
\hspace{0.4truecm}\displaystyle{-(1-1/n)\int_{\overline M}
(H^s_B)^{-\alpha}(\psi_{\delta})_B^{\beta-1}
\vert\vert\,dH^s\,\,\vert\vert_B^2d\bar v}\\
\hspace{0.4truecm}\displaystyle{+3\int_{\overline M}(H^s_B)^{-\alpha}
(\psi_{\delta})_B^{\beta-1}{\rm Tr}_{g_{\mathcal H}}
({\mathcal A}^{\phi}_{\xi}\circ{\mathcal A}^{\phi}_{{\rm grad}H^s})_Bd\bar v.}
\end{array}\leqno{(6.11)}$$
Also, by using $\int_{\overline M}(\triangle_{\mathcal H}\psi_{\delta}^{\beta})_Bd\bar v=0$, we can show 
$$\int_{\overline M}(\psi_{\delta})_B^{\beta-1}(\triangle_{\mathcal H}\psi_{\delta})_Bd\bar v
=-(\beta-1)\int_{\overline M}(\psi_{\delta})_B^{\beta-2}\vert\vert d\psi_{\delta}\vert\vert_B^2d\bar v
\leqno{(6.12)}$$
and hence 
$$\begin{array}{l}
\displaystyle{\int_{\overline M}(H^s_B)^{-1}(\psi_{\delta})_B^{\beta}
(\triangle_{\mathcal H}H^s)_Bd\bar v}\\
\hspace{0truecm}\displaystyle{=-2\beta\int_{\overline M}
(H^s_B)^{-1}(\psi_{\delta})_B^{\beta-1}
g_{\mathcal H}(dH^s,d\psi_{\delta})_Bd\bar v}\\
\hspace{0.5truecm}\displaystyle{+2\int_{\overline M}
(H^s_B)^{-2}(\psi_{\delta})_B^{\beta}
\vert\vert\,dH^s\,\,\vert\vert_B^2d\bar v.}
\end{array}\leqno{(6.13)}$$
By multiplying $\psi_{\delta}^{\beta-1}$ to both sides of the inequality in Lemma 6.2.6 and 
integrating the functions over $\overline M$ associated with both sides and using $(6.11),(6.12)$ and 
$(6.13)$, we can derive 
$$\begin{array}{l}
\displaystyle{\int_{\overline M}(H^s_B)^{1-\alpha}(\psi_{\delta})_B^{\beta-1}
{\rm Tr}((A_{\mathcal H})^3)_Bd\bar v
-\int_{\overline M}(H^s_B)^{-\alpha}(\psi_{\delta})_B^{\beta-1}
\vert\vert A_{\mathcal H}\vert\vert_B^4d\bar v}\\
\displaystyle{\leq
\int_{\overline M}(H^s_B)^{1-\alpha}(\psi_{\delta})_B^{\beta-1}
({\rm Tr}(({\mathcal A}^{\phi}_{\xi})^2\circ A_{\mathcal H}))_Bd\bar v}\\
\hspace{0.5truecm}\displaystyle{-\int_{\overline M}(H^s_B)^{-\alpha}(\psi_{\delta})_B^{\beta-1}
({\rm Tr}({\mathcal A}^{\phi}_{\xi})^2\vert_{\mathcal H})_B(\vert\vert A_{\mathcal H}\vert\vert^2)_Bd\bar v}\\
\hspace{0.5truecm}\displaystyle{-\frac{\beta-1}{2}\int_{\overline M}(\psi_{\delta})_B^{\beta-2}
\vert\vert d\psi_{\delta}\vert\vert_B^2d\bar v}\\
\hspace{0.5truecm}\displaystyle{-(\alpha\beta-\alpha+1)\int_{\overline M}
(H^s_B)^{-1}(\psi_{\delta})_B^{\beta-1}
g_{\mathcal H}(dH^s,d\psi_{\delta})_Bd\bar v}\\
\hspace{0.5truecm}\displaystyle{
-\alpha\int_{\overline M}(H^s_B)^{-\alpha-1}(\psi_{\delta})_B^{\beta-1}
g_{\mathcal H}((dH^s\otimes dH^s,h_{\mathcal H}-(H^s/n)
g_{\mathcal H})_Bd\bar v}\\
\hspace{0.5truecm}\displaystyle{+(\beta-1)\int_{\overline M}
(H^s_B)^{-\alpha}(\psi_{\delta})_B^{\beta-2}
g_{\mathcal H}((dH^s\otimes d\psi_{\delta},h_{\mathcal H}-(H^s/n)
g_{\mathcal H})_Bd\bar v}\\
\hspace{0.5truecm}\displaystyle{+(1-1/n)\int_{\overline M}
(H^s_B)^{-\alpha}(\psi_{\delta})_B^{\beta-1}
\vert\vert\,dH^s\,\,\vert\vert_B^2d\bar v}\\
\hspace{0.5truecm}\displaystyle{-3\int_{\overline M}(H^s_B)^{-\alpha}
(\psi_{\delta})_B^{\beta-1}{\rm Tr}_{g_{\mathcal H}}
({\mathcal A}^{\phi}_{\xi}\circ{\mathcal A}^{\phi}_{{\rm grad}H^s})_Bd\bar v}\\
\hspace{0.5truecm}\displaystyle{-\alpha\int_{\overline M}
(H^s_B)^{-2}(\psi_{\delta})_B^{\beta}
\vert\vert\,dH^s\,\,\vert\vert_B^2d\bar v}\\
\hspace{0.5truecm}\displaystyle{-\int_{\overline M}(\psi_{\delta})_B^{\beta-1}
(H^s_B)^{-\alpha}
({\rm Tr}_{g_{\mathcal H}}^{\bullet}{\mathcal R}(A_{\mathcal H}\bullet,\bullet))_Bd\bar v.}
\end{array}\leqno{(6.14)}$$
Denote by $\sharp_1$ the sum of the first term, the second one, the eight one and the last one 
in the right-hand side of $(6.14)$, and $\sharp_2$ the sum of the remained terms in the right-hand side of 
$(6.14)$.  Then, by simple calculations, we can derive 
$$\begin{array}{l}
\displaystyle{\sharp_1\leq 2\sqrt n\int_{\overline M}(\psi_{\delta})_B^{\beta-1}
(H^s_B)^{-\alpha}
\vert\vert({\mathcal A}^{\phi}_{\xi})^2\vert_{\mathcal H}\vert\vert_B\cdot\vert\vert A_{\mathcal H}\vert\vert^2_B}\\
\hspace{1truecm}\displaystyle{+3\sqrt n\int_{\overline M}(\psi_{\delta})_B^{\beta-1}
(H^s_B)^{-\alpha}
\vert\vert{\mathcal A}^{\phi}_{\xi}\circ{\mathcal A}^{\phi}_{{\rm grad}H^s}\vert\vert_Bd\bar v}\\
\hspace{1truecm}\displaystyle{
-\int_{\overline M}(\psi_{\delta})_B^{\beta-1}(H^s_B)^{-\alpha}
\left({\rm Tr}^{\bullet}_{g_{\mathcal H}}{\mathcal R}(A_{\mathcal H}\bullet,\bullet)\right)_Bd\bar v,}
\end{array}\leqno{(6.15)}$$
where we use $-{\rm Tr}(({\mathcal A}^{\phi}_{\xi})^2\vert_{\mathcal H})\leq\sqrt n
\vert\vert({\mathcal A}^{\phi}_{\xi})^2\vert_{\mathcal H}\vert\vert$ and 
$H^s\leq\sqrt{n\vert\vert A_{\mathcal H}\vert\vert^2}$.  
Also, by simple calculations, we can derive 
$$\begin{array}{l}
\displaystyle{\sharp_2\leq\alpha\int_{\overline M}(H^s_B)^{-\alpha/2-1}
(\psi_{\delta})_B^{\beta-1/2}\vert\vert dH^s\,\vert\vert_B^2d\bar v}\\
\hspace{1truecm}\displaystyle{+(\alpha\beta+\alpha-1)\int_{\overline M}
(H^s_B)^{-1}\cdot(\psi_{\delta})_B^{\beta-1}\cdot
\vert\vert dH^s\,\vert\vert_B\cdot
\vert\vert d\psi_{\delta}\vert\vert_Bd\bar v}\\
\hspace{1truecm}\displaystyle{+(\beta-1)\int_{\overline M}
(H^s_B)^{-\alpha/2}(\psi_{\delta})_B^{\beta-3/2}\vert\vert\,dH^s\,
\vert\vert_B\cdot\vert\vert\,d\psi_{\delta}\vert\vert_Bd\bar v}\\
\hspace{1truecm}\displaystyle{+\frac{n-1}{n}\int_{\overline M}
(H^s_B)^{-\alpha}(\psi_{\delta})_B^{\beta-1}
\vert\vert\,dH^s\,\,\vert\vert_B^2d\bar v,}
\end{array}\leqno{(6.16)}$$
where we use 
$\vert\vert\,dH^s\otimes dH^s\,\vert\vert
=\vert\vert\,dH^s\,\vert\vert^2$, 
$\vert\vert\,dH^s\otimes d\psi_{\delta}\vert\vert
=\vert\vert\,dH^s\,\vert\vert\cdot\vert\vert\,d\psi_{\delta}\vert\vert$ and 
$\vert\vert h_{\mathcal H}-\frac{H^s}{n}g_{\mathcal H}\vert\vert^2
=\psi_{\delta}(H^s)^{\alpha}$.  
By noticing $ab\leq\frac{\eta}{2}a^2+\frac{1}{2\eta}b^2$ for any $a,b,\eta>0$ and 
$\psi_{\delta}\leq(H^s)^{\delta}\,\,(0<\delta<1)$, we have 
$$\begin{array}{l}
\hspace{0.5truecm}\displaystyle{\int_{\overline M}
(H^s_B)^{-1}\cdot(\psi_{\delta})_B^{\beta-1}\cdot
\vert\vert dH^s\,\vert\vert_B\cdot\vert\vert d\psi_{\delta}\vert\vert_Bd\bar v}\\
\displaystyle{=\int_{\overline M}
((H^s_B)^{-\alpha/2}(\psi_{\delta})_B^{(\beta-1)/2}
\vert\vert dH^s\,\vert\vert_B)\cdot
((H^s_B)^{\alpha/2-1}(\psi_{\delta})_B^{(\beta-1)/2}
\vert\vert d\psi_{\delta}\vert\vert_B)d\bar v}\\
\displaystyle{\leq\frac{\eta}{2}\int_{\overline M}(H^s_B)^{-\alpha}
(\psi_{\delta})_B^{\beta-1}\vert\vert dH^s\,\vert\vert_B^2\,d\bar v
+\frac{1}{2\eta}\int_{\overline M}
(\psi_{\delta})_B^{\beta-2}\vert\vert d\psi_{\delta}\vert\vert_B^2d\bar v}
\end{array}\leqno{(6.17)}$$
and 
$$\begin{array}{l}
\hspace{0.5truecm}\displaystyle{\int_{\overline M}
(H^s_B)^{-\alpha/2}(\psi_{\delta})_B^{\beta-3/2}\vert\vert dH^s\,\vert\vert_B
\cdot\vert\vert d\psi_{\delta}\vert\vert_Bd\bar v}\\
\displaystyle{=\int_{\overline M}
\left((H^s_B)^{-\alpha/2}(\psi_{\delta})_B^{(\beta-1)/2}
\vert\vert dH^s\,\vert\vert_B\right)
\left((\psi_{\delta})_B^{(\beta-2)/2}\vert\vert d\psi_{\delta}\vert\vert_B\right)d\bar v}\\
\displaystyle{\leq\frac{\eta}{2}\int_{\overline M}(H^s_B)^{-\alpha}
(\psi_{\delta})_B^{\beta-1}\vert\vert dH^s\,\vert\vert_B^2\,d\bar v
+\frac{1}{2\eta}\int_{\overline M}
(\psi_{\delta})_B^{\beta-2}\vert\vert d\psi_{\delta}\vert\vert_B^2d\bar v.}
\end{array}\leqno{(6.18)}$$
From $(6.4)$ and $(6.14)-(6.18)$, we can derive 
$$\begin{array}{l}
\hspace{0.5truecm}\displaystyle{
n\varepsilon^2\int_{\overline M}(H^s_B)^2(\psi_{\delta})_B^{\beta}\,d\bar v}\\
\displaystyle{\leq\frac{(\alpha\beta+\alpha+\beta-2)\eta}{2}
\int_{\overline M}(H^s_B)^{-\alpha}(\psi_{\delta})_B^{\beta-1}
\vert\vert dH^s\,\vert\vert_B^2\,d\bar v}\\
\hspace{0.5truecm}\displaystyle{+\frac{\alpha\beta+\alpha+\beta-2}{2\eta}\int_{\overline M}
(\psi_{\delta})_B^{\beta-2}\vert\vert d\psi_{\delta}\vert\vert_B^2d\bar v}\\
\hspace{0.5truecm}\displaystyle{+\left(\alpha+\frac{n-1}{n}\right)
\int_{\overline M}(H^s_B)^{-\alpha}(\psi_{\delta})_B^{\beta-1}
\vert\vert dH^s\,\vert\vert_B^2\,d\bar v}\\
\hspace{0.5truecm}\displaystyle{+2\sqrt n\int_{\overline M}(\psi_{\delta})_B^{\beta-1}
(H^s_B)^{-\alpha}\cdot\vert\vert({\mathcal A}^{\phi}_{\xi})^2\vert_{\mathcal H}\vert\vert_B
\cdot\vert\vert A_{\mathcal H}\vert\vert^2_B\,d\bar v}\\
\hspace{0.5truecm}\displaystyle{+3\sqrt n\int_{\overline M}(\psi_{\delta})_B^{\beta-1}
(H^s_B)^{-\alpha}
\cdot\vert\vert{\mathcal A}^{\phi}_{\xi}\circ{\mathcal A}^{\phi}_{{\rm grad}H^s}\vert\vert_B\,d\bar v}\\
\hspace{0.5truecm}\displaystyle{-\int_{\overline M}(\psi_{\delta})_B^{\beta-1}
(H^s_B)^{-\alpha}
{\rm Tr}^{\bullet}_{g_{\mathcal H}}{\mathcal R}(A_{\mathcal H}\bullet,\bullet)_B\,d\bar v.}
\end{array}\leqno{(6.19)}$$
Since $0\leq\delta\leq\frac{1}{2}$ (hence $\frac{2}{3}\leq\alpha\leq 2$), we can derive the desired inequality.  
\qed

\vspace{0.5truecm}

Also, we can derive the following inequality.  

\vspace{0.5truecm}

\noindent
{\bf Lemma 6.2.8.} {\sl Assume that $0\leq\delta\leq\frac{1}{2}$.  
Then, for any $\beta\geq 100\varepsilon^{-2}$, we have 
$$\begin{array}{l}
\hspace{0.5truecm}\displaystyle{
\frac{\partial}{\partial t}\int_{\overline M}(\psi_{\delta})_B^{\beta}d\bar v
+2\int_{\overline M}(\psi_{\delta})_B^{\beta}(H^s_B)^2d\bar v
+\frac{\beta(\beta-1)}{2}\int_{\overline M}(\psi_{\delta})_B^{\beta-2}\vert\vert d\psi_{\delta}\vert\vert_B^2
d\bar v}\\
\hspace{0.5truecm}\displaystyle{+\frac{\beta\varepsilon^2}{8}\int_{\overline M}(\psi_{\delta})_B^{\beta-1}
(H^s_B)^{-\alpha}\cdot\vert\vert dH^s\,\vert\vert_B^2d\bar v}\\
\displaystyle{\leq\beta\delta\int_{\overline M}(\psi_{\delta})_B^{\beta}(H^s_B)^2d\bar v
+16\beta\varepsilon^{-2}K\int_{\overline M}(\psi_{\delta})_B^{\beta-1}(H^s_B)^{-\alpha}d\bar v}\\
\hspace{0.5truecm}\displaystyle{
-3\beta\delta\int_{\overline M}(\psi_{\delta})_B^{\beta}
{\rm Tr}(({\mathcal A}^{\phi}_{\xi})^2)_{\mathcal H})_Bd\bar v
-6\beta\int_{\overline M}(\psi_{\delta})_B^{\beta-1}(H^s_B)^{\delta}
{\rm Tr}(({\mathcal A}^{\phi}_{\xi})^2\circ A_{\mathcal H})_Bd\bar v}\\
\hspace{0.5truecm}\displaystyle{
+\frac{6\beta}{n}\int_{\overline M}(\psi_{\delta})_B^{\beta-1}(H^s_B)^{\delta}
{\rm Tr}(({\mathcal A}^{\phi}_{\xi})^2)_Bd\bar v}\\
\hspace{0.5truecm}\displaystyle{
-4\beta\int_{\overline M}(\psi_{\delta})_B^{\beta-1}(H^s_B)^{-\alpha}
{\rm Tr}^{\cdot}_{g_{\mathcal H}}{\rm Tr}^{\bullet}_{g_{\mathcal H}}
h(((\nabla_{\bullet}{\mathcal A})_{\bullet}\circ A_{\mathcal H})\cdot,\cdot)_Bd\bar v}\\
\hspace{0.5truecm}\displaystyle{
-4\beta\int_{\overline M}(\psi_{\delta})_B^{\beta-1}(H^s_B)^{-\alpha}
{\rm Tr}^{\cdot}_{g_{\mathcal H}}{\rm Tr}^{\bullet}_{g_{\mathcal H}}
h(({\mathcal A}_{\bullet}\circ A_{\mathcal H})\cdot,{\mathcal A}_{\bullet}\cdot)_Bd\bar v}\\
\hspace{0.5truecm}\displaystyle{
+4\beta\int_{\overline M}(\psi_{\delta})_B^{\beta-1}(H^s_B)^{-\alpha}
{\rm Tr}^{\cdot}_{g_{\mathcal H}}{\rm Tr}^{\bullet}_{g_{\mathcal H}}
h(({\mathcal A}_{\bullet}\circ{\mathcal A}_{\bullet})\cdot,A_{\mathcal H}\cdot)_Bd\bar v.}
\end{array}$$
}

\vspace{0.5truecm}

\noindent
{\it Proof.} By multiplying $\beta\psi_{\delta}^{\beta-1}$ to both sides of the inequality in Lemma 6.2.4 
and integrating over $\overline M$, we obtain 
$$\begin{array}{l}
\hspace{0.5truecm}\displaystyle{
\int_{\overline M}\left(\frac{\partial\psi_{\delta}^{\beta}}{\partial t}\right)_Bd\bar v
+\beta(\beta-1)\int_{\overline M}(\psi_{\delta})_B^{\beta-2}\vert\vert d\psi_{\delta}\vert\vert_B^2d\bar v}\\
\hspace{0.5truecm}\displaystyle{+\frac{\beta\varepsilon^2}{4}\int_{\overline M}(\psi_{\delta})_B^{\beta-1}
(H^s_B)^{-\alpha}\cdot\vert\vert dH^s\,\vert\vert_B^2d\bar v}\\
\displaystyle{\leq\beta\delta\int_{\overline M}(\psi_{\delta})_B^{\beta}(H^s_B)^2d\bar v
+2\beta(\alpha-1)\int_{\overline M}(\psi_{\delta})_B^{\beta-1}(H^s_B)^{-1}\cdot
\vert\vert\,dH^s\,\vert\vert_B\cdot\vert\vert d\psi_{\delta}\vert\vert_Bd\bar v}\\
\hspace{0.5truecm}\displaystyle{
+16\beta\varepsilon^{-2}K\int_{\overline M}(\psi_{\delta})_B^{\beta-1}(H^s_B)^{-\alpha}d\bar v
+\frac{6\beta}{n}\int_{\overline M}(\psi_{\delta})_B^{\beta-1}(H^s_B)^{\delta}
{\rm Tr}(({\mathcal A}^{\phi}_{\xi})^2)_Bd\bar v}\\
\hspace{0.5truecm}\displaystyle{
-3\beta\delta\int_{\overline M}(\psi_{\delta})_B^{\beta}{\rm Tr}(({\mathcal A}^{\phi}_{\xi})^2)_{\mathcal H})_B
d\bar v-6\beta\int_{\overline M}(\psi_{\delta})_B^{\beta-1}(H^s_B)^{\delta}
{\rm Tr}(({\mathcal A}^{\phi}_{\xi})^2\circ A_{\mathcal H})_Bd\bar v}\\
\hspace{0.5truecm}\displaystyle{
-4\beta\int_{\overline M}(\psi_{\delta})_B^{\beta-1}(H^s_B)^{-\alpha}
{\rm Tr}^{\cdot}_{g_{\mathcal H}}{\rm Tr}^{\bullet}_{g_{\mathcal H}}
h(((\nabla_{\bullet}{\mathcal A})_{\bullet}\circ A_{\mathcal H})\cdot,\cdot)_Bd\bar v}\\
\hspace{0.5truecm}\displaystyle{
-4\beta\int_{\overline M}(\psi_{\delta})_B^{\beta-1}(H^s_B)^{-\alpha}
{\rm Tr}^{\cdot}_{g_{\mathcal H}}{\rm Tr}^{\bullet}_{g_{\mathcal H}}
h(({\mathcal A}_{\bullet}\circ A_{\mathcal H})\cdot,{\mathcal A}_{\bullet}\cdot)_Bd\bar v}\\
\hspace{0.5truecm}\displaystyle{
+4\beta\int_{\overline M}(\psi_{\delta})_B^{\beta-1}(H^s_B)^{-\alpha}
{\rm Tr}^{\cdot}_{g_{\mathcal H}}{\rm Tr}^{\bullet}_{g_{\mathcal H}}
h(({\mathcal A}_{\bullet}\circ{\mathcal A}_{\bullet})\cdot,A_{\mathcal H}\cdot)_Bd\bar v,}
\end{array}$$
where we use $\int_{\overline M}\triangle_{\mathcal H}(\psi_{\delta}^{\beta})_Bd\bar v=0$ and 
$\vert\vert A_{\mathcal H}\vert\vert^2\leq(H^s)^2$.  
From this inequality, $\frac{\partial}{\partial t}(d\bar v)=-2(H^s_B)^2d\bar v,\,\,\,\,
\vert\vert\,dH^s\,\vert\vert\cdot
\vert\vert d\psi_{\delta}\vert\vert\leq\frac{\beta-1}{4(H^s)^{1-\alpha}}
\vert\vert d\psi_{\delta}\vert\vert^2+\frac{(H^s)^{1-\alpha}}{\beta-1}
\vert\vert\,dH^s\,\vert\vert^2$, $\alpha\leq2$, 
$\vert\vert A_{\mathcal H}\vert\vert^2\leq(H^s)^2$, 
$\psi_{\delta}\leq(H^s)^{\delta}$ and 
$\beta-1\geq 100\varepsilon^{-2}-1\geq16\varepsilon^{-2}$ (which holds because of $\varepsilon\leq 1$), 
we can derive the desired inequality.  \qed

\vspace{0.5truecm}

For a function $\overline{\rho}$ over $\overline M\times[0,T)$, denote by 
$\vert\vert\overline{\rho}(\cdot,t)\vert\vert_{L^{\beta},\bar g_t}$ the $L^{\beta}$-norm of with respect to 
$\bar g_t$ and $\vert\vert\overline{\rho}\vert\vert_{L^{\beta},\bar g}$ the function over $[0,T)$ defined by 
assigning $\vert\vert\overline{\rho}(\cdot,t)\vert\vert_{L^{\beta},\bar g_t}$ to each $t\in[0,T)$.  

By using Lemmas 6.2.7 and 6.2.8, we can derive the fact.  

\vspace{0.5truecm}

\noindent
{\bf Lemma 6.2.9.} {\sl There exists a positive constant $C$ depending only on $K,L$ and $f$ such that, for any 
$\delta$ and $\beta$ satisfying 
$$0\leq\delta\leq\min\left\{\frac{1}{2},\frac{n\varepsilon^2\eta}{3},\frac{n\varepsilon^4}{24(\eta+1)}
\right\}\,\,\,{\rm and}\,\,\,\beta\geq\max\left\{100\varepsilon^{-2},\frac{n\varepsilon^2\eta}{n\varepsilon^2
\eta-3\delta}\right\},
\leqno{(6.20)}$$
the following inequality holds:
$$\mathop{\sup}_{t\in[0,T)}\,\vert\vert(\psi_{\delta})_B(\cdot,t)\vert\vert_{L^{\beta},\bar g_t}\,<\,C.$$
}

\vspace{0.5truecm}

\noindent
{\it Proof.} Set 
$$C_1:=({\rm Vol}_{g_0}(M)+1)\mathop{\sup}_{\delta\in[0,1/2]}\,\mathop{\max}_M\,\psi_{\delta}(\cdot,0).$$
Then we have $\vert\vert\psi_{\delta}(\cdot,0)_B\vert\vert_{L^{\beta},g_0}\leq C_1$.  
By using the inequalities in Lemmas 6.2.7 and 6.2.8, 
$\vert\vert A_{\mathcal H}\vert\vert^2\leq(H^s)^2$ and the Young's inequality, 
we can show that 
$$\begin{array}{l}
\hspace{0.5truecm}\displaystyle{\frac{\partial}{\partial t}
\left(\vert\vert(\psi_{\delta})_B\vert\vert_{L^{\beta},\bar g}^{\beta}\right)}\\
\displaystyle{\leq\frac{\beta((3\delta-n\varepsilon^2\eta)\beta+n\varepsilon^2\eta)}{2n\varepsilon^2\eta}
\int_{\overline M}(\psi_{\delta})_B^{\beta-2}\vert\vert d\psi_{\delta}\vert\vert_B^2\,d\bar v}\\
\hspace{0.5truecm}\displaystyle{+\frac{\beta(12\eta\delta\beta+24\delta-n\varepsilon^4)}{8n\varepsilon^2}
\int_{\overline M}(\psi_{\delta})_B^{\beta-1}(H^s_B)^{-\alpha}
\vert\vert\,dH^s\,\,\vert\vert^2d\bar v}\\
\hspace{0.5truecm}\displaystyle{+C_2\vert\vert(\psi_{\delta})_B\vert\vert_{L^{\beta},\bar g}^{\beta}+C_3}
\end{array}\leqno{(6.21)}$$
holds for some positive constants $C_2$ and $C_3$ depending only on $K$ and $L$.  
Hence we can derive 
$$\begin{array}{l}
\hspace{0.5truecm}\displaystyle{
\mathop{\sup}_{t\in[0,T)}\vert\vert(\psi_{\delta})_B(\cdot,t)\vert\vert_{L^{\beta},\bar g_t}}\\
\displaystyle{\leq\left(\left(
\frac{C_3}{C_2}+\vert\vert(\psi_{\delta})_B(\cdot,0)\vert\vert_{L^{\beta},\bar g_0}^{\beta}
\right)e^{C_2T}-\frac{C_3}{C_2}\right)^{1/\beta}}\\
\displaystyle{\leq\left(\left(\frac{C_3}{C_2}+C_1\right)e^{C_2T}-\frac{C_3}{C_2}\right)^{1/\beta}.}
\end{array}$$
\qed

By using this lemma, we can derive the following inequality.  

\vspace{0.5truecm}

\noindent
{\bf Lemma 6.2.10.} {\sl Take any positive constant $k$.  Assume that 
$$0\leq\delta\leq\min\left\{\frac{1}{2}-\frac{k}{\beta},\frac{n\varepsilon^2\eta}{3}-\frac{k}{\beta},
\frac{n\varepsilon^4}{24(\eta+1)}-\frac{k}{\beta}\right\}\leqno{(6.22)}$$
and
$$\beta\geq\max\left\{100\varepsilon^{-2},\frac{n\varepsilon^2\eta}{n\varepsilon^2\eta-3\delta}\right\}.
\leqno{(6.23)}$$
Then the following inequality holds:
$$\mathop{\sup}_{t\in[0,T)}
\left(\int_{\overline M}(H^s_t)_B^k(\psi_{\delta}(\cdot,t))_B^{\beta}d\bar v\right)^{1/\beta}
\leq C,$$
where $C$ is as in Lemma 6.2.9.}

\vspace{0.5truecm}

\noindent
{\it Proof.} Set $\delta':=\delta+\frac{k}{\beta}$.  
Clearly we have $(H^s_t)_B^k(\psi_{\delta}(\cdot,t))_B^{\beta}=\psi_{\delta'}^{\beta}$.  
From the assumption for $\delta$ and $\beta$, $\delta'$ satisfies $(6.20)$.  
Hence, from Lemma 6.2.9, we have 
$$
\left(\int_{\overline M}(H^s_t)_B^k
(\psi_{\delta}(\cdot,t))_B^{\beta}d\bar v\right)^{1/\beta}
=\left(\int_{\overline M}(\psi_{\delta'}(\cdot,t))_B^{\beta}d\bar v\right)^{1/\beta}\leq C.
$$
\qed

\vspace{0.5truecm}

By using Lemmas 6.2.9, 6.2.10 and Theorem 5.2, we shall prove the statement of Proposition 6.2.  

\vspace{0.5truecm}

\noindent
{\it Proof of Proposition 6.2.} 
(Step I) First we shall show $T<\infty$.  According to Lemma 4.10, we have 
$$\frac{\partial H^s}{\partial t}\geq\triangle_{\mathcal H}H^s
+\frac{(H^s)^3}{n}.$$
Let $\rho$ be the solution of the ordinary differential equation 
$\displaystyle{\frac{dy}{\partial t}=\frac{1}{n}y^3}$ with the initial condition 
$\displaystyle{y(0)=\mathop{\min}_{M}\,H^s_0}$.  
This solution $\rho$ is given by 
$$\rho(t)=\frac{\min_{M}\,H^s_0}{\sqrt{1-(2/n)\min_{M}\,(H^s_0)^2\cdot t}}.$$
We regard $\rho$ as a function over $M\times[0,T)$.  Then we have 
$$\frac{\partial(H^s-\rho)}{\partial t}\geq\triangle_{\mathcal H}(H^s-\rho)
+\frac{(H^s)^3-\rho^3}{n}.$$
Furthermore, by the maximum principle, we can derive that 
$H^s\geq\rho$ holds over $M\times[0,T)$.  
Therefore we obtain 
$$H^s\geq
\frac{\min_{M}\,H^s_0}{\sqrt{1-(2/n)\min_{M}\,(H^s_0)^2\cdot t}}.$$
This implies that $T\leq\frac{1}{(2/n)\min_{M}\,(H^s_0)^2}(<\infty)$.  

\noindent
(Step II) Take positive constants $\delta$ and $\beta$ satisfying $(6.22)$ and $(6.23)$.  
Define a function $\psi_{\delta,k}$ by $\psi_{\delta,k}:=\max\{0,\psi_{\delta}(\cdot,t)-k\}$, where 
$k$ is any positive number with $\displaystyle{k\geq\mathop{\sup}_{M}\,\psi_{\delta}(\cdot,0)}$.  
Set $A_t(k):=\{\phi(u)\,\vert\,\psi_{\delta}(u,t)\geq k\}$ and 
$\displaystyle{\tilde A(k):=\mathop{\cup}_{t\in[0,T)}(A_t(k)\times\{t\})}$, which is finite because of 
$T<\infty$.  For a function $\bar{\rho}$ over $\overline M\times[0,T)$, denote by 
$\displaystyle{\int_{A(k)}\bar{\rho}\,d\bar v}$ the function over $[0,T)$ defined by 
assigning $\displaystyle{\int_{A_t(k)}\bar{\rho}(\cdot,t)\,d\bar v_t}$ to each $t\in[0,T)$.  
By multiplying the inequality in Lemma 6.2.4 by $\beta\psi_{\delta,k}^{\beta-1}$, 
we can show that the inequality in Lemma 6.2.8 holds for $\psi_{\delta,k}$ instead of $\psi_{\delta}$.  
From the inequality, the following inequality is derived directly:
$$
\frac{\partial}{\partial t}\int_{A(k)}(\psi_{\delta,k})_B^{\beta}d\bar v
+\frac{\beta(\beta-1)}{2}\int_{A(k)}(\psi_{\delta,k})_B^{\beta-2}\vert\vert d\psi_{\delta,k}\vert\vert_B^2
d\bar v
\leq\beta\delta\int_{A(k)}(\psi_{\delta,k})_B^{\beta}(H^s_B)^2d\bar v.
$$
Set $\hat{\psi}:=\psi_{\delta,k}^{\beta/2}$.  
On $A_t(k)$, we have 
$$\frac{\beta(\beta-1)}{2}(\psi_{\delta,k})^{\beta-2}_B(\cdot,t)\vert\vert d(\psi_{\delta,k})_B(\cdot,t)
\vert\vert^2\geq\vert\vert d\hat{\psi}_B(\cdot,t)\vert\vert^2$$
and hence 
$$
\frac{\partial}{\partial t}\int_{A(k)}\hat{\psi}_B^2d\bar v
+\int_{A(k)}\vert\vert d\hat{\psi}_B\vert\vert^2\,d\bar v
\leq\beta\delta\int_{A(k)}\hat{\psi}_B^2(H^s_B)^2d\bar v.
$$
By integrating both sides of this inequality from $0$ to any $t_0(\in[0,T))$, we have 
$$
\int_{A_{t_0}(k)}\hat{\psi}_B^2(\cdot,t_0)d\bar v_{t_0}
+\int_0^{t_0}\left(\int_{A(k)}\vert\vert d\hat{\psi}_B\vert\vert^2\,d\bar v\right)dt
\leq\beta\delta\int_0^{t_0}\left(\int_{A(k)}\hat{\psi}_B^2(H^s_B)^2d\bar v\right)dt,
$$
where we use $\displaystyle{k\geq\mathop{\sup}_{M}\,\psi_{\delta}(\cdot,0)}$.
By the arbitrariness of $t_0$, we have 
$$\begin{array}{l}
\displaystyle{\mathop{\sup}_{t\in[0,T)}\int_{A_t(k)}\hat{\psi}_B^2(\cdot,t)d\bar v_t
+\int_0^T\left(\int_{A(k)}\vert\vert d\hat{\psi}_B\vert\vert^2\,d\bar v\right)dt}\\
\displaystyle{\leq 2\beta\delta\int_0^T\left(\int_{A(k)}\hat{\psi}_B^2(H^s_B)^2
d\bar v\right)dt.}
\end{array}\leqno{(6.24)}$$
From $\displaystyle{k\geq\mathop{\sup}_{M}\,\psi_{\delta}(\cdot,0)}$, we have $A_0(k)=\emptyset$.  
Since $f$ satisfies the conditions $(\ast_1)$ and $(\ast_2)$, so is also $f_t$ ($0\leq t <T$) 
because ${\rm Vol}_{\bar g_t}(\overline M)$ decreases with respect to $t$ by Lemma 4.12.  
Hence we can apply the Sobolev's inequality in Theorem 5.2 to $f_t$ ($0\leq t<T$).  
By using the Sobolev's inequality in Theorem 5.2 and the H$\ddot{\rm o}$lder's inequality, we can derive 
$$\begin{array}{l}
\hspace{0.5truecm}\displaystyle{\left(\int_{\overline M}\hat{\psi}_B^{\frac{2n}{n-2}}(\cdot,t)\,d\bar v_t
\right)^{\frac{n-1}{n}}}\\
\displaystyle{\leq C(n)\left(\int_{\overline M}
\vert\vert d(\hat{\psi}^{\frac{2(n-1)}{n-2}})(\cdot,t)\vert\vert_B\,d\bar v_t
+\int_{\overline M}
\hat{\psi}_B^{\frac{2(n-1)}{n-2}}(\cdot,t)\cdot(H^s_t)_B\,d\bar v_t\right)}\\
\displaystyle{=C(n)\left(\frac{2(n-1)}{n-2}\int_{\overline M}
\hat{\psi}_B^{\frac{n}{n-2}}(\cdot,t)\cdot\vert\vert d\hat{\psi}(\cdot,t)\vert\vert_B\,d\bar v_t
+\int_{\overline M}
\hat{\psi}_B^{\frac{2(n-1)}{n-2}}(\cdot,t)\cdot(H^s_t)_B\,d\bar v_t\right)}\\
\displaystyle{\leq C(n)\left\{\frac{2(n-1)}{n-2}\left(\int_{\overline M}
\hat{\psi}_B^{\frac{2n}{n-2}}(\cdot,t)\,d\bar v_t\right)^{1/2}\cdot\left(\int_{\overline M}
\vert\vert d\hat{\psi}(\cdot,t)\vert\vert^2_B\,d\bar v_t\right)^{1/2}\right.}\\
\hspace{1.75truecm}\displaystyle{\left.+\left(\int_{\overline M}
\hat{\psi}_B^{\frac{2n}{n-2}}(\cdot,t)\,d\bar v_t\right)^{\frac{n-1}{n}}\cdot
\left(\int_{\overline M}(H^s_t)_B^n\,d\bar v_t\right)^{1/n}\right\}.}
\end{array}$$
Also, since $\psi_{\delta}(\cdot,t)\geq k$ on $A_t(k)$, it follows from Lemma 6.2.10 that 
$$\left(\int_{\overline M}(H^s_t)_B\,d\bar v_t\right)^{1/n}
\leq k^{-\beta/n}\left(\int_{\overline M}(H^s_t)_B\psi_{\delta}^{\beta}\,d\bar v_t\right)^{1/n}
\leq k^{-\beta/n}\cdot C^{\beta/n},$$
where $C$ is as in Lemma 6.2.9.  
Hence we obtain 
$$\begin{array}{l}
\hspace{0.5truecm}\displaystyle{\left(\int_{\overline M}\hat{\psi}_B^{\frac{2n}{n-2}}(\cdot,t)\,d\bar v_t
\right)^{\frac{n-1}{n}}}\\
\displaystyle{\leq C(n)\left\{\frac{2(n-1)}{n-2}\left(\int_{\overline M}
\hat{\psi}_B^{\frac{2n}{n-2}}(\cdot,t)\,d\bar v_t\right)^{1/2}\cdot\left(\int_{\overline M}
\vert\vert d\hat{\psi}(\cdot,t)\vert\vert^2_B\,d\bar v_t\right)^{1/2}\right.}\\
\hspace{1.75truecm}\displaystyle{\left.+\left(\int_{\overline M}
\hat{\psi}_B^{\frac{2n}{n-2}}(\cdot,t)\,d\bar v_t\right)^{\frac{n-1}{n}}\cdot k^{-\beta/n}\cdot C^{\beta/n}
\right\},}
\end{array}$$
that is, 
$$\begin{array}{l}
\hspace{0.5truecm}\displaystyle{\left(\int_{\overline M}
\vert\vert d\hat{\psi}(\cdot,t)\vert\vert^2_B\,d\bar v_t\right)^{1/2}}\\
\displaystyle{\geq\frac{n-2}{2C(n)(n-1)}
\left(\int_{\overline M}\hat{\psi}_B^{\frac{2n}{n-2}}(\cdot,t)\,d\bar v_t\right)^{\frac{n-2}{2n}}
\left(1-C(n)\cdot\left(\frac{C}{k}\right)^{\beta/n}\right).}
\end{array}$$
Set 
$$k_1:=\max\left\{\mathop{\sup}_M\psi_{\delta}(\cdot,0),\,\,C(n)^{n/\beta}\cdot C\right\}.$$
Assume that $k\geq k_1$.  Then we have 
$$\begin{array}{l}
\hspace{0.5truecm}\displaystyle{\int_{\overline M}
\vert\vert d\hat{\psi}(\cdot,t)\vert\vert^2_B\,d\bar v_t}\\
\displaystyle{\geq\left(\frac{n-2}{2C(n)(n-1)}\right)^2
\left(\int_{\overline M}\hat{\psi}_B^{\frac{2n}{n-2}}(\cdot,t)\,d\bar v_t\right)^{\frac{n-2}{n}}
\left(1-C(n)\cdot\left(\frac{C}{k}\right)^{\beta/n}\right)^2.}
\end{array}\leqno{(6.25)}$$
From $(6.24)$ and $(6.25)$, we obtain 
$$\begin{array}{l}
\hspace{0.5truecm}\displaystyle{\mathop{\sup}_{t\in[0,T]}\int_{\overline M}\hat{\psi}_B^2(\cdot,t)d\bar v_t
+\hat C(n,k)\int_0^{T}\left(\int_{\overline M}\hat{\psi}_B^{\frac{2n}{n-2}}(\cdot,t)\,d\bar v_t
\right)^{\frac{n-2}{n}}dt}\\
\displaystyle{\leq 2\beta\delta\int_0^{T}\left(\int_{\overline M}\hat{\psi}_B^2(H^s_t)_B^2d\bar v_t\right)dt,}
\end{array}\leqno{(6.26)}$$
where $\hat C(n,k):=\left(\frac{(n-2)(1-C(n)\cdot(C/k)^{\beta/n})}{2C(n)(n-1)}\right)^2$.  
Set 
$$q:=\left\{\begin{array}{cc}
\displaystyle{\frac{n}{n-2}} & (n\geq 3)\\
{\rm any}\,\,{\rm positive}\,\,{\rm number} & (n=2)
\end{array}\right.
$$
and $q_0:=2-1/q$ and 
$$\vert\vert A_t(k)\vert\vert_{T}:=\int_0^{T}\left(\int_{A_t(k)}d\bar v_t\right)dt.$$
By using the interpolation inequality, we can derive 
$$\left(\int_{\overline M}\hat{\psi}_B^{2q_0}\,d\bar v_t\right)^{1/q_0}\leq
\left(\int_{\overline M}\hat{\psi}_B^2\,d\bar v_t\right)^{1-1/q_0}\cdot
\left(\int_{\overline M}\hat{\psi}_B^{2q}\,d\bar v_t\right)^{1/qq_0}.$$
By using this inequality and the Young inequality, we  can derive 
$$\begin{array}{l}
\hspace{0.5truecm}\displaystyle{
\left(\int_0^{T}\left(\int_{\overline M}\hat{\psi}_B^{2q_0}(\cdot,t)\,d\bar v_t\right)dt\right)^{1/q_0}}\\
\displaystyle{\leq\left(\int_0^{T}\left(
\left(\int_{\overline M}\hat{\psi}_B^2(\cdot,t)\,d\bar v_t\right)^{q_0-1}
\cdot\left(\int_{\overline M}\hat{\psi}_B^{2q}(\cdot,t)\,d\bar v_t\right)^{1/q}\right)dt\right)^{1/q_0}}\\
\displaystyle{\leq\left(\mathop{\sup}_{t\in[0,T]}\int_{\overline M}\hat{\psi}_B^2(\cdot,t)\,d\bar v_t
\right)^{\frac{q_0-1}{q_0}}
\cdot\left(\int_0^{T}\left(\int_{\overline M}\hat{\psi}_B^{2q}(\cdot,t)\,d\bar v_t\right)^{1/q}dt
\right)^{1/q_0}}\\
\displaystyle{\leq\mathop{\sup}_{t\in[0,T]}\int_{\overline M}\hat{\psi}_B^2(\cdot,t)\,d\bar v_t
+\int_0^{T}\left(\int_{\overline M}\hat{\psi}_B^{2q}(\cdot,t)\,d\bar v_t\right)^{1/q}dt.}
\end{array}\leqno{(6.27)}$$
We may assume that $\hat C(n,k)<1$ holds by replacing $C(n)$ to a bigger positive number and furthemore 
$k$ to a positive number bigger such that 
$1-C(n)\cdot\left(\frac{C}{k}\right)^{\beta/n}>0$ holds for the replaced number $C(n)$.  
Then, from $(6.26)$ and $(6.27)$, we obtain 
$$\begin{array}{l}
\hspace{0.5truecm}\displaystyle{\hat C(n,k)
\left(\int_0^{T}\left(\int_{\overline M}\hat{\psi}_B^{2q_0}(\cdot,t)\,d\bar v_t\right)dt
\right)^{1/q_0}}\\
\displaystyle{\leq 2\beta\delta\int_0^{T}\left(\int_{\overline M}\hat{\psi}_B^2
(H^s_t)_B^2d\bar v_t\right)dt.}
\end{array}\leqno{(6.28)}$$
On the other hand, by using the H$\ddot{\rm o}$lder's inequality, we obtain 
$$
\int_0^{T}\left(\int_{\overline M}\hat{\psi}_B^2(H^s_t)_B^2d\bar v_t\right)dt
\leq\vert\vert A_t(k)\vert\vert_{T}^{\frac{r-1}{r}}\cdot
\left(\int_0^{T}\left(\int_{\overline M}\hat{\psi}_B^{2r}(H^s_t)_B^{2r}d\bar v_t\right)dt
\right)^{1/r},
$$
where $r$ is any positive constant with $r>1$.  
From $(6.28)$ and this inequality, we obtain 
$$\begin{array}{l}
\hspace{0.5truecm}\displaystyle{\left(\int_0^{T}\left(\int_{\overline M}\hat{\psi}_B^{2q_0}(\cdot,t)\,
d\bar v_t\right)dt\right)^{1/q_0}}\\
\displaystyle{\leq 2\hat C(n,k)^{-1}\beta\delta\vert\vert A_t(k)\vert\vert_{T}^{\frac{r-1}{r}}
\cdot\left(\int_0^{T}\left(\int_{\overline M}\hat{\psi}_B^{2r}(H^s_t)_B^{2r}d\bar v_t\right)dt
\right)^{1/r}.}
\end{array}\leqno{(6.29)}$$
On the other hand, according to Lemma 6.2.10, we have 
$$\int_{\overline M}\widehat{\psi}^{2r}_B(H^s_t)_B^{2r}\,d\bar v_t\leq C^{2r}\leqno{(6.30)}$$
for some positive constant $C$ (depending only on $K,L$ and $f$) by replacing $r$ to a bigger positive number 
if necessary.  Also, by using the H$\ddot{\rm o}$lder inequality, we obtain 
$$\begin{array}{l}
\hspace{0.5truecm}\displaystyle{\int_0^{T}\left(\int_{\overline M}\hat{\psi}_B^2(\cdot,t)\,d\bar v_t\right)dt}\\
\displaystyle{\leq\vert\vert A_t(k)\vert\vert_{T}^{\frac{q_0-1}{q_0}}\cdot
\left(\int_0^{T}\left(\int_{\overline M}\hat{\psi}_B^{2q_0}(\cdot,t)\,d\bar v_t\right)dt\right)^{1/q_0}.}
\end{array}$$
From $(6.29),\,(6.30)$ and this inequality, we obtain 
$$
\int_0^{T}\left(\int_{\overline M}\hat{\psi}_B^2(\cdot,t)\,d\bar v_t\right)dt
\leq\vert\vert A_t(k)\vert\vert_{T}^{2-1/q_0-1/r}\cdot C^2\cdot\hat C(n,k)^{-1}\cdot 2\beta\delta.
\leqno{(6.31)}$$
We may assume that $2-1/{q_0}-1/r>1$ holds by replacing $r$ to a bigger positive number if necessary.  
Take any positive constants $h$ and $k$ with $h>k\geq k_1$.  Then we have 
$$\begin{array}{l}
\displaystyle{\int_0^T\left(\int_{\overline M}\psi^{\beta}_{\delta,k}d\overline v_t\right)dt
\geq\int_0^T\left(\int_{\overline M}(\psi_{\delta,k}-\psi_{\delta,h})^{\beta}d\overline v_t\right)dt}\\
\displaystyle{\geq\int_0^T\left(\int_{A_t(h)}\vert h-k\vert^{\beta}\,d\overline v_t\right)dt
=\vert h-k\vert^{\beta}\cdot\vert\vert A_t(h)\vert\vert_T.}
\end{array}$$
From this inequality and $(6.31)$, we obtain 
$$\vert h-k\vert^{\beta}\cdot\vert\vert A_t(h)\vert\vert_T
\leq\vert\vert A_t(k)\vert\vert_{T}^{2-1/q_0-1/r}\cdot C^2\cdot\hat C(n,k)^{-1}\cdot 2\beta\delta.
\leqno{(6.32)}$$
Since $\bullet\mapsto\vert\vert A_t(\bullet)\vert\vert_T$ is a non-increasing and non-negative function and 
$(6.32)$ holds for any $h>k\geq k_1$, it follows from the Stambaccha's iteration lemma that 
$\vert\vert A_t(k_1+d)\vert\vert_T=0$, where $d$ is a positive constant depending only on 
$\beta,\delta,q_0,r,C,\widehat C(n,k)$ and $\vert\vert A_t(k_1)\vert\vert_T$.  
This implies that 
$\displaystyle{\mathop{\sup}_{t\in[0,T)}\,\mathop{\max}_{\overline M}\,
\psi_{\delta}(\cdot,t)\leq}$\newline
$k_1+d<\infty$.  This completes the proof.  \qed

\section{Estimate of the gradient of the mean curvature from above} 
In this section, we shall derive the following estimate of ${\rm grad}H^s$ from above 
by using Proposition 6.2.  

\vspace{0.5truecm}

\noindent
{\bf Proposition 7.1.} {\sl For any positive constant $b$, there exists a constant $C(b,f_0)$ 
(depending only on $b$ and $f_0$) satisfying 
$$\vert\vert{\rm grad}\,H^s\,\vert\vert^2\leq b\cdot(H^s)^4
+C(b,f_0)\,\,\,\,{\rm on}\,\,\,\,M\times[0,T).$$
}

\vspace{0.5truecm}

We prepare some lemmas to prove this proposition.  

\vspace{0.5truecm}

\noindent
{\bf Lemma 7.1.1.} {\sl The family $\{\vert\vert{\rm grad}_tH^s_t\,\vert\vert^2\}_{t\in[0,T)}$ 
satisfies the following equation:
$$\begin{array}{l}
\hspace{0.5truecm}\displaystyle{\frac{\partial\vert\vert{\rm grad}\,H^s\,\vert\vert^2}
{\partial t}-\triangle_{\mathcal H}(\vert\vert{\rm grad}\,H^s\,\vert\vert^2)}\\
\displaystyle{=-2\vert\vert\nabla^{\mathcal H}{\rm grad}\,H^s\,\vert\vert^2
+2\vert\vert A_{\mathcal H}\vert\vert^2\cdot\vert\vert{\rm grad}\,H^s\,\vert\vert^2}\\
\hspace{0.5truecm}\displaystyle{+2H^s\cdot g_{\mathcal H}({\rm grad}(\vert\vert A_{\mathcal H}\vert\vert^2),
{\rm grad}\,H^s)}\\
\hspace{0.5truecm}\displaystyle{+2g_{\mathcal H}((A_{\mathcal H})^2({\rm grad}\,H^s),{\rm grad}\,
H^s)}\\
\hspace{0.5truecm}\displaystyle{-6H^s\cdot g_{\mathcal H}
({\rm grad}({\rm Tr}(({\mathcal A}^{\phi}_{\xi})^2)_{\mathcal H},{\rm grad}\,H^s)}\\
\hspace{0.5truecm}\displaystyle{-6{\rm Tr}(({\mathcal A}^{\phi}_{\xi})^2)_{\mathcal H}\cdot
\vert\vert{\rm grad}\,H^s\,\vert\vert^2.}
\end{array}\leqno{(7.1)}$$
Hence we have the following inequality:
$$\begin{array}{l}
\hspace{0.5truecm}\displaystyle{\frac{\partial\vert\vert{\rm grad}\,H^s\,\vert\vert^2}
{\partial t}-\triangle_{\mathcal H}(\vert\vert{\rm grad}\,H^s\,\vert\vert^2)}\\
\displaystyle{\leq-2\vert\vert\nabla^{\mathcal H}{\rm grad}\,H^s\,\vert\vert^2
+4\vert\vert A_{\mathcal H}\vert\vert^2\cdot\vert\vert{\rm grad}\,H^s\,\vert\vert^2}\\
\hspace{0.5truecm}\displaystyle{+2H^s\cdot g_{\mathcal H}
({\rm grad}(\vert\vert A_{\mathcal H}\vert\vert^2),{\rm grad}\,H^s)}\\
\hspace{0.5truecm}\displaystyle{+6H^s\cdot
\vert\vert{\rm grad}({\rm Tr}(({\mathcal A}^{\phi}_{\xi})^2)_{\mathcal H})\vert\vert
\cdot\vert\vert{\rm grad}\,H^s\,\vert\vert}\\
\hspace{0.5truecm}\displaystyle{-6{\rm Tr}(({\mathcal A}^{\phi}_{\xi})^2)_{\mathcal H}\cdot
\vert\vert{\rm grad}\,H^s\,\vert\vert^2.}
\end{array}\leqno{(7.2)}$$
}

\vspace{0.5truecm}

\noindent
{\it Proof.} By using Lemmas 4.2 and 4.9, we have 
$$\begin{array}{l}
\displaystyle{\frac{\partial\vert\vert{\rm grad}\,H^s\,\vert\vert^2}
{\partial t}=\frac{\partial g_{\mathcal H}}{\partial t}
({\rm grad}\,H^s,{\rm grad}\,H^s)
+2g_{\mathcal H}\left({\rm grad}\left(\frac{\partial H^s}{\partial t}\right),
{\rm grad}\,H^s\right)}\\
\displaystyle{=-2H^s\cdot h_{\mathcal H}
({\rm grad}\,H^s,{\rm grad}\,H^s)
+2g_{\mathcal H}({\rm grad}(\triangle_{\mathcal H}H^s),\,{\rm grad}\,H^s)}\\
\hspace{0.5truecm}\displaystyle{
+2g_{\mathcal H}({\rm grad}(H^s\cdot\vert\vert A_{\mathcal H}\vert\vert^2)),\,
{\rm grad}\,H^s)}\\
\hspace{0.5truecm}\displaystyle{-6g_{\mathcal H}({\rm grad}(H^s\cdot
{\rm Tr}(({\mathcal A}^{\phi}_{\xi})^2)_{\mathcal H},\,{\rm grad}\,H^s).}
\end{array}$$
Also we have 
$$\begin{array}{l}
\displaystyle{\triangle_{\mathcal H}(\vert\vert{\rm grad}\,H^s\,\vert\vert^2)
=2g_{\mathcal H}(\triangle_{\mathcal H}^{\mathcal H}({\rm grad}\,H^s),\,
{\rm grad}\,H^s)}\\
\hspace{3.8truecm}\displaystyle{+2g_{\mathcal H}(\nabla^{\mathcal H}{\rm grad}\,H^s,\,
\nabla^{\mathcal H}{\rm grad}\,H^s)}
\end{array}$$
and 
$$\triangle_{\mathcal H}^{\mathcal H}({\rm grad}\,H^s)
={\rm grad}(\triangle_{\mathcal H}H^s)+H^s\cdot A_{\mathcal H}
({\rm grad}\,H^s)-(A_{\mathcal H})^2({\rm grad}\,H^s).$$
By using these relations and noticing $g_{\mathcal H}(A_{\mathcal H}(\bullet),\cdot)=-h_{\mathcal H}(\bullet,\cdot)$, 
we can derive the desired evolution equation $(7.1)$.  
The inequality $(7.2)$ is derived from $(7.1)$ and 
$$g_{\mathcal H}((A_{\mathcal H})^2({\rm grad}\,H^s),\,{\rm grad}\,H^s)
\leq\vert\vert A_{\mathcal H}\vert\vert^2\cdot\vert\vert{\rm grad}\,H^s\,\vert\vert^2.$$
\qed

\vspace{0.5truecm}

\noindent
{\bf Lemma 7.1.2.} {\sl The family 
$\displaystyle{\left\{\frac{\vert\vert{\rm grad}_tH^s_t\,\vert\vert^2}
{H^s_t}\right\}_{t\in[0,T)}}$ satisfies the following inequality:
$$\begin{array}{l}
\hspace{0.5truecm}\displaystyle{\frac{\partial}{\partial t}\left(
\frac{\vert\vert{\rm grad}\,H^s\,\vert\vert^2}
{H^s}\right)
-\triangle_{\mathcal H}\left(\frac{\vert\vert{\rm grad}\,H^s\,\vert\vert^2}
{H^s}\right)}\\
\displaystyle{\leq\frac{3\vert\vert{\rm grad}\,H^s\,\vert\vert^2}
{H^s}\cdot\vert\vert A_{\mathcal H}\vert\vert^2+2g_{\mathcal H}
({\rm grad}(\vert\vert A_{\mathcal H}\vert\vert^2,\,{\rm grad}\,H^s)}\\
\hspace{0.5truecm}\displaystyle{+6\vert\vert{\rm grad}({\rm Tr}(({\mathcal A}^{\phi}_{\xi})^2)_{\mathcal H}\vert\vert\cdot
\vert\vert{\rm grad}\,H^s\,\vert\vert
-\frac{3}{H^s}{\rm Tr}(({\mathcal A}^{\phi}_{\xi})^2)_{\mathcal H}\cdot
\vert\vert{\rm grad}\,H^s\,\vert\vert^2}
\end{array}\leqno{(7.3)}$$
}

\vspace{0.5truecm}

\noindent
{\it Proof.} By a simple calculation, we have 
$$\begin{array}{l}
\hspace{0.5truecm}\displaystyle{\frac{\partial}{\partial t}\left(
\frac{\vert\vert{\rm grad}\,H^s\,\vert\vert^2}
{H^s}\right)
-\triangle_{\mathcal H}\left(\frac{\vert\vert{\rm grad}\,H^s\,\vert\vert^2}
{H^s}\right)}\\
\displaystyle{=\frac{1}{H^s}
\left(\frac{\partial\vert\vert{\rm grad}\,H^s\,\vert\vert^2}{\partial t}
-\triangle_{\mathcal H}(\vert\vert{\rm grad}\,H^s\,\vert\vert^2)\right)}\\
\hspace{0.5truecm}\displaystyle{-\frac{\vert\vert{\rm grad}\,\vert\vert H\vert^\vert\,\vert\vert^2}
{(H^s)^2}\left(\frac{\partial H^s}{\partial t}
-\triangle_{\mathcal H}H^s\right)}\\
\hspace{0.5truecm}\displaystyle{+\frac{2}{(H^s)^2}g_{\mathcal H}
({\rm grad}\,H^s,\,{\rm grad}(\vert\vert{\rm grad}\,H^s\,\vert\vert^2)).}
\end{array}$$
From this relation, Lemmas 4.9 and (7.2), we can derive the desired inequality.  \qed

\vspace{0.5truecm}

From Lemma 4.10, we can derive the following evolution equation directly.  

\vspace{0.5truecm}

\noindent
{\bf Lemma 7.1.3.} {\sl The family 
$\displaystyle{\left\{(H^s_t)^3\right\}_{t\in[0,T)}}$ satisfies the following evolution 
equation:
$$\begin{array}{l}
\hspace{0.5truecm}\displaystyle{\frac{\partial(H^s)^3}{\partial t}
-\triangle_{\mathcal H}((H^s)^3)}\\
\displaystyle{=3(H^s)^3\cdot\vert\vert A_{\mathcal H}\vert\vert^2-6H^s\cdot
\vert\vert{\rm grad}\,H^s\,\vert\vert^2-9(H^s)^3\cdot
{\rm Tr}(({\mathcal A}^{\phi}_{\xi})^2)_{\mathcal H}.}
\end{array}$$
}

By using Lemmas 4.9, 4.11 and Proposition 6.2, we can derive the following evolution inequality.  

\vspace{0.5truecm}

\noindent
{\bf Lemma 7.1.4.} {\sl The family 
$\displaystyle{\left\{\left(\vert\vert(A_{\mathcal H})_t\vert\vert^2-\frac{(H^s_t)^2}{n}\right)
\cdot H^s_t\right\}_{t\in[0,T)}}$ satisfies the following evolution inequality:
$$\begin{array}{l}
\hspace{0.5truecm}\displaystyle{\frac{\partial}{\partial t}
\left(\left(\vert\vert A_{\mathcal H}\vert\vert^2-\frac{(H^s)^2}{n}\right)
\cdot H^s\right)
-\triangle_{\mathcal H}\left(\left(\vert\vert A_{\mathcal H}\vert\vert^2-\frac{(H^s)^2}{n}\right)
\cdot H^s\right)}\\
\displaystyle{\leq-\frac{2(n-1)}{3n}H^s\cdot\vert\vert\nabla^{\mathcal H}A_{\mathcal H}
\vert\vert^2+\breve C(n,C_0,\delta)\cdot\vert\vert\nabla^{\mathcal H}A_{\mathcal H}\vert\vert^2}\\
\hspace{0.5truecm}\displaystyle{+3H^s\cdot\vert\vert A_{\mathcal H}\vert\vert^2\cdot
\left(\vert\vert A_{\mathcal H}\vert\vert^2-\frac{(H^s)^2}{n}\right)
-2H^s\cdot{\rm Tr}(({\mathcal A}^{\phi}_{\xi})^2)_{\mathcal H}
\cdot\left(\vert\vert A_{\mathcal H}\vert\vert^2-\frac{(H^s)^2}{n}\right)}\\
\hspace{0.5truecm}\displaystyle{-4(H^s)^2\cdot{\rm Tr}\left(({\mathcal A}^{\phi}_{\xi})^2)\circ
\left(A_{\mathcal H}-\frac{H^s}{n}\cdot{\rm id}\right)\right)}\\
\hspace{0.5truecm}\displaystyle{
-2H^s\cdot{\rm Tr}_{g_{\mathcal H}}^{\bullet}
{\mathcal R}\left(\left(A_{\mathcal H}-\frac{H^s}{n}\cdot{\rm id}\right)(\bullet),\bullet\right)}\\
\hspace{0.5truecm}\displaystyle{-3\left(\vert\vert A_{\mathcal H}\vert\vert^2
-\frac{(H^s)^2}{n}\right)
\cdot H^s\cdot{\rm Tr}(({\mathcal A}^{\phi}_{\xi})^2)_{\mathcal H}.}
\end{array}$$
}

\vspace{0.5truecm}

\noindent
{\it Proof.} By using Lemmas 4.9 and 4.11, we can derive 
$$\begin{array}{l}
\hspace{0.5truecm}\displaystyle{\frac{\partial}{\partial t}
\left(\left(\vert\vert A_{\mathcal H}\vert\vert^2-\frac{(H^s)^2}{n}\right)\cdot H^s\right)
-\triangle_{\mathcal H}\left(\left(\vert\vert A_{\mathcal H}\vert\vert^2-\frac{(H^s)^2}{n}\right)\cdot H^s\right)}\\
\displaystyle{=\frac{2H^s}{n}\cdot\vert\vert{\rm grad}\,H^s\,\vert\vert^2
+2H^s\cdot\vert\vert A_{\mathcal H}\vert\vert^2\cdot
\left({\rm Tr}(((A_{\mathcal H})^2)-\frac{(H^s)^2}{n}\right)}\\
\hspace{0.5truecm}\displaystyle{-2H^s\cdot\vert\vert\nabla^{\mathcal H}A_{\mathcal H}\vert\vert^2
+H^s\cdot\vert\vert A_{\mathcal H}\vert\vert^2\cdot
\left(\vert\vert A_{\mathcal H}\vert\vert^2-\frac{(H^s)^2}{n}\right)}\\
\hspace{0.5truecm}\displaystyle{-g_{\mathcal H}\left({\rm grad}\left(
\vert\vert A_{\mathcal H}\vert\vert^2-\frac{H^s}{n}\cdot{\rm id}\right),\,
{\rm grad}\,H^s\right)}\\
\hspace{0.5truecm}\displaystyle{-2H^s\cdot{\rm Tr}(({\mathcal A}^{\phi}_{\xi})^2)_{\mathcal H}
\cdot\left({\rm Tr}(((A_{\mathcal H})^2)-\frac{(H^s)^2}{n}\right)}\\
\hspace{0.5truecm}\displaystyle{-4(H^s)^2\cdot{\rm Tr}
\left(({\mathcal A}^{\phi}_{\xi})^2)\circ
\left(A_{\mathcal H}-\frac{H^s}{n}\cdot{\rm id}\right)\right)}\\
\hspace{0.5truecm}\displaystyle{-2H^s\cdot{\rm Tr}_{g_{\mathcal H}}^{\bullet}
{\mathcal R}\left(\left(A_{\mathcal H}-\frac{H^s}{n}\cdot{\rm id}\right)(\bullet),\bullet\right)}\\
\hspace{0.5truecm}\displaystyle{
-3\left(\vert\vert A_{\mathcal H}\vert\vert^2-\frac{(H^s)^2}{n}\right)
\cdot H^s\cdot{\rm Tr}(({\mathcal A}^{\phi}_{\xi})^2)_{\mathcal H}.}
\end{array}\leqno{(7.4)}$$
On the othe hand, by using 
$\vert\vert A_{\mathcal H}\vert\vert^2-\frac{(H^s)^2}{n}=\vert\vert 
A_{\mathcal H}-\frac{H^s}{n}\cdot{\rm id}\,\vert\vert$, 
we can derive 
$$\begin{array}{l}
\hspace{0.5truecm}\displaystyle{
\left\vert g_{\mathcal H}\left({\rm grad}\left(\vert\vert A_{\mathcal H}\vert\vert^2-\frac{(H^s)^2}{n}
\right),\,{\rm grad}\,H^s\right)\right\vert}\\
\displaystyle{=\left\vert d\left(\vert\vert A_{\mathcal H}\vert\vert^2-\frac{(H^s)^2}{n}\right)
({\rm grad}\,H^s)\right\vert}\\
\hspace{0truecm}\displaystyle{=2\left\vert g_{\mathcal H}\left(\nabla^{\mathcal H}_{{\rm grad}H^s}
\left(A_{\mathcal H}-\frac{H^s}{n}\cdot{\rm id}\right),\,
A_{\mathcal H}-\frac{H^s}{n}\cdot{\rm id}\right)\right\vert}
\end{array}$$

$$\begin{array}{l}
\hspace{0truecm}\displaystyle{\leq 2\vert\vert{\rm grad}H^s\,\vert\vert
\cdot\vert\vert\nabla^{\mathcal H}A_{\mathcal H}\vert\vert\cdot
\left\vert\left\vert A_{\mathcal H}-\frac{H^s}{n}\cdot{\rm id}\right\vert\right\vert}\\
\hspace{0truecm}\displaystyle{\leq 2n\vert\vert\nabla^{\mathcal H}A_{\mathcal H}\vert\vert^2\cdot
\left\vert\left\vert A_{\mathcal H}-\frac{H^s}{n}\cdot{\rm id}\right\vert\right\vert,}
\end{array}$$
where we use $\frac{1}{n}\vert\vert{\rm grad}H^s\,\vert\vert^2
\leq\vert\vert\nabla^{\mathcal H}A_{\mathcal H}\vert\vert^2$.  
Also, according to Proposition 6.2, we have 
$$\left\vert\left\vert A_{\mathcal H}-\frac{H^s}{n}\cdot{\rm id}\right\vert\right\vert
\leq\sqrt{C_0}\cdot(H^s)^{1-\delta/2}.$$
Hence we have 
$$\begin{array}{l}
\hspace{0.5truecm}\displaystyle{
\left\vert g_{\mathcal H}\left({\rm grad}\left(\vert\vert A_{\mathcal H}\vert\vert^2-\frac{(H^s)^2}{n}
\right),\,{\rm grad}\,H^s\right)\right\vert}\\
\displaystyle{\leq 2n\sqrt{C_0}\vert\vert\nabla^{\mathcal H}A_{\mathcal H}\vert\vert^2\cdot
(H^s)^{1-\delta/2}.}
\end{array}\leqno{(7.5)}$$
Furthermore, according to the Young's inequality:
$$ab\leq\varepsilon\cdot a^p+\varepsilon^{-1/(p-1)}\cdot b^q\,\,\,(\forall\,a>0,\,\,b>0)\leqno{(7.6)}.$$
(where $p$ and $q$ are any positive constants with $\frac{1}{p}+\frac{1}{q}=1$ and $\varepsilon$ is any positive 
constant), we have 
$$2n\sqrt{C_0}(H^s)^{1-\delta/2}\leq\frac{2(n-1)}{3n}\cdot H^s
+\breve C(n,C_0,\delta),\leqno{(7.7)}$$
where $\breve C(n,C_0,\delta)$ is a positive constant only on $n,C_0$ and $\delta$.  
Also, we have 
$$\vert\vert\nabla^{\mathcal H}A_{\mathcal H}\vert\vert^2
\geq\frac{3}{n+2}\vert\vert\nabla^{\mathcal H}H\vert\vert^2.$$
From $(7.4)$ and these inequalities, we can derive the desired evolution inequality.  \qed

\vspace{0.5truecm}

By using Lemmas 4.10, 7.1.2, 7.1.3 and 7.1.4, we shall prove Theorem Proposition 7.1.  

\vspace{0.5truecm}

\noindent
{\it Proof of Proposition 7.1.} Define a function $\rho$ over $M\times[0,T)$ by 
$$\rho:=\frac{\vert\vert{\rm grad}\,H^s\,\vert\vert^2}{H^s}
+C_1H^s\left(\vert\vert A_{\mathcal H}\vert\vert^2-\frac{(H^s)^2}{n}\right)
+C_1\cdot\breve C(n,C_0,\delta)\vert\vert A_{\mathcal H}\vert\vert^2-b(H^s)^3,$$
where $b$ is any positive constant and $C_1$ is a positive constant which is sufficiently big compared to 
$n$ and $b$.  
By using Lemmas 4.10, 7.1.2, 7.1.3 and 7.1.4, we can derive 
$$\begin{array}{l}
\hspace{0.5truecm}\displaystyle{\frac{\partial\rho}{\partial t}-\triangle_{\mathcal H}\rho}\\
\displaystyle{\leq\frac{3\vert\vert{\rm grad}\,H^s\,\vert\vert^2}
{H^s}\cdot\vert\vert A_{\mathcal H}\vert\vert^2
+2g_{\mathcal H}({\rm grad}(\vert\vert A_{\mathcal H}\vert\vert^2),\,{\rm grad}\,H^s)}\\
\hspace{0.5truecm}\displaystyle{-\frac{2(n-1)}{3n}\cdot C_1\cdot H^s\cdot
\vert\vert\nabla^{\mathcal H}A_{\mathcal H}\vert\vert^2}\\
\hspace{0.5truecm}\displaystyle{+3C_1\cdot H^s\cdot\vert\vert A_{\mathcal H}\vert\vert^2
\left(\vert\vert A_{\mathcal H}\vert\vert^2-\frac{(H^s)^2}{n}\right)}\\
\hspace{0.5truecm}\displaystyle{+2C_1\cdot\breve C(n,C_0,\delta)\cdot\vert\vert A_{\mathcal H}\vert\vert^4
-3b(H^s)^3\cdot\vert\vert A_{\mathcal H}\vert\vert^2
+6bH^s\cdot\vert\vert{\rm grad}\,H^s\,\vert\vert^2}\\
\hspace{0.5truecm}\displaystyle{+6\vert\vert{\rm grad}\,H^s\,\vert\vert\cdot
\vert\vert{\rm grad}({\rm Tr}(({\mathcal A}^{\phi}_{\xi})^2)_{\mathcal H})\vert\vert
-\frac{3}{H^s}\cdot\vert\vert{\rm grad}\,H^s\,\vert\vert^2\cdot
{\rm Tr}(({\mathcal A}^{\phi}_{\xi})^2)_{\mathcal H}}\\
\hspace{0.5truecm}\displaystyle{
-2C_1H^s\cdot{\rm Tr}(({\mathcal A}^{\phi}_{\xi})^2)_{\mathcal H}
\cdot\left(\vert\vert A_{\mathcal H}\vert\vert^2-\frac{(H^s)^2}{n}\right)}\\
\end{array}\leqno{(7.8)}$$
\newpage
$$\begin{array}{l}
\hspace{0.5truecm}\displaystyle{-4C_1(H^s)^2\cdot
{\rm Tr}\left(({\mathcal A}^{\phi}_{\xi})^2)\circ
\left(A_{\mathcal H}-\frac{H^s}{n}\cdot{\rm id}\right)\right)}\\
\hspace{0.5truecm}\displaystyle{-2C_1H^s\cdot{\rm Tr}_{g_{\mathcal H}}^{\bullet}
{\mathcal R}\left(\left(A_{\mathcal H}-\frac{H^s}{n}\cdot{\rm id}\right)(\bullet),\bullet\right)}\\
\hspace{0.5truecm}\displaystyle{-3C_1\left(\vert\vert A_{\mathcal H}\vert\vert^2
-\frac{(H^s)^2}{n}\right)\cdot H^s\cdot
{\rm Tr}(({\mathcal A}^{\phi}_{\xi})^2)_{\mathcal H}}\\
\hspace{0.5truecm}\displaystyle{
-2C_1\cdot\breve C(n,C_0,\delta)\vert\vert A_{\mathcal H}\vert\vert^2\cdot
{\rm Tr}(({\mathcal A}^{\phi}_{\xi})^2)_{\mathcal H}}\\
\hspace{0.5truecm}\displaystyle{-4C_1\cdot\breve C(n,C_0,\delta)H^s\cdot
{\rm Tr}\left((({\mathcal A}^{\phi}_{\xi})^2)_{\mathcal H}\circ A_{\mathcal H}\right)}\\
\hspace{0.5truecm}\displaystyle{
-2C_1\cdot\breve C(n,C_0,\delta){\rm Tr}^{\bullet}_{g_{\mathcal H}}{\mathcal R}(A_{\mathcal H}\bullet,\bullet)
+9b\cdot(H^s)^3\cdot{\rm Tr}(({\mathcal A}^{\phi}_{\xi})^2)_{\mathcal H}.}
\end{array}$$
Also, in similar to $(7.5)$, we obtain 
$$\begin{array}{l}
\hspace{0.5truecm}\displaystyle{
\vert g_{\mathcal H}({\rm grad}(\vert\vert A_{\mathcal H}\vert\vert^2),\,{\rm grad}\,H^s)\vert}\\
\displaystyle{\leq 2n\sqrt{C_0}\vert\vert\nabla^{\mathcal H}A_{\mathcal H}\vert\vert^2\cdot
(H^s)^{1-\delta/2}.}
\end{array}$$
This implies together with $(7.7)$ that 
$$\vert g_{\mathcal H}({\rm grad}(\vert\vert A_{\mathcal H}\vert\vert^2),\,{\rm grad}\,H^s)\vert
\leq\left(\frac{2(n-1)}{3n}\cdot H^s+\breve C(n,C_0,\delta)\right)
\vert\vert\nabla^{\mathcal H}A_{\mathcal H}\vert\vert^2.\leqno{(7.9)}$$
Denote by $T^1V$ the unit tangent bundle of $V$.  Define a function $\Psi$ over $T^1V$ by 
$$\Psi(X):=\vert\vert\,\,d({\rm Tr}({\mathcal A}^{\phi}_X)^2)_{\widetilde{\mathcal  H}})\,\vert\vert
\quad\,\,(X\in T^1V).$$
It is clear that $\Psi$ is continuous.  
Set $\displaystyle{\widehat K_1:=\mathop{\sup}_{t\in[0,T)}\max_M\vert\vert{\rm grad}
({\rm Tr}(({\mathcal A}^{\phi}_{\xi})^2)_{\mathcal H})\vert\vert}$, which is finite because 
$\Psi$ is continuous and the closure of $\displaystyle{\mathop{\cup}_{t\in[0,T)}\phi(f_t(M))}$ is compact.  
Also, we have 
$${\rm Tr}^{\bullet}_{g_{\mathcal H}}{\mathcal R}
\left(\left(A_{\mathcal H}-\frac{H^s}{n}\cdot{\rm id}\right)(\bullet),\bullet\right)
\leq\widehat K_2\cdot\left\vert\left\vert A_{\mathcal H}-\frac{H^s}{n}\cdot{\rm id}\right\vert
\right\vert\leqno{(7.10)}$$
for some positive constant $\widehat K_2$ because of the homogeneity of $N$.  
By using $(7.7), (7.9), (7.10)$, $\vert\vert A_{\mathcal H}\vert\vert\leq H^s$, 
$\frac{1}{n}\vert\vert{\rm grad}H^s\,\vert\vert^2
\leq\vert\vert\nabla^{\mathcal H}A_{\mathcal H}\vert\vert^2$ and Proposition 6.2, 
we can derive 
$$\begin{array}{l}
\hspace{0.5truecm}\displaystyle{\frac{\partial\rho}{\partial t}-\triangle_{\mathcal H}\rho}\\
\displaystyle{\leq\left(3n+\frac{4(n-1)}{3n}-\frac{2(n-1)C_1}{3n}+6nb\right)H^s\cdot
\vert\vert\nabla^{\mathcal H}A_{\mathcal H}\vert\vert^2}\\
\hspace{0.5truecm}\displaystyle{+2\breve C(n,C_0,\delta)\cdot\vert\vert\nabla^{\mathcal H}A_{\mathcal H}\vert\vert^2
+3C_0\cdot C_1(H^s)^{5-\delta}+2C_1\cdot\breve C(n,C_0,\delta)(H^s)^4}\\
\hspace{0.5truecm}\displaystyle{-3b(H^s)^5
+6\widehat K_1\vert\vert{\rm grad}\,H^s\,\vert\vert
+\frac{3\widehat K_1}{H^s}\cdot\vert\vert{\rm grad}\,H^s\,\vert\vert^2}\\
\hspace{0.5truecm}\displaystyle{
+2C_0\cdot C_1\cdot\widehat K_1(H^s)^{3-\delta}
+4C_1\cdot\sqrt{C_0}\cdot\widehat K_1(H^s)^{3-\delta/2}}\\
\hspace{0.5truecm}\displaystyle{+2C_1\cdot\widehat K_2\cdot\sqrt{C_0}\cdot(H^s)^{2-\delta/2}
+3C_0\cdot C_1\cdot\widehat K_1\cdot(H^s)^{3-\delta}}\\
\hspace{0.5truecm}\displaystyle{+2C_1\cdot\breve C(n,C_0,\delta)\cdot\widehat K_1\cdot(H^s)^2
+4C_1\cdot\breve C(n,C_0,\delta)\cdot\widehat K_1\cdot(H^s)^3}\\
\hspace{0.5truecm}\displaystyle{+2C_1\cdot\breve C(n,C_0,\delta)\cdot\widehat K_2\cdot H^s
+9b\cdot\widehat K_1\cdot(H^s)^3.}
\end{array}\leqno{(7.11)}$$
Furthermore, by using the Young's inequality $(7.6)$ and the fact that $C_1$ is sufficiently big compared to 
$n$ and $b$, we can derive that 
$$\frac{\partial\rho}{\partial t}-\triangle_{\mathcal H}\rho
\leq C_3(n,C_0,C_1,b,\delta,\widehat K_1,\widehat K_2)$$
holds for some positive constant $C_3(n,C_0,C_1,b,\delta,\widehat K_1,\widehat K_2)$ only on 
$n,C_0,C_1,b,\delta,\widehat K_1$ and $\widehat K_2$.  
This together with $T<\infty$ implies that 
$$\begin{array}{l}
\displaystyle{
\mathop{\max}_M\,\rho_t\leq\mathop{\max}_M\,\rho_0+C_3(n,C_0,C_1,b,\delta,\widehat K_1,\widehat K_2)t}\\
\hspace{1.15truecm}\displaystyle{\leq\mathop{\max}_M\,\rho_0+C_3(n,C_0,C_1,b,\delta,\widehat K_1,\widehat K_2)
\cdot T}
\\\end{array}
$$
$(0\leq t<T)$.  Therefore, we obtain 
$$
\vert\vert{\rm grad}\,H^s\,\vert\vert^2
\leq b(H^s)^4+\mathop{\max}_M\,\rho_0\cdot H^s
+C_3(n,C_0,C_1,b,\delta,\widehat K_1,\widehat K_2)\cdot T\cdot H^s.
$$
Furthermore, by using the Young inequality $(7.6)$, we obtain 
$$
\vert\vert{\rm grad}\,H^s\,\vert\vert^2
\leq 2b(H^s)^4+C_4(n,C_0,C_1,b,\delta,\widehat K_1,\widehat K_2,T)$$
holds for some positive constant $C_4(n,C_0,C_1,b,\delta,\widehat K_1,\widehat K_2,T)$ only on 
$n,C_0,C_1,b,\delta,\widehat K_1$\newline
$\widehat K_2$ and $T$.  
Since $b$ is any positive constant and $C_4(n,C_0,C_1,b,\delta,\widehat K_1,\widehat K_2,T)$ essentially 
depends only on $n$ and $f_0$, we obtain the statement of Proposition 7.1.  \qed

\section{Proof of Theorem A.} 
In this section, we shall prove Theorem A.  
G. Huisken (\cite{Hu}) obtained the evolution inequality for the squared norm of all iterated 
covariant derivatives of the shape operators of the mean curvature flow in a complete Riemannian manifold 
satisfying curvature-pinching conditions in Theorem 1.1 of \cite{Hu}.  
See the proof of Lemma 7.2 (Page 478) of \cite{Hu} about this evolution inequality.  
In similar to this evolution inequality, we obtain the following evolution inequality.  

\vspace{0.5truecm}

\noindent
{\bf Lemma 8.1.} {\sl For any positive integer $m$, the family 
$\{\vert\vert(\nabla^{\mathcal H})^mA_{\mathcal H}\vert\vert^2\}_{t\in[0,T)}$ satisfies the following 
evolution inequality:
$$\begin{array}{l}
\hspace{0.5truecm}\displaystyle{\frac{\partial\vert\vert(\nabla^{\mathcal H})^mA_{\mathcal H}\vert\vert^2}{\partial t}
-\triangle_{\mathcal H}\vert\vert(\nabla^{\mathcal H})^mA_{\mathcal H}\vert\vert^2}\\
\displaystyle{\leq-2\vert\vert(\nabla^{\mathcal H})^mA_{\mathcal H}\vert\vert^2+C_4(n,m)}\\
\hspace{0.5truecm}\displaystyle{\times\left(\sum_{i+j+k=m}\vert\vert(\nabla^{\mathcal H})^iA_{\mathcal H}\vert\vert\cdot
\vert\vert(\nabla^{\mathcal H})^jA_{\mathcal H}\vert\vert\cdot\vert\vert(\nabla^{\mathcal H})^kA_{\mathcal H}\vert\vert\cdot
\vert\vert(\nabla^{\mathcal H})^mA_{\mathcal H}\vert\vert\right.}\\
\hspace{1truecm}\displaystyle{\left.+C_5(m)\sum_{i\leq m}\vert\vert(\nabla^{\mathcal H})^iA_{\mathcal H}\vert\vert\cdot
\vert\vert(\nabla^{\mathcal H})^mA_{\mathcal H}\vert\vert+C_6(m)\vert\vert(\nabla^{\mathcal H})^mA_{\mathcal H}\vert\vert\right),}
\end{array}\leqno{(8.1)}$$
where $C_4(n,m)$ is a positive constant depending only on $n,m$ and $C_i(m)$ ($i=5,6$) are positive constants 
depending only on $m$.}

\vspace{0.5truecm}

In similar to Corollary 12.6 of \cite{Ha}, we can derive the following interpolation inequality.  

\vspace{0.5truecm}

\noindent
{\bf Lemma 8.2.} {\sl Let $S$ be an element of $\Gamma(\pi_M^{\ast}(T^{(1,1)}M))$ such that, 
for any $t\in[0,T)$, $S_t$ is a $\mathcal G$-invariant $(1,1)$-tensor field on $M$.  For any positive integer $m$, 
the following inequality holds:
$$\int_{\overline M}\vert\vert(\nabla^{\mathcal H})^iS_{\mathcal H})\vert\vert_B^{2m/i}\,d\bar v\leq C(n,m)
\cdot\mathop{\max}_M\,\vert\vert S_{\mathcal H}\vert\vert^{2(m/i-1)}\cdot
\int_{\overline M}\vert\vert(\nabla^{\mathcal H})^mS_{\mathcal H})\vert\vert_B^2\,d\bar v,$$
where $C(n,m)$ is a positive constant depending only on $n$ and $m$.}

\vspace{0.5truecm}

From these lemmas, we can derive the following inequality.  

\vspace{0.5truecm}

\noindent
{\bf Lemma 8.3.} {\sl For any positive integer $m$, the following inequality holds:
$$\begin{array}{l}
\hspace{0.5truecm}\displaystyle{
\frac{d}{dt}\int_{\overline M}\vert\vert(\nabla^{\mathcal H})^mA_{\mathcal H}\vert\vert_B^2\,d\bar v
+2\int_{\overline M}\vert\vert(\nabla^{\mathcal H})^{m+1}A_{\mathcal H}\vert\vert_B^2\,d\bar v}\\
\displaystyle{\leq C_7(n,m,C_6(m),{\rm Vol}(M_0))
\cdot\left(\mathop{\max}_M\,\vert\vert A_{\mathcal H}\vert\vert^2+1\right)}\\
\hspace{0.5truecm}\displaystyle{\times
\left(\int_{\overline M}\vert\vert(\nabla^{\mathcal H})^mA_{\mathcal H}\vert\vert_B^2\,d\bar v
+\left(\int_{\overline M}\vert\vert(\nabla^{\mathcal H})^mA_{\mathcal H}\vert\vert_B^2\,d\bar v\right)^{1/2}\right),}
\end{array}\leqno{(8.2)}$$
where $C_7(n,m,C_6(m),{\rm Vol}(M_0))$ is a positive constant depending only on $n,m,C_6(m)$ and 
the volume ${\rm Vol}(M_0)$ of $M_0=f_0(M)$.}

\vspace{0.5truecm}

\noindent
{\it Proof.} By using $(8.1)$ and the generalized H$\ddot{\rm o}$lder inequality, we can derive 
{\small
$$\begin{array}{l}
\hspace{0.5truecm}\displaystyle{
\frac{d}{dt}\int_{\overline M}\vert\vert(\nabla^{\mathcal H})^mA_{\mathcal H}\vert\vert_B^2\,d\bar v
+2\int_{\overline M}\vert\vert(\nabla^{\mathcal H})^{m+1}A_{\mathcal H}\vert\vert_B^2\,d\bar v}\\
\displaystyle{\leq C_4(n,m)\cdot\left(\sum_{i+j+k=m}
\int_{\overline M}\vert\vert(\nabla^{\mathcal H})^iA_{\mathcal H}\vert\vert_B^{\frac{2m}{i}}\,d\bar v
\right)^{\frac{i}{2m}}
\cdot
\left(\int_{\overline M}\vert\vert(\nabla^{\mathcal H})^jA_{\mathcal H}\vert\vert_B^{\frac{2m}{j}}\,d\bar v
\right)^{\frac{j}{2m}}}\\
\hspace{1.5truecm}\displaystyle{\times
\left(\int_{\overline M}\vert\vert(\nabla^{\mathcal H})^kA_{\mathcal H}\vert\vert_B^{\frac{2m}{k}}\,d\bar v
\right)^{\frac{k}{2m}}
\cdot
\left(\int_{\overline M}\vert\vert(\nabla^{\mathcal H})^mA_{\mathcal H}\vert\vert_B^2\,d\bar v
\right)^{\frac{1}{2}}}\\
\hspace{0.5truecm}\displaystyle{+C(n,m)\widetilde C(m)\sum_{i\leq m}
\left(\int_{\overline M}\vert\vert(\nabla^{\mathcal H})^iA_{\mathcal H}\vert\vert_B^{\frac{2m}{i}}\,d\bar v
\right)^{\frac{i}{2m}}
\cdot
\left(\int_{\overline M}\vert\vert(\nabla^{\mathcal H})^mA_{\mathcal H}\vert\vert_B^{\frac{2m}{2m-i}}\,d\bar v
\right)^{\frac{2m-i}{2m}}}\\
\hspace{0.5truecm}\displaystyle{+C(n,m)\widetilde C(m+1)\cdot
\left(\int_{\overline M}\vert\vert(\nabla^{\mathcal H})^mA_{\mathcal H}\vert\vert_B^2\,d\bar v\right)^{\frac{1}{2}}
\cdot\left(\int_{\overline M}\,d\bar v\right)^{\frac{1}{2}}.}
\end{array}$$
}
From this inequality and Lemma 8.2, we can derive the desired inequality.  \qed

\vspace{0.5truecm}

From this lemma, we can derive the following statement.  

\vspace{0.5truecm}

\noindent
{\bf Proposition 8.4.} {\sl The family $\{\vert\vert A_{\mathcal H}\vert\vert^2\}_{t\in[0,T)}$ 
is not uniform bounded.}

\vspace{0.5truecm}

\noindent
{\it Proof.} Suppose that 
$\displaystyle{\mathop{\sup}_{t\in[0,T)}\mathop{\max}_M\,\vert\vert A_{\mathcal H}\vert\vert^2<\infty}$.  
Denote by $C_A$ this supremum.  Define a function $\Phi$ over $[0,T)$ by 
$$\Phi(t):=\int_{\overline M}\vert\vert(\nabla^{\mathcal H})^m(A_{\mathcal H})_t\vert\vert_B^2\,d\bar v_t
\qquad(t\in[0,T)).$$
Then, according to $(8.2)$, we have 
$$\frac{d\Phi}{dt}\leq C_7(n,m,C_6(m),{\rm Vol}(M_0))\cdot(C_A+1)\cdot(\Phi+\Phi^{1/2}).$$
Assume that $\mathop{\sup}_{t\in[0,T)}\,\Phi>1$.  Set $E:=\{t\in[0,T)\,\vert\,\Phi(t)>1\}$.  
Take any $t_0\in E$.  
Then $\Phi\geq 1$ holds over $[t_0,t_0+\varepsilon)$ for some a sufficiently small 
positive number $\varepsilon$.  Hence we have 
$$\frac{d\Phi}{dt}\leq 2C_7(n,m,C_6(m),{\rm Vol}(M_0))\cdot(C_A+1)\cdot\Phi$$
on $[t_0,t_0+\varepsilon)$.  From this inequality, we can derive 
$$\Phi(t)\leq\Phi(t_0)e^{2C_7(n,m,C_6(m),{\rm Vol}(M_0))\cdot(C_A+1)(t-t_0)}\quad\,\,(t\in[t_0,t_0+\varepsilon))$$
and hence 
$$\Phi(t)\leq\Phi(t_0)e^{2C_7(n,m,C_6(m),{\rm Vol}(M_0))\cdot(C_A+1)T}\quad\,\,(t\in[t_0,t_0+\varepsilon)).$$
This fact together with the arbitrariness of $t_0$ implies that $\Phi_t$ is uniform bounded.  
Thus, we see that 
$$\mathop{\sup}_{t\in[0,T)}\,\int_{\overline M}\vert\vert(\nabla^{\mathcal H})^m(A_{\mathcal H})_t\vert\vert_B^2\,
d\bar v_t<\infty$$
holds in general.  
Furthermore, since this inequality holds for any positive integar $m$, it follows from Lemma 8.2 that 
$$\mathop{\sup}_{t\in[0,T)}\int_{\overline M}\vert\vert(\nabla^{\mathcal H})^m(A_{\mathcal H})_t\vert\vert_B^{\it l}\,
d\bar v_t<\infty$$
holds for any positive integar $m$ and any positive constant ${\it l}$.  
Hence, by the Sobolev's embedding theorem, we obtain 
$$\mathop{\sup}_{t\in[0,T)}\max_M\,\vert\vert(\nabla^{\mathcal H})^m(A_{\mathcal H})_t\vert\vert<\infty.$$
Since this fact holds for any positive integer $m$, $f_t$ converges to a $C^{\infty}$-embedding $f_T$ as 
$t\to T$ in $C^{\infty}$-topology.  This implies that the mean curvature flow $f_t$ extends after $T$ 
because of the short time existence of the mean curvature flow starting from $f_T$.  
This contradicts the definition of $T$.  Therefore we obtain 
$$\mathop{\sup}_{t\in[0,T)}\mathop{\max}_M\,\vert\vert A_{\mathcal H}\vert\vert^2=\infty.$$
\qed

\vspace{0.5truecm}

By imitating the proof of Theorem 4.1 of \cite{A1,A2}, we can show the folowing fact, where we note that 
more general curvature flows (including mean curvature flows as special case) is treated in \cite{A1,A2}.  

\vspace{0.5truecm}

\noindent
{\bf Lemma 8.5.} {\sl The following uniform boundedness holds:
$$\mathop{\inf}_{t\in[0,T)}\max\,
\{\varepsilon>0\,\vert\,(A_{\mathcal H})_t\geq\varepsilon H^s_t\cdot{\rm id}\,\,{\rm on}\,\,M\}>0$$
and hence 
$$\mathop{\sup}_{(x,t)\in M\times[0,T)}\frac{\lambda_{\max}(x,t)}{\lambda_{\min}(x,t)}
\leq\frac{1}{\varepsilon_0},$$
where $\lambda_{\max}(x,t)$ (resp. $\lambda_{\min}(x,t)$) denotes 
the maximum (resp. minimum) eigenvalue of $(A_{\mathcal H})_{(x,t)}$ and $\varepsilon_0$ denotes 
the above infimum.
}

\vspace{0.5truecm}

\noindent
{\it Proof.} Since 
$$\begin{array}{l}
\hspace{0.5truecm}\displaystyle{\left(\frac{\partial h_{\mathcal H}}{\partial t}-\triangle_{\mathcal H}^{\mathcal H}
h_{\mathcal H}\right)(X,Y)}\\
\displaystyle{=-2H^s\cdot h_{\mathcal H}(A_{\mathcal H}(X),Y)
+g_{\mathcal H}
\left(\left(\frac{\partial A_{\mathcal H}}{\partial t}-\triangle_{\mathcal H}^{\mathcal H}A_{\mathcal H}\right)(X),Y\right).}
\end{array}$$
From this relation, Lemmas 4.6 and 4.9, we can derive 
$$\begin{array}{l}
\hspace{0.5truecm}\displaystyle{\frac{\partial A_{\mathcal H}}{\partial t}-\triangle_{\mathcal H}^{\mathcal H}A_{\mathcal H}}\\
\displaystyle{=-2H^s(({\mathcal A}^{\phi}_{\xi})^2)_{\mathcal H}
+{\rm Tr}\left((A_{\mathcal H})^2-(({\mathcal A}^{\phi})^2)_{\mathcal H}\right)\cdot A_{\mathcal H}-{\mathcal R}^{\sharp}.}
\end{array}$$
Furthermore, from this evolution equation and Lemma 4.10, we can derive 
$$\begin{array}{l}
\hspace{0.5truecm}\displaystyle{\frac{\partial}{\partial t}\left(\frac{A_{\mathcal H}}{H^s}\right)
-\triangle_{\mathcal H}^{\mathcal H}\left(\frac{A_{\mathcal H}}{H^s}\right)}\\
\displaystyle{=\frac{1}{H^s}\nabla^{\mathcal H}_{{\rm grad}\,H^s}
\left(\frac{A_{\mathcal H}}{H^s}\right)
+\frac{\vert\vert{\rm grad}\,H^s\,\vert\vert^3}{(H^s)^3}\cdot A_{\mathcal H}
-2(({\mathcal A}^{\phi}_{\xi})^2)_{\mathcal H}}\\
\hspace{0.5truecm}\displaystyle{+\frac{2}{H^s}\cdot
{\rm Tr}(({\mathcal A}^{\phi}_{\xi})^2)_{\mathcal H}\cdot A_{\mathcal H}-\frac{1}{H^s}{\mathcal R}^{\sharp}.}
\end{array}$$
For simplicity, we set 
$$S_{\mathcal H}:=g_{\mathcal H}\left(\frac{1}{H^s}A_{\mathcal H}(\bullet),\bullet\right)$$
and 
$$\begin{array}{l}
\displaystyle{P(S)_{\mathcal H}
:=\frac{\vert\vert{\rm grad}\,H^s\,\vert\vert^3}{(H^s)^3}\cdot h_{\mathcal H}
-2((({\mathcal A}^{\phi}_{\xi})^2)_{\mathcal H})_{\flat}}\\
\hspace{2truecm}\displaystyle{
+\frac{2}{H^s}\cdot{\rm Tr}(({\mathcal A}^{\phi}_{\xi})^2)_{\mathcal H}\cdot h_{\mathcal H}
-\frac{1}{H^s}{\mathcal R},}
\end{array}$$
where $((({\mathcal A}^{\phi}_{\xi})^2)_{\mathcal H})_{\flat}$ is defined by 
$((({\mathcal A}^{\phi}_{\xi})^2)_{\mathcal H})_{\flat}(\bullet.\bullet)
:=g_{\mathcal H}((({\mathcal A}^{\phi}_{\xi})^2)_{\mathcal H}(\bullet),\bullet)$.  
Also, set 
$$\varepsilon_0:=\max\{\varepsilon>0\,\vert\,(S_{\mathcal H})_0\geq\varepsilon g_{\mathcal H}\}.$$
Then, for any $(x,t)\in M\times[0,T)$, any $\varepsilon>0$ and any 
$X\in{\rm Ker}(S_{\mathcal H}+\varepsilon g_{\mathcal H})_{(x,t)}$, we can show 
$P(S_{\mathcal H}+\varepsilon g_{\mathcal H})_{(x,t)}(X,X)\geq 0$.  
Hence, by the maximum principle (Theorem 5.1 of \cite{Koi2}), we can derive that 
$(S_{\mathcal H})_t\geq\varepsilon_0 g_{\mathcal H}$, that is, 
$(A_{\mathcal H})_t\geq\varepsilon_0H^s_t g_{\mathcal H}$ holds for all $t\in[0,T)$.  
From this fact, it follows that 
$\lambda_{\min}(x,t)\geq\varepsilon_0\vert\vert H_{(x,t)}\vert\vert$ 
holds for all $(x,t)\in M\times[0,T)$.  Hence we obtain 
$$\mathop{\sup}_{(x,t)\in M\times[0,T)}\frac{\lambda_{\max}(x,t)}{\lambda_{\min}(x,t)}
\leq\mathop{\sup}_{(x,t)\in M\times[0,T)}\frac{\lambda_{\max}(x,t)}{\varepsilon_0\vert\vert H_{(x,t)}\vert\vert}
\leq\frac{1}{\varepsilon_0}.$$
\qed

\vspace{0.5truecm}

According to this lemma, we see that such a case as in Figure 8.1 does not happen.  


\vspace{0.25truecm}

{\small
\centerline{
\unitlength 0.1in
\begin{picture}( 35.4000, 11.0000)( 12.0000,-16.9000)
%
\special{pn 8}%
\special{pa 1600 1000}%
\special{pa 2190 1000}%
\special{fp}%
%
\special{pn 8}%
\special{ar 1600 1170 60 170  3.2081608 4.7123890}%
%
\special{pn 8}%
\special{ar 2180 1180 60 170  3.2081608 4.7123890}%
%
\special{pn 8}%
\special{ar 1870 1290 670 320  5.1826283 6.2831853}%
\special{ar 1870 1290 670 320  0.0000000 4.2669212}%
%
\special{pn 8}%
\special{ar 1850 1010 800 150  0.7965225 2.3009202}%
%
\special{pn 8}%
\special{ar 1870 1310 670 100  6.2831853 6.2831853}%
\special{ar 1870 1310 670 100  0.0000000 3.1415927}%
%
\special{pn 8}%
\special{ar 1870 1320 670 100  3.1415927 3.1727615}%
\special{ar 1870 1320 670 100  3.2662680 3.2974368}%
\special{ar 1870 1320 670 100  3.3909433 3.4221121}%
\special{ar 1870 1320 670 100  3.5156186 3.5467875}%
\special{ar 1870 1320 670 100  3.6402940 3.6714628}%
\special{ar 1870 1320 670 100  3.7649693 3.7961381}%
\special{ar 1870 1320 670 100  3.8896446 3.9208134}%
\special{ar 1870 1320 670 100  4.0143199 4.0454888}%
\special{ar 1870 1320 670 100  4.1389953 4.1701641}%
\special{ar 1870 1320 670 100  4.2636706 4.2948394}%
\special{ar 1870 1320 670 100  4.3883459 4.4195147}%
\special{ar 1870 1320 670 100  4.5130212 4.5441901}%
\special{ar 1870 1320 670 100  4.6376965 4.6688654}%
\special{ar 1870 1320 670 100  4.7623719 4.7935407}%
\special{ar 1870 1320 670 100  4.8870472 4.9182160}%
\special{ar 1870 1320 670 100  5.0117225 5.0428914}%
\special{ar 1870 1320 670 100  5.1363978 5.1675667}%
\special{ar 1870 1320 670 100  5.2610732 5.2922420}%
\special{ar 1870 1320 670 100  5.3857485 5.4169173}%
\special{ar 1870 1320 670 100  5.5104238 5.5415927}%
\special{ar 1870 1320 670 100  5.6350991 5.6662680}%
\special{ar 1870 1320 670 100  5.7597745 5.7909433}%
\special{ar 1870 1320 670 100  5.8844498 5.9156186}%
\special{ar 1870 1320 670 100  6.0091251 6.0402940}%
\special{ar 1870 1320 670 100  6.1338004 6.1649693}%
\special{ar 1870 1320 670 100  6.2584758 6.2831853}%
%
\special{pn 8}%
\special{ar 1590 1220 60 210  4.7123890 4.8012779}%
\special{ar 1590 1220 60 210  5.0679445 5.1568334}%
\special{ar 1590 1220 60 210  5.4235001 5.5123890}%
%
\special{pn 8}%
\special{ar 2180 1220 60 210  4.7123890 4.8012779}%
\special{ar 2180 1220 60 210  5.0679445 5.1568334}%
\special{ar 2180 1220 60 210  5.4235001 5.5123890}%
%
\special{pn 8}%
\special{ar 1860 1200 760 130  3.9197094 3.9466757}%
\special{ar 1860 1200 760 130  4.0275746 4.0545409}%
\special{ar 1860 1200 760 130  4.1354397 4.1624060}%
\special{ar 1860 1200 760 130  4.2433049 4.2702712}%
\special{ar 1860 1200 760 130  4.3511701 4.3781364}%
\special{ar 1860 1200 760 130  4.4590352 4.4860015}%
\special{ar 1860 1200 760 130  4.5669004 4.5938667}%
\special{ar 1860 1200 760 130  4.6747656 4.7017319}%
\special{ar 1860 1200 760 130  4.7826307 4.8095970}%
\special{ar 1860 1200 760 130  4.8904959 4.9174622}%
\special{ar 1860 1200 760 130  4.9983611 5.0253274}%
\special{ar 1860 1200 760 130  5.1062262 5.1331925}%
\special{ar 1860 1200 760 130  5.2140914 5.2410577}%
\special{ar 1860 1200 760 130  5.3219566 5.3489229}%
\special{ar 1860 1200 760 130  5.4298218 5.4567880}%
\special{ar 1860 1200 760 130  5.5376869 5.5439500}%
%
\special{pn 8}%
\special{pa 2810 1140}%
\special{pa 3030 1140}%
\special{fp}%
\special{sh 1}%
\special{pa 3030 1140}%
\special{pa 2964 1120}%
\special{pa 2978 1140}%
\special{pa 2964 1160}%
\special{pa 3030 1140}%
\special{fp}%
%
\special{pn 8}%
\special{pa 3348 1040}%
\special{pa 3538 1040}%
\special{fp}%
%
\special{pn 8}%
\special{ar 3348 1104 20 64  3.2160814 4.7123890}%
%
\special{pn 8}%
\special{ar 3534 1108 20 66  3.2160814 4.7123890}%
%
\special{pn 8}%
\special{ar 3436 1150 216 122  5.1851897 6.2831853}%
\special{ar 3436 1150 216 122  0.0000000 4.2657852}%
%
\special{pn 8}%
\special{ar 3430 1044 256 58  0.7944302 2.3002306}%
%
\special{pn 8}%
\special{ar 3436 1158 216 38  6.2831853 6.2831853}%
\special{ar 3436 1158 216 38  0.0000000 3.1415927}%
%
\special{pn 8}%
\special{ar 3436 1162 216 38  3.1415927 3.2364543}%
\special{ar 3436 1162 216 38  3.5210393 3.6159010}%
\special{ar 3436 1162 216 38  3.9004859 3.9953476}%
\special{ar 3436 1162 216 38  4.2799326 4.3747942}%
\special{ar 3436 1162 216 38  4.6593792 4.7542409}%
\special{ar 3436 1162 216 38  5.0388259 5.1336875}%
\special{ar 3436 1162 216 38  5.4182725 5.5131342}%
\special{ar 3436 1162 216 38  5.7977191 5.8925808}%
\special{ar 3436 1162 216 38  6.1771658 6.2720274}%
%
\special{pn 8}%
\special{ar 3346 1124 20 80  4.7123890 4.9548132}%
\special{ar 3346 1124 20 80  5.6820860 5.7031825}%
%
\special{pn 8}%
\special{ar 3534 1124 20 80  4.7123890 4.9548132}%
\special{ar 3534 1124 20 80  5.6820860 5.7031825}%
%
\special{pn 8}%
\special{ar 3432 1116 244 50  3.9325982 4.0145095}%
\special{ar 3432 1116 244 50  4.2602433 4.3421545}%
\special{ar 3432 1116 244 50  4.5878883 4.6697996}%
\special{ar 3432 1116 244 50  4.9155334 4.9974447}%
\special{ar 3432 1116 244 50  5.2431784 5.3250897}%
%
\special{pn 8}%
\special{pa 3870 1140}%
\special{pa 4390 1140}%
\special{fp}%
\special{sh 1}%
\special{pa 4390 1140}%
\special{pa 4324 1120}%
\special{pa 4338 1140}%
\special{pa 4324 1160}%
\special{pa 4390 1140}%
\special{fp}%
\put(39.4000,-12.4000){\makebox(0,0)[lt]{$t\to T$}}%
%
\special{pn 8}%
\special{pa 4650 1070}%
\special{pa 4740 1070}%
\special{fp}%
\put(17.3000,-16.9000){\makebox(0,0)[lt]{$\overline M_0$}}%
\put(33.3000,-13.9000){\makebox(0,0)[lt]{$\overline M_t$}}%
%
\special{pn 8}%
\special{pa 1970 830}%
\special{pa 1850 1040}%
\special{dt 0.045}%
\special{sh 1}%
\special{pa 1850 1040}%
\special{pa 1900 992}%
\special{pa 1876 994}%
\special{pa 1866 972}%
\special{pa 1850 1040}%
\special{fp}%
\put(14.8000,-7.6000){\makebox(0,0)[lb]{almost cylindrical part}}%
\end{picture}%
\hspace{0.5truecm}}
}

\vspace{0.5truecm}

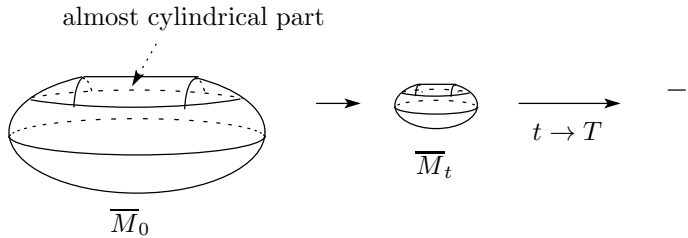
\centerline{{\bf Figure 8.1$\,\,:\,\,$ The case where $\lim\limits_{t\to T}\overline M_t$ is not a round point}}

\vspace{0.5truecm}

By using Proposition 8.4 and Lemma 8.5, we shall prove the statement (i) of Theorem A.  

\vspace{0.5truecm}

\noindent
{\it Proof of (i) of Theorem A.} 
According to Proposition 8.4 and Lemma 8.5, we have 
$$\lim_{t\to T}\min_{x\in M}\,\lambda_{\min}(x,t)=\infty.$$
Set $\Lambda_{\min}(t):=\min_{x\in M}\,\lambda_{\min}(x,t)$.  
Let $x_{\min}(t)$ be a point of $\overline M$ with $\lambda_{\min}(x_{\min}(t),t)$\newline
$=\Lambda_{\min}(t)$ and set $\bar x_{\min}(t):=\phi_M(x_{\min}(t))$.  
Denote by $\gamma_{\bar f_t(\bar x_{\min}(t))}$ the normal geodesic of $\bar f_t(\overline M)$ starting from 
$\bar f_t(\bar x_{\min}(t))$.  Set $p_t:=\gamma_{\bar f_t(\bar x_{\min}(t))}(1/\Lambda_{\min}(t))$.  
Since $N$ is of non-negative curvature, 
the focal radii of $\overline M_t$ along any normal geodesic are smaller than or equal to 
$\frac{1}{\Lambda_{\min}(t)}$.  This implies that $\bar f_t(\overline M)$ is included by the geodesic sphere of 
radius $\frac{1}{\Lambda_{\min}(t)}$ centered at $p_t$ in $N$.  Hence, since 
$\lim\limits_{t\to T}\frac{1}{\Lambda_{\min}(t)}=0$, we see that, as $t\to T$, 
$\overline M_t$ collapses to a one-point set, that is, $M_t$ collapses to a $\mathcal G$-orbit.  \qed

\vspace{0.5truecm}

Denote by $({\rm Ric}_{\overline M})_t$ the Ricci tensor of $\overline g_t$ and let 
${\rm Ric}_{\overline M}$ be the element of \newline
$\Gamma(\pi_{\overline M}^{\ast}(T^{(0,2)}\overline M))$ defined by $({\rm Ric}_{\overline M})_t$'s.  
To show the statement (ii) of Theorem A, we prepare the following some lemmas.  

\vspace{0.3truecm}

\noindent
{\bf Lemma 8.6.} {\sl {\rm (i)} For the section ${\rm Ric}_{\overline M}$, the following relation holds:
$${\rm Ric}_{\overline M}(X,Y)=-3{\rm Tr}({\mathcal A}^{\phi}_{X^L}\circ{\mathcal A}^{\phi}_{Y^L})_{\mathcal H}
-\overline g(\overline A^2X,Y)+\vert\vert\overline H\vert\vert\cdot\overline g(\overline AX,Y)\leqno{(8.3)}$$
($X,Y\in\Gamma(\pi_{\overline M}^{\ast}(T\overline M))$), where 
where $X^L$ (resp. $Y^L$) is the horizontal lift of $X$ (resp. $Y$) to $V$.

{\rm (ii)} Let $\lambda_1$ be the smallest eigenvalue of $\overline A_{(x,t)}$.  Then we have 
$$({\rm Ric}_{\overline M})_{(x,t)}(v,v)\geq(n-1)\lambda_1^2\overline g_{(x,t)}(v,v)
\quad\,\,(v\in T_x\overline M).\leqno{(8.4)}$$
}

\vspace{0.3truecm}

\noindent
{\it Proof.} Denote by $\overline{\rm Ric}$ the Ricci tensor of $N$.  
By the Gauss equation, we have 
$$\begin{array}{r}
\displaystyle{{\rm Ric}_{\overline M}(X,Y)=\overline{\rm Ric}(X,Y)-\overline g(\overline A^2X,Y)
+\vert\vert\overline H\vert\vert\overline g(\overline AX,Y)-\overline R(\xi,X,Y,\xi)}\\
(X,Y\in T\overline M).
\end{array}$$
Also, by a simple calculation, we have 
$$\overline{\rm Ric}(X,Y)=-3{\rm Tr}({\mathcal A}^{\phi}_{X^L}\circ{\mathcal A}^{\phi}_{Y^L})_{\mathcal H}
+3g_{\mathcal H}(({\mathcal A}^{\phi}_{X^L}\circ{\mathcal A}^{\phi}_{Y^L})(\xi),\,\xi)$$
and 
$$\overline R(\xi,X,Y,\xi)=3g_{\mathcal H}(({\mathcal A}^{\phi}_{X^L}\circ{\mathcal A}^{\phi}_{Y^L})(\xi),\,\xi)$$
$(X,Y\in\Gamma(\pi_{\overline M}^{\ast}(T\overline M)))$.  
From these relations, we obtain the relation $(8.3)$.  

Next we show the inequality in the statement (ii).  
Since ${\mathcal A}^{\phi}_{v^L}$ is skew-symmetric, we have ${\rm Tr}(({\mathcal A}^{\phi}_{v^L})^2)\leq 0$.  
Also we have 
$$-\overline g_{(x,t)}(\overline A_{(x,t)}^2(v),v)+\vert\vert\overline H_{(x,t)}\vert\vert
\cdot\overline g_{(x,t)}(\overline A_{(x,t)}(v),v)
\geq(n-1)\lambda_1^2\cdot\overline g_{(x,t)}(v,v).$$
Hence, from the relation in (i), we can derive the inequality $(8.4)$.  \qed

\vspace{0.3truecm}

By remarking the behavior of geodesic rays reaching the singular set of a compact Riemannian orbifold (see Figure 8.2) 
and using the discussion in the proof of Myers's theorem (\cite{M}), we can show the following Myers-type theorem for 
Riemannian orbifolds.  

\vspace{0.3truecm}

\noindent
{\bf Theorem 8.7.} {\sl 
Let $(N,g)$ be an $n$-dimensional compact (connected) Riemannian orbifold.  
If its Ricci curvature $Ric$ of $(N,g)$ satisfies $Ric\geq(n-1)K$ for some positive constant $K$, then 
the first conjugate radius along any geodesic in $(N,g)$ is smaller than or equal to 
$\displaystyle{\frac{\pi}{\sqrt{K}}}$ and hence so is also the diameter of $(N,g)$.}

\vspace{0.3truecm}

By using Propositions 7.1, 8.4, Lemmas 8.6 and Theorem 8.7, we prove the statement (ii) of Theorem A.  

\vspace{0.3truecm}

{\small
\centerline{
\unitlength 0.1in
%
\hspace{3.25truecm}}
}

\vspace{0.3truecm}

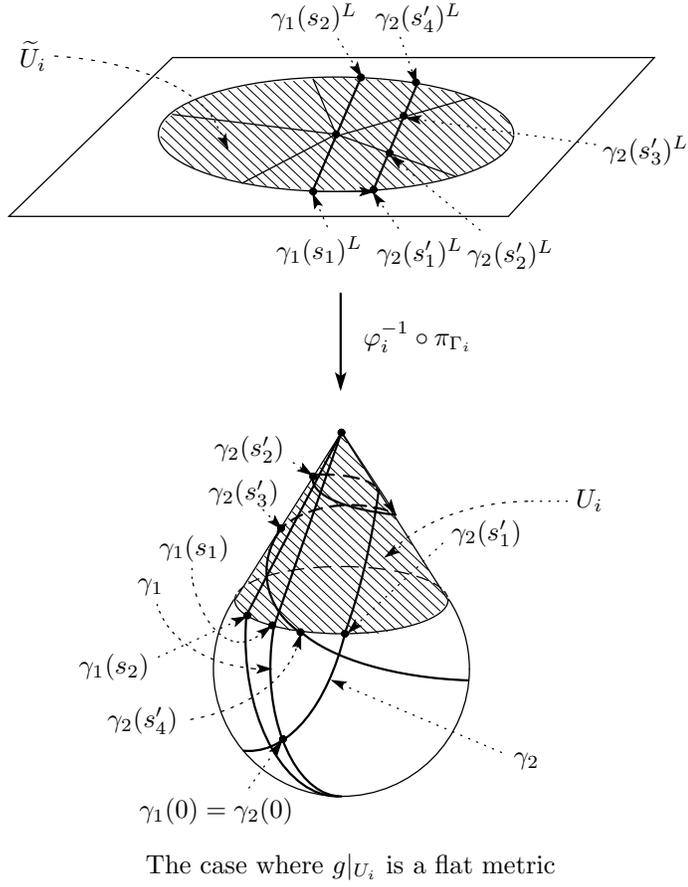
\centerline{{\bf Figure 8.2$\,\,:\,\,$ The behavior of geodesic rays around the singular set}}

\vspace{0.5truecm}

\noindent
{\it Proof of (ii) of Theorem A.} 
(Step I) According to Proposition 7.1, 
for any positive constant $b$, there exists a constant $C(b,f_0)$ (depending only on $b$ and $f_0$) satisfying 
$$\vert\vert{\rm grad}\,H^s\,\vert\vert^2\leq b\cdot(H^s)^4
+C(b,f_0)\,\,\,\,{\rm on}\,\,\,\,M\times[0,T).$$
According to Proposition 8.4, we have $\lim_{t\to T}(H^s_t)_{\max}=\infty$.  
Hence there exists a positive constant $t(b)$ with 
$\displaystyle{(H^s_t)_{\max}\geq\left(\frac{C(b,f_0)}{b}\right)^{1/4}}$ for any 
$t\in[t(b),T)$.  Then we have 
$$\vert\vert{\rm grad}\,H^s_t\,\vert\vert\leq\sqrt{2b}(H^s_t)_{\max}^2
\leqno{(8.5)}$$
for any $t\in[t(b),T)$.  Fix $t_0\in[t(b),T)$.  Let $x_{t_0}$ be a maximal point of 
$\vert\vert H_{t_0}\vert\vert$.  
Take any geodesic $\gamma$ of length $\frac{1}{\sqrt 2\vert\vert H_{t_0}\vert\vert_{\max}\cdot b^{1/4}}$ 
starting from $x_{t_0}$.  According to $(8.5)$, we have 
$$\vert\vert H_{t_0}\vert\vert\geq(1-b^{1/4})\vert\vert H_{t_0}\vert\vert_{\max}$$
along $\gamma$.  From the arbitrariness of $t_0$, this fact holds for any $t\in[t(b),T)$.  

(Step II) 
For any $x\in \overline M$, denote by $\gamma_{\overline f_t(x)}$ the normal geodesic of 
$\overline f_t(\overline M)$ starting from $\overline f_t(x)$.  
Set $p_t:=\gamma_{\overline f_t(x)}\left(1/\lambda_{\min}(x,t)\right)$ and 
$q_t(s):=\gamma_{\overline f_t(x)}\left(s/\lambda_{\max}(x,t)\right)$.  
Since $N$ is of non-negative curvature, the focal radii of $\overline f_t(M)$ at $x$ are smaller than or equal to 
$1/\lambda_{\min}(x,t)$.  
Denote by $G_2(TN)$ the Grassmann bundle of $N$ of $2$-planes and ${\rm Sec}:G_2(TN)\to{\Bbb R}$ 
the function defined by assigning the sectional curvature of $\Pi$ to each element $\Pi$ of $G_2(TN)$.  
Since $\displaystyle{\overline{\mathop{\cup}_{t\in[0,T)}\overline f_t(\overline M)}}$ is compact, 
there exists the maximum of ${\rm Sec}$ over 
$\displaystyle{\overline{\mathop{\cup}_{t\in[0,T)}\overline f_t(\overline M)}}$.  
Denote by $\kappa_{\max}$ this maximum.  It is easy to show that the focal radii of $\overline f_t(\overline M)$ 
at $x$ are bigger than or equal to $\widehat c/\lambda_{\max}(x,t)$ for some positive constant $\widehat c$ 
depending only on $\kappa_{\max}$.  
Hence a sufficiently small neighborhood of $\overline f_t(x)$ in $\overline f_t(\overline M)$ 
is included by the closed domain surrounded by the geodesic spheres of radius $1/\lambda_{\min}(x,t)$ centered at 
$p_t$ and that of radius $\widehat c/\lambda_{\max}(x,t)$ centered at $q_t(\widehat c)$.  
On the other hand, according to Lemma 8.5, we have 
$$\mathop{\sup}_{(x,t)\in M\times[0,T)}\frac{\lambda_{\max}(x,t)}{\lambda_{\min}(x,t)}<\infty.$$
By using these facts, we can show 
$$\mathop{\sup}_{t\in[0,T)}\,\frac{(H^s_t)_{\max}}{(H^s_t)_{\min}}<\infty$$
and 
$$\mathop{\inf}_{t\in[0,T)}\max\,\{\varepsilon>0\,\vert\,
(A_{\mathcal H})_t\geq\varepsilon H^s_t\cdot{\rm id}\,\,{\rm on}\,\,M\}>0.$$
Set 
$$C_0:=\mathop{\sup}_{t\in[0,T)}\,
\frac{(H^s_t)_{\max}}{(H^s_t)_{\min}}$$
and 
$$\varepsilon_0:=\mathop{\inf}_{t\in[0,T)}\max\,
\{\varepsilon>0\,\vert\,(A_{\mathcal H})_t\geq\varepsilon H^s_t\cdot{\rm id}\,\,{\rm on}\,\,M\}.$$
Then, since $A_{\mathcal H}\geq\varepsilon_0H^s_{\min}\cdot{\rm id}$ on $M\times[0,T)$, 
it follows from (ii) of Lemma 8.6 that 
$$({\rm Ric}_{\overline M})_{(x,t)}(v,v)\geq(n-1)\varepsilon_0^2\cdot(H^s_t)_{\min}^2
\cdot\overline g_{(x,t)}(v,v)$$
for any $(x,t)\in M\times[0,T)$ and any $v\in T_x\overline M$.  
Hence, according to Theorem 8.7, the first conjugate radius along any geodesic $\gamma$ in 
$(\overline M,\overline g_t)$ is smaller than or equal to 
$\frac{\pi}{\varepsilon_0(H^s_t)_{\min}}$ for any $t\in[0,T)$.  
This implies that 
$$\exp_{\overline f_t(x)}\left(B_{\overline f_t(x)}
\left(\frac{\pi}{\varepsilon_0(H^s_t)_{\min}}\right)\right)=\overline M$$ 
holds for any $t\in[0,T)$, where 
$\exp_{f_t(x)}$ denotes the exponential map of $(\overline M,\overline g_t)$ at $\overline f_t(x)$ 
and $\displaystyle{B_{\overline f_t(x)}\left(\frac{\pi}{\varepsilon_0(H^s_t)_{\min}}\right)}$ 
denotes the closed ball of radius $\frac{\pi}{\varepsilon_0(H^s_t)_{\min}}$ in 
$T_{\overline f_t(x)}\overline M$ centered at the zero vector {\bf 0}.  
By the arbitrariness of $b$ (in (Step I)), we may assume that 
$b\leq\frac{\varepsilon_0^4}{4\pi^4C_0^4}$.  Then we have 
$$\frac{1}{\sqrt 2(H^s_t)_{\max}\cdot b^{1/4}}
\geq\frac{\pi}{\varepsilon_0(H^s_t)_{\min}}$$
$(t\in[0,T)$).  Let $t_0$ be as in Step I.  Then it follows from the above facts that 
$$\vert\vert H_{t_0}\vert\vert\geq(1-b^{1/4})\vert\vert H_{t_0}\vert\vert_{\max}$$ 
holds on $\overline M$.  From the arbitariness of $t_0$, it follows that 
$$H^s\geq(1-b^{1/4})H^s_{\max}$$ 
holds on $\overline M\times[t(b),T)$.  In particular, we obtain 
$$\frac{H^s_{\max}}{H^s_{\min}}\leq\frac{1}{1-b^{1/4}}$$
on $\overline M\times[t(b),T)$.  Therefore, by approaching $b$ to $0$, we can derive 
$$\lim_{t\to T}\frac{(H^s_t)_{\max}}{(H^s_t)_{\min}}=1.$$
\qed

\vspace{1truecm}

\section{Towards application to the Guage theory} 
In this section, we shall state the vision for applying the study of regularized mean curvature flows to the Gauge theory.  
In the future, we plan to find interesting Riemannian submanifolds and interesting flows (of Riemannian submanifolds) in the Yang-Milles moduli space 
or the self-dual moduli space.  To state the strategy of this plan, we first recall some basic notions in the theory of the connections of 
the principal bundles.  
Let $\pi:P\to B$ be a principal bundle over a compact manifold $B$ having a compact semi-simple Lie group $G$ 
as the structure group.  Fix an ${\rm Ad}(G)$-invariant inner product $\langle\,\,,\,\,\rangle_{\mathfrak g}$ 
(for example, the $(-1)$-multiple of the Killing form) of the Lie algebra $\mathfrak g$ of $G$, where ${\rm Ad}$ 
denotes the adjoint representation of $G$.  Denote by $g_G$ the bi-invariant metric of $G$ induced from $\langle\,\,,\,\,\rangle_{\mathfrak g}$.  
Set 
$$\Omega^{\infty}_{\mathcal T,1}(P,\mathfrak g):=
\{\widehat A\in\Omega_1^{\infty}(P,\mathfrak g)\,\vert\,
R_g^{\ast}\widehat A={\rm Ad}(g^{-1})\circ\widehat A\,\,\,(\forall\,g\in G)\,\,\,\,\&\,\,\,\,
\widehat A|_{\mathcal V}=0\},$$
where $\mathcal V$ denotes the vertical distribution of the bundle $P$.  
Each element of $\Omega^{\infty}_{\mathcal T,1}(P,\mathfrak g)$ is called 
a {\it $\mathfrak g$-valued tensorial $1$-form of class $C^{\infty}$} on $P$.  
Also, let 
$\Omega_1^{\infty}(B,{\rm Ad}(P))(=\Gamma^{\infty}(T^{\ast}B\otimes{\rm Ad}(P)))$ be the space of 
all ${\rm Ad}(P)$-valued $1$-forms of class $C^{\infty}$ over $B$, where ${\rm Ad}(P)$ denotes the adjoint bundle 
$P\times_{\rm Ad}\mathfrak g$.  
The space ${\mathcal A}_P^{\infty}$ is the affine space having 
$\Omega^{\infty}_{\mathcal T,1}(P,\mathfrak g)$ as the associated vector space.  
Furthermore, $\Omega^{\infty}_{\mathcal T,1}(P,\mathfrak g)$ is 
identified with $\Omega_1^{\infty}(B,{\rm Ad}(P))$ under the correspondence 
$\widehat A\leftrightarrow A$ defined by 
$u\cdot\widehat A_u(X)=A_{\pi(u)}(\pi_{\ast}X)$ ($u\in P,\,\,X\in T_uP$).  

Denote by $\mathcal A_P^{w,s}$ the space of all $s$-times weak differentiable connections of $P$ and 
$\Omega_i^{w,s}(B,{\rm Ad}(P))$ the space of all $s$-times weak differentiable ${\rm Ad}(P)$-valued $i$-form on $P$.  
Fix a $C^{\infty}$-connection $\omega_0$ of $P$.  Define an operator 
$\square_{\omega_0}:\Omega_i^{w,s}(B,{\rm Ad}(P))\to\Omega_i^{w,s-2}(B,{\rm Ad}(P))$ by 
$$
\square_{\omega_0}:=\left\{\begin{array}{ll}
d_{\omega_0}\circ d_{\omega_0}^{\ast}+d_{\omega_0}^{\ast}\circ d_{\omega_0}+{\rm id} & (i\geq 1)\\
d_{\omega_0}^{\ast}\circ d_{\omega_0}+{\rm id} & (i=0),
\end{array}\right.$$
where $d_{\omega_0}$ denotes the covariant exterior derivative with respect to $\omega_0$ and 
$d_{\omega_0}^{\ast}$ deontes the adjoint operator of $d_{\omega_0}$ with respect to the $L^0$-inner products of 
$\Omega_i^{w,j}(B,{\rm Ad}(P))$ ($j\geq0$).  
The $H^s$-inner product $\langle\,\,,\,\,\rangle^{\omega_0}_s$ of 
$T_{\omega}\mathcal A_P^{w,s}(\approx\Omega_1^{w,s}(B,{\rm Ad}(P))\approx\Gamma^{w,s}(T^{\ast}B\otimes{\rm Ad}(P)))$ 
is defined by 
$$\begin{array}{r}
\displaystyle{\langle A_1,A_2\rangle^{\omega_0}_s
:=\int_{x\in B}\langle(A_1)_x,(\square_{\omega_0}^s(A_2))_x\rangle_{B,\mathfrak g}\,dv_B}\\
\displaystyle{(A_1,A_2\in\Omega_1^{w,s}(B,{\rm Ad}(P))),}
\end{array}\leqno{(9.1)}$$
where $\square_{\omega_0}^s(A_2)$ denotes the element of $\Omega_1^{w,0}(B,{\rm Ad}(P))$ corresponding to 
$\square^s(\widehat A_2)$, $\langle\,\,,\,\,\rangle_{B,\mathfrak g}$ denotes the fibre metric of 
$T^{\ast}B\otimes{\rm Ad}(P)$ defined by the the Riemannian metric of $B$ and 
$\langle\,\,,\,\,\rangle_{\mathfrak g}$ and $dv_B$ denotes the volume element of the Riemannian metric of $B$.  
Let $\Omega_1^{H^s}(B,{\rm Ad}(P))$ be the completion of $\Omega_1^{\infty}(B,{\rm Ad}(P))$ with respect to 
$\langle\,\,,\,\,\rangle^{\omega_0}_s$, that is, 
$$\Omega_1^{H^s}(B,{\rm Ad}(P)):=\{A\in\Omega_1^{w,s}(B,{\rm Ad}(P))\,|\,\langle A,A\rangle^{\omega_0}_s<\infty\}$$
and ${\mathcal A}_P^{H^s}$ the completion of $\mathcal A_P^{\infty}$ with respect to 
$\langle\,\,,\,\,\rangle^{\omega_0}_s$, that is, 
$${\mathcal A}_P^{H^s}:=\{\omega_0+A\,|\,A\in\Omega_1^{H^s}(B,{\rm Ad}(P))\}.$$
Let $\Omega_{\mathcal T,1}^{H^s}(P,\mathfrak g)$ be the completion of $\Omega_{\mathcal T,1}^{\infty}(P,\mathfrak g)$ corresponding to 
$\Omega_1^{H^s}(B,{\rm Ad}(P))$.  

\vspace{0.25truecm}

{\small
\centerline{
\unitlength 0.1in
\begin{picture}( 48.2000,  6.4000)( -9.0000,-15.6000)
\put(43.6000,-11.4000){\makebox(0,0)[rb]{$\Omega_1^{H^s}(M,{\rm Ad}(P))$}}%
\put(14.5000,-11.6000){\makebox(0,0)[lb]{$T_{\omega_0}\mathcal A_P^{H^s}=\Omega_{\mathcal T,1}^{H^s}(P,\mathfrak g)$}}%
\put(10.7000,-10.9000){\makebox(0,0)[lb]{$\approx$}}%
%
\special{pn 8}%
\special{ar 3860 1210 60 90  6.2831853 6.2831853}%
\special{ar 3860 1210 60 90  0.0000000 3.1415927}%
%
\special{pn 8}%
\special{pa 3860 1300}%
\special{pa 3860 1210}%
\special{fp}%
\put(38.0000,-13.9000){\makebox(0,0)[lt]{$A$}}%
\put(17.4000,-13.6000){\makebox(0,0)[lt]{$\widehat A(:=\omega-\omega_0)$}}%
%
\special{pn 8}%
\special{ar 1800 1210 60 90  6.2831853 6.2831853}%
\special{ar 1800 1210 60 90  0.0000000 3.1415927}%
%
\special{pn 8}%
\special{pa 1800 1300}%
\special{pa 1800 1210}%
\special{fp}%
\put(9.0000,-11.5000){\makebox(0,0)[rb]{$\mathcal A_P^{H^s}$}}%
\put(6.9000,-13.9000){\makebox(0,0)[lt]{$\omega$}}%
%
\special{pn 8}%
\special{ar 740 1210 60 90  6.2831853 6.2831853}%
\special{ar 740 1210 60 90  0.0000000 3.1415927}%
%
\special{pn 8}%
\special{pa 740 1300}%
\special{pa 740 1210}%
\special{fp}%
%
\special{pn 8}%
\special{pa 1292 1446}%
\special{pa 1040 1446}%
\special{fp}%
\special{sh 1}%
\special{pa 1040 1446}%
\special{pa 1108 1466}%
\special{pa 1094 1446}%
\special{pa 1108 1426}%
\special{pa 1040 1446}%
\special{fp}%
%
\special{pn 8}%
\special{pa 1300 1446}%
\special{pa 1550 1446}%
\special{fp}%
\special{sh 1}%
\special{pa 1550 1446}%
\special{pa 1484 1426}%
\special{pa 1498 1446}%
\special{pa 1484 1466}%
\special{pa 1550 1446}%
\special{fp}%
%
\special{pn 8}%
\special{pa 3360 1450}%
\special{pa 2784 1450}%
\special{fp}%
\special{sh 1}%
\special{pa 2784 1450}%
\special{pa 2852 1470}%
\special{pa 2838 1450}%
\special{pa 2852 1430}%
\special{pa 2784 1450}%
\special{fp}%
%
\special{pn 8}%
\special{pa 3340 1450}%
\special{pa 3530 1450}%
\special{fp}%
\special{sh 1}%
\special{pa 3530 1450}%
\special{pa 3464 1430}%
\special{pa 3478 1450}%
\special{pa 3464 1470}%
\special{pa 3530 1450}%
\special{fp}%
\put(29.9000,-10.9000){\makebox(0,0)[lb]{$\approx$}}%
\put(15.8000,-15.6000){\makebox(0,0)[lt]{($\omega_0:$ The base point of $\mathcal A_P^{H^s}$)}}%
\end{picture}%
\hspace{4truecm}}
}

\vspace{0.85truecm}

\noindent
Let $\mathcal G_P^{\infty}$ be the group of all $C^{\infty}$-gauge transformations ${\bf g}$'s of $P$ with $\pi\circ{\bf g}=\pi$.  
For each ${\bf g}\in\mathcal G_P^{\infty}$, 
$\widehat{\bf g}\in C^{\infty}(P,G)$ is defined by ${\bf g}(u)=u\widehat{\bf g}(u)\,\,\,(u\in P)$.  
This element $\widehat{\bf g}$ satisfies 
$$\widehat{\bf g}(ug)={\rm Ad}(g^{-1})(\widehat{\bf g}(u))\,\,\,\,(\forall u\in P,\,\,\forall g\in G),$$
where ${\rm Ad}$ denotes the homomorphism of $G$ to ${\rm Aut}(G)$ defined by 
${\rm Ad}(g_1)(g_2):=g_1\cdot g_2\cdot g_1^{-1}$ ($g_1,g_2\in G$).  
Under the correspondence ${\bf g}\leftrightarrow\widehat{\bf g}$, $\mathcal G_P^{\infty}$ is identified with 
$$\widehat{\mathcal G}_P^{\infty}:=\{\widehat{\bf g}\in C^{\infty}(P,G)\,|\,
\widehat{\bf g}(ug)={\rm Ad}(g^{-1})(\widehat{\bf g}(u))\,\,\,(\forall\,u\in P,\,\,\forall\,g\in G)\}.$$
For $\widehat{\bf g}\in\widehat{\mathcal G}_P^{\infty}$, 
the $C^{\infty}$-section $\breve{\bf g}$ of the associated $G$-bundle $P\times_{{\rm Ad}}G$ is 
defined by 
$\breve{\bf g}(x):=u\cdot\widehat{\bf g}(u)\,\,\,(x\in M)$, where $u$ is any element of $\pi^{-1}(x)$.  
Under the correspondence $\widehat{\bf g}\leftrightarrow\breve{\bf g}$, 
$\widehat{\mathcal G}_P^{\infty}(=\mathcal G_P^{\infty})$ is identified with the space 
$\Gamma^{\infty}(P\times_{{\rm Ad}}G)$ of all $C^{\infty}$-sections of $P\times_{{\rm Ad}}G$.  
The $H^{s+1}$-completion of $\Gamma^{\infty}(P\times_{{\rm Ad}}G)$ was defined by Groisser and Parker 
(see Section 1 (P668) of \cite{GP1}).  
Denote by $\Gamma^{H^{s+1}}(P\times_{{\rm Ad}}G)$ this completion.  
Also, denote by $\mathcal G_P^{H^{s+1}}$ (resp. $\widehat{\mathcal G}_P^{H^{s+1}}$) the completion of 
$\mathcal G_P^{\infty}$ (resp. $\widehat{\mathcal G}_P^{\infty}$) corresponding to 
$\Gamma^{H^{s+1}}(P\times_{{\rm Ad}}G)$.  
Assume that $s>\frac{1}{2}\,{\rm dim}\,M-1$.  Then, according to Lemma 1.2 of \cite{U}, 
this $H^{s+1}$-completion $\mathcal G_P^{H^{s+1}}$ is a $C^{\infty}$-Hilbert Lie group and 
the gauge action $\mathcal G_P^{H^{s+1}}\curvearrowright\mathcal A_P^{H^s}$ is of class $C^{\infty}$.  
However, this action does not act isometrically on the Hilbert space 
$(\mathcal A_P^{H^s},\langle\,\,,\,\,\rangle^{\omega_0}_s)$.  
Define a Riemannian metric ${\it g}_s$ on $\mathcal A_P^{H^s}$ by 
$({\it g}_s)_{\omega}:=\langle\,\,\,\,\,\rangle^{\omega}_s$ ($\omega\in\mathcal A_P^{H^s}$).  
This Riemannian metric ${\it g}_s$ is non-flat and translation-invariant.  
The gauge action $\mathcal G_P^{H^{s+1}}\curvearrowright\mathcal A_P^{H^s}$ acts isometrically on 
the Riemannian Hilbert manifold $(\mathcal A_P^{H^s},{\it g}_s)$.  
Note that the Hilbert space $(\mathcal A_P^{H^s},\langle\,\,,\,\,\rangle^{\omega_0}_s)$ is regarded as 
the tangent space of $(\mathcal A_P^{H^s},{\it g}_s)$ at $\omega_0$.  
Give the moduli space $\mathcal M_P^{H^s}:=\mathcal A_P^{H^s}/\mathcal G_P^{s+1}$ the Riemannian orbimetric 
$\overline{\it g}_s$ such that the orbit map 
$\pi_{\mathcal M_P}:(\mathcal A_P^{H^s},{\it g}_s)\to(\mathcal M_P^{H^s},\overline{\it g}_s)$ is 
a Riemannian orbisubmersion.  

\vspace{0.25truecm}

{\small
\centerline{
\unitlength 0.1in
\begin{picture}( 23.6000,  4.8500)( 10.3000,-16.1500)
\put(10.3000,-13.2000){\makebox(0,0)[lb]{$\mathcal G^{H^{s+1}}_P$}}%
\put(15.4000,-13.1000){\makebox(0,0)[lb]{$\approx$}}%
\put(20.2000,-13.3000){\makebox(0,0)[lb]{$\widehat{\mathcal G}^{H^{s+1}}_P$}}%
\put(25.5000,-13.0000){\makebox(0,0)[lb]{$\approx$}}%
\put(29.7000,-13.4000){\makebox(0,0)[lb]{$\Gamma^{H^{s+1}}(P\times_{\rm Ad}G)$}}%
\put(11.1000,-15.3000){\makebox(0,0)[lt]{${\bf g}$}}%
\put(20.7000,-15.2000){\makebox(0,0)[lt]{$\widehat{\bf g}$}}%
\put(32.9000,-15.2000){\makebox(0,0)[lt]{$\breve{\bf g}$}}%
%
\special{pn 8}%
\special{ar 1150 1360 60 90  6.2831853 6.2831853}%
\special{ar 1150 1360 60 90  0.0000000 3.1415927}%
%
\special{pn 8}%
\special{pa 1150 1450}%
\special{pa 1150 1360}%
\special{fp}%
%
\special{pn 8}%
\special{ar 2110 1350 60 90  6.2831853 6.2831853}%
\special{ar 2110 1350 60 90  0.0000000 3.1415927}%
%
\special{pn 8}%
\special{pa 2110 1440}%
\special{pa 2110 1350}%
\special{fp}%
%
\special{pn 8}%
\special{ar 3330 1360 60 90  6.2831853 6.2831853}%
\special{ar 3330 1360 60 90  0.0000000 3.1415927}%
%
\special{pn 8}%
\special{pa 3330 1450}%
\special{pa 3330 1360}%
\special{fp}%
%
\special{pn 8}%
\special{pa 1692 1606}%
\special{pa 1430 1606}%
\special{fp}%
\special{sh 1}%
\special{pa 1430 1606}%
\special{pa 1498 1626}%
\special{pa 1484 1606}%
\special{pa 1498 1586}%
\special{pa 1430 1606}%
\special{fp}%
%
\special{pn 8}%
\special{pa 1700 1606}%
\special{pa 1920 1606}%
\special{fp}%
\special{sh 1}%
\special{pa 1920 1606}%
\special{pa 1854 1586}%
\special{pa 1868 1606}%
\special{pa 1854 1626}%
\special{pa 1920 1606}%
\special{fp}%
%
\special{pn 8}%
\special{pa 2810 1610}%
\special{pa 2476 1610}%
\special{fp}%
\special{sh 1}%
\special{pa 2476 1610}%
\special{pa 2542 1630}%
\special{pa 2528 1610}%
\special{pa 2542 1590}%
\special{pa 2476 1610}%
\special{fp}%
%
\special{pn 8}%
\special{pa 2820 1610}%
\special{pa 3030 1610}%
\special{fp}%
\special{sh 1}%
\special{pa 3030 1610}%
\special{pa 2964 1590}%
\special{pa 2978 1610}%
\special{pa 2964 1630}%
\special{pa 3030 1610}%
\special{fp}%
\end{picture}%
\hspace{0.75truecm}}
}

\vspace{0.5truecm}

Let $\displaystyle{S^1:=\{e^{2\pi\sqrt{-1}t}\,|\,t\in[0,1]\}}$ and define $z:[0,1]\to S^1$ by $z(t):=e^{2\pi\sqrt{-1}t}$ ($t\in[0,1]$).  
Fix $x_0\in B$ and $u_0\in\pi^{-1}(x_0)$.  
Take $c:S^1\to B$ be a $C^{\infty}$-loop with $c(1)=x_0$.  
Denote by $\pi^c:c^{\ast}P\to S^1$ the induced bundle of $P$ by $c$, which is identified with the trivial $G$-bundle 
$P_o:=S^1\times G$ over $S^1$ by $\sigma$.  
Define an immersion $\iota_c$ of the induced bundle $c^{\ast}P$ into $P$ by 
$\iota_c(z(t),u)=u\,\,\,((z(t),u)\in c^{\ast}P)$.  

\vspace{0.25truecm}

\noindent
{\bf Definition 9.1.}\ 
Define a map ${\rm hol}_c:{\mathcal A}_P^{H^s}\to G$ by 
$$(\mathcal P_{c\circ z}^{\omega}\circ(P_{c\circ z}^{\omega_0})^{-1})(u_0)=u_0\cdot{\rm hol}_c(\omega),$$
where 
$\mathcal P_{c\circ z}^{\omega}$ (resp. $\mathcal P_{c\circ z}^{\omega_0}$) denotes the parallel translation along $c\circ z$ 
with respect to $\omega$ (resp. $\omega_0$).  We call this map ${\rm hol}_c$ the {\it holonomy map along} $c$.  

\vspace{0.25truecm}

In particular, in the case where $P$ is the trivial $G$-bundle $\pi^o:S^1\times G\to S^1$, 
$\mathcal A_P^{H^s}$ is identified with the Hilbert space $H^s([0,1],\mathfrak g)$ of all $H^s$-curves in the Lie algebra $\mathfrak g$ of $G$ 
and the holonomy map ${\rm hol}_c$ along $c(t)=t$ ($t\in [0,1]$) coincides with the parallel 
transport map $\phi:H^s([0,1],\mathfrak g)\to G$ for $G$ stated in Example 4.1, where $s$ may be any non-negative integer because $[0,1]$ is 
of one-dimension.  

The {\it based gauge group} $(\mathcal G_P^{H^{s+1}})_x$ {\it at} $x\in M$ is defined by 
$$(\mathcal G_P^{H^{s+1}})_x:=\{{\bf g}\in\mathcal G_P^{H^{s+1}}\,|\,
\widehat{\bf g}(\pi^{-1}(x)))=\{e\}\},$$
where $e$ denotes the identity element of $G$.  
For a finite subgroup $\Gamma$ of $G$ and $u\in P$, we define a closed subgroup 
$(\mathcal G_P^{H^{s+1}})_{u,\Gamma}$ of $\mathcal G_P^{H^{s+1}}$ by 
$$(\mathcal G_P^{H^{s+1}})_{u,\Gamma}:=\{{\bf g}\in\mathcal G_P^{H^{s+1}}\,|\,\widehat{\bf g}(u)\in\Gamma\},$$
where $\widehat{\bf g}$ is the element of $H^{s+1}(P,G)$ corresponding to ${\bf g}$.  

\vspace{0.5truecm}

\noindent
{\it Example 9.1.}\ 
Let $G$ be a compact semi-simple Lie group and consider the trivial $G$-bundle $P_o:=S^1\times G$ over $S^1$.  
Also let $s$ be any non-negative integer.  
Since ${\rm Ad}(P_o)$ is identified with the trivial $\mathfrak g$-bundle $P'_o:=S^1\times\mathfrak g$ over $S^1$, 
$\mathcal A_{P_o}^{H^s}(\approx\Omega_1^{H^s}([0,1],{\rm Ad}(P_o)))$ is identified with $H^s([0,1],\mathfrak g)$.  
Also, $\mathcal G_{P_o}^{H^{s+1}}$ is identified with the Hilbert Lie group $H^{s+1}([0,1],G)$ of all $H^{s+1}$-paths in $G$.  
Here we note that $H^{s+1}([0,1],G)$ is a $C^{\infty}$-Hilbert Lie group and $H^{s+1}([0,1],G)\curvearrowright H^s([0,1],\mathfrak g)$ is a $C^{\infty}$-action 
because $[0,1]$ is of one-dimension.  For closed subgroup $H$ of $G\times G$, we define a closed subgroup $P(G,H)$ of $H^{s+1}[0,1],G)$ by 
$$P(G,H):=\{{\bf g}\in H^{s+1}([0,1],G)\,|\,({\bf g}(0),{\bf g}(a))\in H\}.$$
In particular, denote by $\Lambda_e^{H^{s+1}}(G)$ the loop group $P(G,\{e\}\times\{e\})$ at the idenitity element $e$ of $G$.  
Let $\Gamma$ be a finite subgroup of $G$ and $K$ a closed subgroup of $G$ such that 
$(G,K)$ is a reductive pair, that is, there exists a subspace $\mathfrak p$ of $\mathfrak g$ satisfying 
$\mathfrak g=\mathfrak k\oplus\mathfrak p$ and $[\mathfrak k,\mathfrak p]\subset\mathfrak p$, where 
$\mathfrak g$ and $\mathfrak k$ denote the Lie algebras of $G$ and $K$, respectively. 
Let $\mathcal B:=\mathcal A_{P_o}^{H^0}$ and $\mathcal G_1:=P(G,\Gamma\times K)$.  
Then the pair $(\mathcal B,\mathcal G_1)$ satisfies the conditions (I) and (II) stated in Introduction, where we note that 
the moduli space $\mathcal B/\mathcal G_1$ is orbi-diffeomorphic to the orbifold $\Gamma\setminus G\,/\,K$.  

\vspace{0.5truecm}

\noindent
{\bf Assumption.}\ \ Assume that $\displaystyle{s\geq\frac{1}{2}\,{\rm dim}\,B\,-\,1}$.  

\vspace{0.5truecm}

Denote by $C^{\infty}([0,1],M)$ the space of all $C^{\infty}$-paths in $M$.  
Take $c\in C^{\infty}([0,1],M)$.  
Denote by $\pi^c:c^{\ast}P\to[0,1]$ the induced bundle of $P$ by $c$, which is isomophic to the trivial $G$-bundle 
$P_o(=[0,1]\times G)\to [0,1]$.  
Define an immersion $\iota_c$ of the induced bundle $c^{\ast}P$ into $P$ by 
$\iota_c(t,u)=u\,\,\,((t,u)\in c^{\ast}P)$.  
Fix a base point $\omega_0$ of $\mathcal A_P^{\infty}$.  
We shall define a map linking $\mathcal A_P^{H^s}$ to $H^0([0,1],\mathfrak g)$.  

Let $x_0,\,\,c,\,\,u_0$ be as above and $\sigma$ the horizontal lift of $c\circ z$ starting from $u_0$ with respect to $\omega_0$.  

Then the pull-back bundle $c^{\ast}P$ of $P$ by $c$ is identified with the trivial $G$-bundle $S^1\times G$ over $S^1$ by $\sigma$ as follows.  
Define a map $\eta:S^1\times G\to c^{\ast}P$ by 
$$\eta(z(t),g):=(z(t),\sigma(t)g)\quad\,\,((t,g)\in[0,1]\times G).$$
It is clear that $\eta$ is a bundle isomorphism.  
Throughout this bundle isomorphism $\eta$, $c^{\ast}P$ is identified with $S^1\times G$.  
Similarly, a bundle isomorphism $\bar{\eta}:[0,1]\times G\to (c\circ z)^{\ast}P$ is defined by 
$$\bar{\eta}(t,g):=(t,\sigma(t)g)\quad\,\,((t,g)\in[0,1]\times G).$$
Throughout this bundle isomorphism $\bar{\eta}$, $(c\circ z)^{\ast}P$ is identified with the trivial bundle $[0,1]\times G$.  
The natural embedding $\iota_{\sigma}:[0,1]\times G\hookrightarrow P$ by 
$$\iota_{\sigma}(t,g):=\sigma(t)g\quad\,\,((t,g)\in[0,1]\times G).$$
For $\omega\in\mathcal A_P^{H^s}$, the pull-back connection $\iota_{\sigma}^{\ast}\omega$ of $(c\circ z)^{\ast}P$ is defined.  
Then we have 
$$(\iota_{\sigma}^{\ast}\omega)_{(t_0,g)}\left(\left(\frac{\partial}{\partial t}\right)_{(t_0,g)}\right)
:=\omega_{\sigma(t_0)}(\sigma'(t_0))\quad\,\,((t_0,g)\in[0,1]\times G)).$$
By the one-to-one correspondence 
$$\iota_{\sigma}^{\ast}\omega\quad\longleftrightarrow\quad t\mapsto\omega_{\sigma(t)}(\sigma'(t)) \,\,\,(t\in[0,1]),$$
$\mathcal A_{(c\circ z)^{\ast}P}^{H^s}$ is identified with the Hilbert space $H^s([0,1],\mathfrak g)$ of all $H^s$-paths in $\mathfrak g$.  

\vspace{0.25truecm}

\noindent
{\bf Definition 9.2.} Define a map $\mu_c:\mathcal A_P^{H^s}\to H^s([0,1],\mathfrak g)$ by 
$$(\mu_c(\omega))(t):=-\omega_{\sigma(t)}(\sigma'(t))\quad\,\,(t\in[0,1],\,\,\,\,\omega\in\mathcal A_P^{H^s}).\leqno{(2.3)}$$
As above, $\omega_{\sigma(z(t_0))}(\sigma'(t_0))$ is equal to 
$\displaystyle{(\sigma^{\ast}\omega)_{(t_0,g)}\left(\left(\frac{\partial}{\partial t}\right)_{(t_0,g)}\right)}$.  
From this fact, we call this map $\mu_c$ the {\it pull-back connection map by $c$}.  
By the defiinitions of ${\rm hol}_c$ and $\mu_c$, we can show the following relation.  

\vspace{0.25truecm}

\noindent
{\bf Lemma 9.1.} {\sl Among ${\rm hol}_c,\,\phi$ and $\mu_c$, the relation ${\rm hol}_c=\phi\circ\mu_c$ holds.}

\vspace{0.25truecm}

\noindent
{\bf Definition 9.3.} Define a map $\lambda_c:\mathcal G_P^{H^{s+1}}\to H^{s+1}([0,1],G)$ by 
$$\lambda_c({\bf g})(t):=\widehat{\bf g}(\sigma(t))\quad\,\,(t\in[0,1],\,\,\,\,{\bf g}\in\mathcal G_P^{H^{s+1}}).$$

\vspace{0.25truecm}

Take a $C^{\infty}$-loop $c:[0,1]\to M$ and $\widetilde c$ the horizontal lift of $c$ with respect to $\omega_0$.  
Easily we can show the following fact from the definitions of $\mu_c$ and $\lambda_c$.  

\vspace{0.5truecm}

\noindent
{\bf Proposition 9.2.} {\sl {\rm (i)} 
Between $\mu_c$ and $\lambda_c$, the following relation holds:
$$\mu_c({\bf g}\cdot\omega)=\lambda_c({\bf g})\cdot\mu_c(\omega)\quad\,\,
({\bf g}\in\mathcal G_P^{H^{s+1}},\,\,\omega\in\mathcal A_P^{H^s}).$$

{\rm(ii)} 
$\mu_c$ maps $(\mathcal G_P^{H^{s+1}})_{\widetilde c(0),\Gamma}$-orbits in 
$\mathcal A_P^{H^s}$ to $P(G,\Gamma^{\omega_0})$-orbits in $H^0([0,1],\mathfrak g)$ and hence 
there uniquely exists the map $\overline{\mu}_c^{\Gamma}$ between the orbit spaces \newline
$\mathcal A_P^{H^s}/(\mathcal G_P^{H^{s+1}})_{\widetilde c(0),\Gamma}$ and 
$$H^0([0,1],\mathfrak g)/P(G,\Gamma^{\omega_0})(=G/(R_{{\rm hol}_c(\omega_0)}\circ{\rm Ad}(\Gamma)\circ 
R_{{\rm hol}_c(\omega_0)}^{-1}))$$
satisfying 
$$\overline{\mu}_c^{\Gamma}\circ\pi_{(\mathcal G_P^{H^{s+1}})_{\widetilde c(0),\Gamma}}
=\pi_{\Gamma^{\omega_0}}\circ\phi\circ\mu_c=\pi_{\Gamma}\circ{\rm hol}_c,$$
where $\pi_{\Gamma}$ (resp. $\pi_{(\mathcal G_P^{H^{s+1}})_{\widetilde c(0),\Gamma}}$) 
denotes the orbit map of the action 
${\rm Ad}(\Gamma)\curvearrowright G$ (resp. $(\mathcal G_P^{H^{s+1}})_{\widetilde c(0),\Gamma}\curvearrowright\mathcal A_P^{H^s}$).
}

\vspace{0.5truecm}

Thus the study of $P(G,\Gamma)$-orbits in $H^1([0,1],\mathfrak g)$ 
leads to that of $(\mathcal G_P^{H^{s+1}})_{\widetilde c(0),\Gamma}$-orbits ($u\in P$) in 
$\mathcal A_P^{H^s}$ through $\mu_c$ (see Figure 9.1).  

\vspace{0.5truecm}

{\small 
\centerline{
\unitlength 0.1in
\begin{picture}( 79.3000, 15.0500)(-43.0000,-21.2500)
\put(18.9000,-13.5000){\makebox(0,0)[rt]{$\mathcal A_P^{H^s}$}}%
\put(31.0000,-13.5000){\makebox(0,0)[lt]{$H^0([0,1],\mathfrak g)$}}%
\put(12.4000,-9.1000){\makebox(0,0)[rb]{$(\mathcal G_P^{H^{s+1}})_{\widetilde c(0),\Gamma}$}}%
\put(30.4000,-9.1000){\makebox(0,0)[rb]{$P(G,\Gamma^{\omega_0})$}}%
%
\special{pn 8}%
\special{pa 1440 850}%
\special{pa 2330 850}%
\special{fp}%
\special{sh 1}%
\special{pa 2330 850}%
\special{pa 2264 830}%
\special{pa 2278 850}%
\special{pa 2264 870}%
\special{pa 2330 850}%
\special{fp}%
\put(17.8000,-7.9000){\makebox(0,0)[lb]{$\lambda_c$}}%
%
\special{pn 8}%
\special{pa 2050 1430}%
\special{pa 2880 1430}%
\special{fp}%
\special{sh 1}%
\special{pa 2880 1430}%
\special{pa 2814 1410}%
\special{pa 2828 1430}%
\special{pa 2814 1450}%
\special{pa 2880 1430}%
\special{fp}%
\put(23.7000,-13.8000){\makebox(0,0)[lb]{$\mu_c$}}%
%
\special{pn 13}%
\special{ar 1100 1310 590 450  5.0767644 5.9136541}%
%
\special{pn 13}%
\special{pa 1660 1150}%
\special{pa 1700 1230}%
\special{fp}%
\special{sh 1}%
\special{pa 1700 1230}%
\special{pa 1688 1162}%
\special{pa 1676 1182}%
\special{pa 1652 1180}%
\special{pa 1700 1230}%
\special{fp}%
%
\special{pn 13}%
\special{ar 2900 1320 590 450  5.0767644 5.9136541}%
%
\special{pn 13}%
\special{pa 3460 1160}%
\special{pa 3500 1240}%
\special{fp}%
\special{sh 1}%
\special{pa 3500 1240}%
\special{pa 3488 1172}%
\special{pa 3476 1192}%
\special{pa 3452 1190}%
\special{pa 3500 1240}%
\special{fp}%
%
\special{pn 8}%
\special{pa 1710 1610}%
\special{pa 1710 1980}%
\special{fp}%
\special{sh 1}%
\special{pa 1710 1980}%
\special{pa 1730 1914}%
\special{pa 1710 1928}%
\special{pa 1690 1914}%
\special{pa 1710 1980}%
\special{fp}%
%
\special{pn 8}%
\special{pa 3550 1580}%
\special{pa 3550 1950}%
\special{fp}%
\special{sh 1}%
\special{pa 3550 1950}%
\special{pa 3570 1884}%
\special{pa 3550 1898}%
\special{pa 3530 1884}%
\special{pa 3550 1950}%
\special{fp}%
\put(20.0000,-20.3000){\makebox(0,0)[rt]{$\mathcal A_P^{H^s}/(\mathcal G_P^{H^{s+1}})_{\widetilde c(0),\Gamma}$}}%
\put(29.6000,-20.4000){\makebox(0,0)[lt]{$H^0([0,1],\mathfrak g)/P(G,\Gamma)$}}%
%
\special{pn 8}%
\special{pa 2110 2120}%
\special{pa 2800 2120}%
\special{fp}%
\special{sh 1}%
\special{pa 2800 2120}%
\special{pa 2734 2100}%
\special{pa 2748 2120}%
\special{pa 2734 2140}%
\special{pa 2800 2120}%
\special{fp}%
\put(23.7000,-20.8000){\makebox(0,0)[lb]{$\overline{\mu}_c^{\Gamma}$}}%
%
\special{pn 8}%
\special{pa 2010 1640}%
\special{pa 2850 1990}%
\special{fp}%
\special{sh 1}%
\special{pa 2850 1990}%
\special{pa 2796 1946}%
\special{pa 2802 1970}%
\special{pa 2782 1984}%
\special{pa 2850 1990}%
\special{fp}%
\put(23.4000,-17.5000){\makebox(0,0)[lb]{{\small$\pi_{\Gamma}\circ{\rm hol}_c$}}}%
\put(36.3000,-17.0000){\makebox(0,0)[lt]{$\pi_{\Gamma}\circ\phi$}}%
%
\special{pn 8}%
\special{ar 2002 1772 72 72  0.6528466 5.7088805}%
%
\special{pn 8}%
\special{pa 2062 1732}%
\special{pa 2082 1782}%
\special{fp}%
\special{sh 1}%
\special{pa 2082 1782}%
\special{pa 2076 1714}%
\special{pa 2062 1732}%
\special{pa 2040 1728}%
\special{pa 2082 1782}%
\special{fp}%
%
\special{pn 8}%
\special{ar 3210 1770 72 72  0.6528466 5.7088805}%
%
\special{pn 8}%
\special{pa 3280 1740}%
\special{pa 3300 1790}%
\special{fp}%
\special{sh 1}%
\special{pa 3300 1790}%
\special{pa 3294 1722}%
\special{pa 3280 1740}%
\special{pa 3258 1736}%
\special{pa 3300 1790}%
\special{fp}%
\put(16.1000,-18.4000){\makebox(0,0)[rb]{$\pi_{(\mathcal G_P^{H^{s+1}})_{\widetilde c(0),\Gamma}}$}}%
\end{picture}%
\hspace{14truecm}}
}

\vspace{0.5truecm}

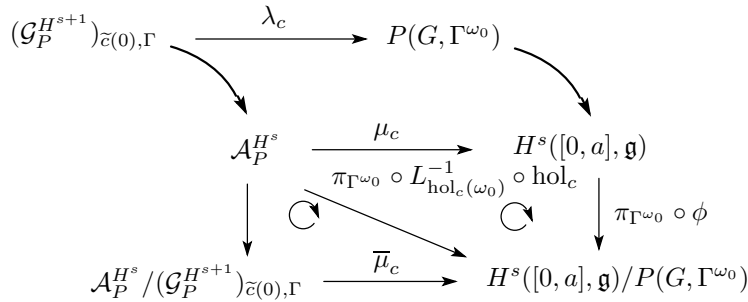
\centerline{{\bf Figure 9.1$\,\,:\,\,$ $P(G,\Gamma)$-orbits and 
$(\mathcal G_P^{H^{s+1}})_{\widetilde c(0),\Gamma}$-orbits}}

\vspace{0.5truecm}

\noindent
{\it Remark 9.2.}\ Since $\mu_c$ is a bounded linear map of $(\mathcal A_P^{H^s},\langle\,\,,\,\,\rangle^{\omega_0}_s)$ (regarded as a Hilbert space) onto 
$H^0([0,1],\mathfrak g)$, the following facts hold:

(i)\ The inverse images $\mu_c^{-1}(u)$'s ($u\in H^0([0,1],\mathfrak g)$) are closed affine subspaces of the affine space $\mathcal A_P^{H^s}$ 
and they are parallel.  Hence, for any a $C^{\infty}$-orbisubmanifold $S$ in $H^0([0,1],\mathfrak g)/P(G,\Gamma)$, 
$(\pi_{\Gamma}\circ{\rm hol}_c)^{-1}(S)$ is a cylindrical $C^{\infty}$-orbisubmanifold of finite codimension in 
the Hilbert space $(\mathcal A_P^{H^s},\langle\,\,,\,\,\rangle^{\omega_0}_s)$;

(ii)\ The map $\mu_c$ is a $C^{\infty}$-submersion as a map between the Riemannian Hilbert manifolds, where we regard 
$(\mathcal A_P^{H^s},\langle\,\,,\,\,\rangle^{\omega_0}_s)$ and $H^0([0,1],\mathfrak g)$ as Riemannian Hilbert manifolds.  

Also, it is easy to show that the following facts holds:

(iii)\ The operator norm $\|(d\mu_c)_{\omega}\|_{\rm op}$ of the differential $(d\mu_c)_{\omega}$ of $\mu_c$ at any point $\omega$ is smaller than 
or equal to one.  

\vspace{0.5truecm}

Denote by $\pi^{H^s}_{\mathcal M_P}$ the orbit map of the action $\mathcal G_P^{H^{s+1}}\curvearrowright\mathcal A_P^{H^s}$.  
Also, let $\exp_{\omega_0}$ be the exponential map of the Riemannian Hilbert manifold $(\mathcal A_P^{H^s},{\it g}_s)$ at $\omega_0$.  
Then, for an orbisubmanifold $S$ in $H^0([0,1],\mathfrak g)/P(G,\Gamma^{\omega_0})$, we can costruct an orbisubmanifolds 
$\pi_{\mathcal M_P}^{H^s}(\exp_{\omega}((\pi_{\Gamma}\circ{\rm hol}_c)^{-1}(S)))$ in the muduli space 
$(\mathcal M_P^{H^s},\overline{\it g}_s)$ (see Figure 9.2).  

Let $\mathcal{YM}_P^{H^s}$ be the Hilbert space of all Yang-Mills $H^s$-connections of the $G$-bundle $\pi:P\to M$.  
Also, in the case of ${\rm dim}\,M=4$, denote by $\mathcal{SD}_P^{H^2}$ the Hilbert space of all self-dual 
$H^2$-connections of the $G$-bundle $\pi:P\to M$, where we note that $s\geq 2$ in this case.  
Note that a self-dual connection means ``instanton'' because $M$ is compact.  
The Yang-Mills moduli space $\mathcal M_P^{\mathcal{YM}}:=\mathcal{YM}_P^{H^s}/\mathcal G_P^{H^{s+1}}$ 
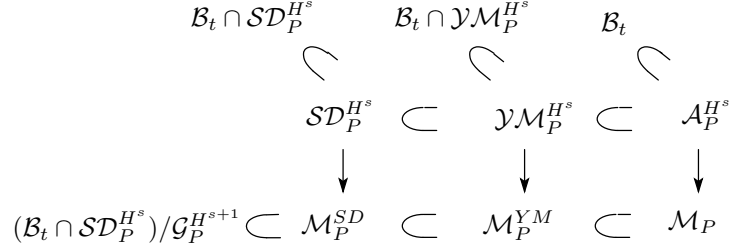
and the self-dual moduli space $\mathcal M_P^{\mathcal{SD}}:=\mathcal{SD}_P^{H^s}/\mathcal G_P^{H^{s+1}}$ are 
finite dimensional manifolds with singularity in general.  
Give these moduli spaces the singular Riemannian metrics induced naturally from the (non-flat) Riemannian metric 
${\it g}_s$ on $\mathcal A_P^{H^s}$.  Denote by $\overline{\it g}_s$ these singular Riemannian metrics.  

\vspace{0.5truecm}

{\small
\centerline{
\unitlength 0.1in
\begin{picture}( 97.6900, 23.2400)(-11.6900,-33.9700)
\put(65.3200,-20.9700){\makebox(0,0)[lb]{$(\mathcal A_P^{H^s},{\it g}_s)$}}%
\put(65.0800,-25.7700){\makebox(0,0)[lt]{$(\mathcal M_P^{H^s},\overline{\it g}_s)$}}%
%
\special{pn 8}%
\special{pa 6780 2202}%
\special{pa 6780 2550}%
\special{fp}%
\special{sh 1}%
\special{pa 6780 2550}%
\special{pa 6800 2482}%
\special{pa 6780 2496}%
\special{pa 6760 2482}%
\special{pa 6780 2550}%
\special{fp}%
%
\special{pn 8}%
\special{pa 7840 1768}%
\special{pa 8218 2106}%
\special{fp}%
\special{sh 1}%
\special{pa 8218 2106}%
\special{pa 8182 2048}%
\special{pa 8178 2070}%
\special{pa 8154 2076}%
\special{pa 8218 2106}%
\special{fp}%
\put(81.5800,-28.2600){\makebox(0,0)[lt]{$H^0([0,1],\mathfrak g)/P(G,\Gamma^{\omega_0})$}}%
%
\special{pn 8}%
\special{pa 8530 2370}%
\special{pa 8530 2720}%
\special{fp}%
\special{sh 1}%
\special{pa 8530 2720}%
\special{pa 8550 2652}%
\special{pa 8530 2666}%
\special{pa 8510 2652}%
\special{pa 8530 2720}%
\special{fp}%
\put(83.1300,-22.9500){\makebox(0,0)[lb]{$H^0([0,1],\mathfrak g)$}}%
\put(80.9600,-19.3700){\makebox(0,0)[lb]{$\mu_c$}}%
\put(86.0000,-24.5600){\makebox(0,0)[lt]{$\pi_{\Gamma}\circ\phi$}}%
\put(81.4200,-23.4200){\makebox(0,0)[rt]{$\pi_{\Gamma}\circ{\rm hol}_c$}}%
\put(68.3400,-23.0400){\makebox(0,0)[lt]{$\pi_{\mathcal M_P}^{H^s}$}}%
\put(70.8200,-17.3000){\makebox(0,0)[rb]{$\exp_{\omega_0}$}}%
\put(80.3000,-31.8100){\makebox(0,0)[rt]{$S$}}%
\put(66.5200,-23.4900){\makebox(0,0)[rt]{{\small $\exp_{\omega_0}((\pi_{\Gamma}\circ{\rm hol}_c)^{-1}(S))$}}}%
\put(71.1100,-30.0900){\makebox(0,0)[rt]{{\small $\pi_{\mathcal M_P}^{H^s}(\exp_{\omega_0}((\pi_{\Gamma}\circ{\rm hol}_c)^{-1}(S)))$}}}%
\put(78.8200,-33.9700){\makebox(0,0)[lt]{$\,$}}%
%
\special{pn 8}%
\special{pa 6474 2182}%
\special{pa 6452 2206}%
\special{pa 6424 2220}%
\special{pa 6398 2208}%
\special{pa 6392 2176}%
\special{pa 6402 2146}%
\special{pa 6416 2122}%
\special{sp}%
%
\special{pn 8}%
\special{pa 6412 2138}%
\special{pa 6442 2098}%
\special{fp}%
%
\special{pn 8}%
\special{pa 6474 2182}%
\special{pa 6504 2142}%
\special{fp}%
%
\special{pn 8}%
\special{pa 6504 2862}%
\special{pa 6482 2886}%
\special{pa 6454 2898}%
\special{pa 6428 2884}%
\special{pa 6424 2854}%
\special{pa 6432 2824}%
\special{pa 6446 2800}%
\special{sp}%
%
\special{pn 8}%
\special{pa 6442 2816}%
\special{pa 6472 2776}%
\special{fp}%
%
\special{pn 8}%
\special{pa 6504 2862}%
\special{pa 6534 2820}%
\special{fp}%
%
\special{pn 8}%
\special{pa 8126 3126}%
\special{pa 8104 3150}%
\special{pa 8076 3162}%
\special{pa 8052 3150}%
\special{pa 8046 3118}%
\special{pa 8054 3090}%
\special{pa 8068 3064}%
\special{sp}%
%
\special{pn 8}%
\special{pa 8064 3080}%
\special{pa 8094 3040}%
\special{fp}%
%
\special{pn 8}%
\special{pa 8126 3126}%
\special{pa 8156 3084}%
\special{fp}%
\put(72.6000,-16.6400){\makebox(0,0)[lb]{$(\mathcal A_P^{H^s},\langle\,\,,\,\,\rangle_s^{\omega_0})(=(T_{\omega_0}\mathcal A_P^{H^s},({\it g}_s)_{\omega_0}))$}}%
\put(76.0300,-12.4300){\makebox(0,0)[rb]{{\small $(\pi_{\Gamma}\circ{\rm hol}_c)^{-1}(S)$}}}%
%
\special{pn 8}%
\special{pa 7326 1352}%
\special{pa 7304 1330}%
\special{pa 7274 1318}%
\special{pa 7250 1332}%
\special{pa 7246 1362}%
\special{pa 7254 1392}%
\special{pa 7268 1416}%
\special{sp}%
%
\special{pn 8}%
\special{pa 7264 1398}%
\special{pa 7294 1438}%
\special{fp}%
%
\special{pn 8}%
\special{pa 7326 1352}%
\special{pa 7356 1394}%
\special{fp}%
%
\special{pn 8}%
\special{pa 7214 1656}%
\special{pa 6942 1880}%
\special{fp}%
\special{sh 1}%
\special{pa 6942 1880}%
\special{pa 7006 1854}%
\special{pa 6984 1846}%
\special{pa 6982 1822}%
\special{pa 6942 1880}%
\special{fp}%
%
\special{pn 8}%
\special{pa 7840 1776}%
\special{pa 8452 2748}%
\special{fp}%
\special{sh 1}%
\special{pa 8452 2748}%
\special{pa 8434 2680}%
\special{pa 8424 2702}%
\special{pa 8400 2702}%
\special{pa 8452 2748}%
\special{fp}%
\end{picture}%
\hspace{21.5truecm}}
}

\vspace{0.15truecm}

\centerline{{\bf Figure 9.2$\,\,:\,\,$ Submanifolds in the moduli space defined by ${\rm hol}_c$}}

\vspace{0.5truecm}

\noindent
{\bf Strategy}\ {\sl (i)\ We plan to find a pair $(S,c)$ of an orbisubmanifold $S$ in the Riemannian orbifold 
$H^0([0,1],\mathfrak g)/P(G,\Gamma)$ and a $C^{\infty}$-loop $c$ in $M$ such that 
$$\begin{array}{c}
\pi^{H^s}_{\mathcal M_P}(\exp_{\omega}((\pi_{\Gamma}\circ{\rm hol}_c)^{-1}(S)))\cap\mathcal{YM}_P^{H^s}\\
({\rm resp.}\,\,\pi^{H^s}_{\mathcal M_P}(\exp_{\omega}((\pi_{\Gamma}\circ{\rm hol}_c)^{-1}(S)))\cap\mathcal{SD}_P^{H^s})
\end{array}$$
gives an interesting submanifold in 
$\mathcal{YM}_P^{H^s}$ (resp. $\mathcal{SD}_P^{H^s}$).  

(ii)\ We plan to find a pair $(\{S_t\}_{t\in[0,T)},c)$ of a mean curvature flow $\{S_t\}_{t\in[0,T)}$ in the Riemannian orbifold 
$H^0([0,1],\mathfrak g)/P(G,\Gamma)$ and a $C^{\infty}$-loop $c$ in $M$ such that 
$$\begin{array}{c}
\{\pi^{H^s}_{\mathcal M_P}(\exp_{\omega}((\pi_{\Gamma}\circ{\rm hol}_c)^{-1}(S_t)))\cap\mathcal{YM}_P^{H^s}\}_{t\in[0,T)}\\
({\rm resp.}\ \ \{\pi^{H^s}_{\mathcal M_P}(\exp_{\omega}((\pi_{\Gamma}\circ{\rm hol}_c)^{-1}(S_t)))\cap\mathcal{SD}_P^{H^s}\}_{t\in[0,T)})
\end{array}$$ 
gives an interesting flow in $\mathcal{YM}_P^{H^s}$ (resp. $\mathcal{SD}_P^{H^s}$) (for example, a good flow collapsing to a singular point of 
$\mathcal{YM}_P^{H^s}$ (resp. $\mathcal{SD}_P^{H^s}$)).  
We will use Theorem A to find a good flow collapsing to a singular point of $\mathcal{YM}_P^{H^s}$ or $\mathcal{SD}_P^{H^s}$.}

\vspace{0.5truecm}


\vspace{0.5truecm}

{\small 
\rightline{Department of Mathematics, Faculty of Science}
\rightline{Tokyo University of Science, 1-3 Kagurazaka}
\rightline{Shinjuku-ku, Tokyo 162-8601 Japan}
\rightline{(koike@rs.tus.ac.jp)}
}

\end{document}